\documentclass[12pt, a4paper]{article}
\usepackage[english]{babel}
\usepackage{amsthm,amsmath,amssymb,enumitem,bm,mathrsfs}
\usepackage{graphicx}
\usepackage{xcolor}
\usepackage{authblk}
\RequirePackage[numbers]{natbib}
\RequirePackage[colorlinks,citecolor=blue,urlcolor=blue]{hyperref}

\newtheorem{theorem}{Theorem}
\newtheorem{lemma}{Lemma}
\newtheorem{corollary}{Corollary}
\newtheorem{proposition}{Proposition}
\theoremstyle{definition}

\newtheorem{remark}{Remark}

\renewenvironment{proof}{\noindent{\bfseries Proof}:}{}

\renewcommand{\epsilon}{\varepsilon}

\begin{document}

\title{\Large\bfseries A semi-parametric Bernstein-von Mises theorem at rate $n^{-3/2}$ for frequency estimation}
\author{\bfseries Paulo Serra}
\affil{Department of Mathematics, Vrije Universiteit Amsterdam}
\date{\today}
\maketitle

\begin{abstract}
We consider Bayesian inference for the frequency of a periodic signal of unknown shape, observed in Gaussian white noise on the interval $[0,n]$.
The information about the frequency grows like $n^3$, so the frequency can be estimated at rate $n^{-3/2}$, while the shape is an infinite-dimensional nuisance parameter.
We endow the frequency and the shape with independent priors, the prior on the shape being a random series in the trigonometric basis.
We give conditions on the priors under which the posterior for the signal contracts at a rate that, for suitable choices, is the minimax rate for Sobolev classes of smoothness at least one, up to a logarithmic factor and adaptively to the smoothness of the shape.
If, moreover, the prior keeps the shape away from shapes with a shorter period, also the marginal posteriors for the frequency and the shape contract around the truth.
Under further conditions, which hold for a Gaussian random series prior conditioned to keep the shape away from shapes with a shorter period, for every true shape with smallest period 1 and Sobolev smoothness larger than $3/2$, the marginal posterior for the frequency satisfies a semi-parametric Bernstein-von Mises theorem at rate $n^{-3/2}$, which we prove by moving the origin of time, a change that leaves the posterior unchanged when the prior on the shape is invariant under time shifts.
\end{abstract}

\textbf{Keywords:}
adaptive estimation,
Bernstein-von Mises theorem,
frequency estimation,
Gaussian white noise model,
local asymptotic normality,
period estimation,
posterior contraction,
random series prior,
semi-parametric Bayes.

\medskip

\textbf{Subject Classifications:} Primary 62G20; secondary 62F15.

\section{Introduction}\label{sec:introduction}

Stochastic processes with periodic features occur in many fields, such as signal processing and communications, mechanics, astronomy, economics, the atmospheric sciences, medicine and biology.
Depending on the application, the periodicity comes from modulation, sampling and coding operations, from the rotation of gears or planets, from the pulsation of stars, or from seasonality.
Statistical tools for such processes are used in ultrasonic tissue characterisation~\cite{Fellingham:1984}, the analysis of surface EMG signals~\cite{Lopresti:1995} and of circadian and ultradian rhythms~\cite{HadjAmar:2020}, non-linear control systems~\cite{Kim:1999}, plasma physics~\cite{Iavors:1996}, the study of variable stars~\cite{Hall:2000,Hall:2006}, laser vibrometry~\cite{Lavielle:2005}, the vibration analysis of rotating machines~\cite{Antoni:2004}, and the modelling of arrivals at call centres~\cite{Belitser:2013a}; cf.~\cite{Gardner:2006,Serpedin:2005} for overviews.

A periodic signal can be parametrised in terms of a \emph{frequency}, the reciprocal of its smallest period, and a \emph{shape}, the signal rescaled to have frequency 1.
The two are usually estimated separately: the frequency is estimated semi-parametrically, treating the shape as a nuisance parameter, and the shape is subsequently estimated with a plug-in approach.
To ensure that the shape can be optimally reconstructed, strong conditions are then imposed on the estimator of the frequency; cf.~\cite{Castillo:2006} for the plug-in approach in the Gaussian white noise model, \cite{Helmers:2003,Helmers:2005} for similar approaches in Poisson models (cf.~\cite{Kutoyants:1984} for kernel estimation of the intensity when the period is known), and~\cite{Klaassen:2005} for plug-in estimators of a nuisance parameter in general semi-parametric models.
It is therefore of interest to study methods where frequency and shape are estimated jointly, with guarantees for both; for least squares estimators in non-parametric regression with random design, this is done in~\cite{Hall:2000}.
In this paper we consider a Bayesian approach, where frequency and shape are endowed with independent priors, in the Gaussian white noise model.
Our main contribution is to show that the respective posterior recovers the signal at the minimax rate up to a logarithmic factor, that it localises the frequency, and that the marginal posterior for the frequency satisfies a semi-parametric Bernstein-von Mises theorem at rate $n^{-3/2}$, where $n$ is the length of the observation window.

The most commonly used tool for estimating a frequency is the \emph{periodogram}, introduced by Schuster~\cite{Schuster:1898}, which is the squared modulus of the Fourier transform of the observations, seen as a function of the frequency.
For parametric signals, such as a sinusoid in Gaussian noise, its maximiser is asymptotically equivalent to the least squares estimator~\cite{Hannan:1973}; cf.~\cite[Section~3.2]{Quinn:2001}.
Its analogue for point processes is the Bartlett periodogram~\cite{Bartlett:1963}, which gives consistent estimators of the frequency in parametric models of the intensity, efficient only in special cases~\cite{VereJones:1982}.

For models with a non-parametric shape, the frequency is also usually estimated by maximising criteria of periodogram type.
In the Gaussian white noise model one observes $X^{(n)}=(Y_t: t\in[0,n])$ with
\begin{equation}\label{eq:gwn_intro}
dY_t = f(\tau t)\,dt + dW_t, \quad t\in[0,n],
\end{equation}
for a standard Brownian motion $W$, a 1-periodic shape $f$ and a frequency $\tau>0$.
We write $\theta=f(\tau\,\cdot)$ for the signal, $P_\theta^{(n)}$ for the law of $X^{(n)}$, and $f_0$, $\tau_0$ and $\theta_0=f_0(\tau_0\,\cdot)$ for the true shape, frequency and signal.
In this model, such criteria give asymptotically efficient estimators of the frequency~\cite{Golubev:1988}.
A second-order efficient estimator, based on a penalised likelihood derived from a Gaussian prior on the Fourier coefficients of the shape, is given in~\cite{Castillo:2007}; cf.~\cite{Dalalyan:2006} for the analogous problem of estimating a shift.
Both works assume we observe the signal on a window centred at the origin, which removes the loss of information discussed below, and parametrise it by the period.
Similar results exist in non-parametric regression~\cite{Lavielle:2005}, also for several periods~\cite{Gassiat:2006,Lavielle:2005}, and with irregularly spaced observations~\cite{Hall:2000,Hall:2006,LevyLeduc:2008}, a trend~\cite{Vogt:2014} or an evolving periodic function~\cite{Genton:2007}.
For Poisson processes with a periodic, non-parametric intensity function, the frequency can be estimated by maximising the Bartlett periodogram~\cite{VereJones:1982}, by estimators designed for this model~\cite{Helmers:2003b,Bebbington:2004}, or by M-estimators~\cite{Belitser:2013a}.
Diffusions with a periodic drift~\cite{Hopfner:2011,Hopfner:2012} and Cox's proportional hazards model with a periodic baseline hazard and a known period~\cite{Pons:1988} have also been considered.

The periodogram has some limitations.
In non-parametric regression with observations on an irregular time grid, a common situation in astronomy, the least squares estimators of~\cite{Hall:2000} have a smaller asymptotic variance than the periodogram~\cite{Hall:2006}.
The periodogram tends to have many local maxima, in particular near multiples and fractions of the true frequency, which complicates both its optimisation and the quantification of the uncertainty of the estimates~\cite{Kim:2008,Wang:2012}.
In models with a non-parametric shape, maximising criteria of periodogram type does not distinguish the period from its multiples~\cite{Golubev:1988,Gassiat:2006}.
The estimators of~\cite{Golubev:1988,Castillo:2007} avoid this through constraints that use some knowledge of the true shape, such as a known separation from shapes with a shorter period in~\cite{Golubev:1988}, and, in~\cite[Section~2]{Castillo:2007}, a known range of periods, within which the criterion is maximised near the smallest period at which it is close to its maximum.

In the models above the frequency can be estimated at a rate much faster than the usual $n^{-1/2}$ one might expect for a real-valued parameter; note that $n$ is proportional to the number of observed periods.
In model~\eqref{eq:gwn_intro}, for instance, the Fisher information about $\tau$ when the shape $f$ is known, not constant and locally absolutely continuous on $\mathbb{R}$ with $f'\in L^2\big([0,1]\big)$, is
\begin{equation}\label{eq:fisher_intro}
\int_0^n t^2\big\{f'(\tau t)\big\}^2\,dt = \frac{n^3}{3}\|f'\|_2^2\{1+o(1)\},
\qquad
\|f'\|_2^2=\int_0^1\big\{f'(u)\big\}^2\,du,
\end{equation}
cf.~\cite[p.~209]{Ibragimov:1981}, so that the frequency can be estimated at rate $n^{-3/2}$.
When the shape is unknown, a change of frequency is partly confounded with a change of phase: with the origin of time at the start of the observation window, the score $tf'(\tau t)$ for $\tau$ is far from orthogonal to the score $f'(\tau t)$ of a time shift of the shape, and the efficient score, the part of $tf'(\tau t)$ orthogonal to the directions in the shape, is $(t-n/2)f'(\tau t)$, up to a term of lower order.
The efficient information for the frequency is therefore $\{1+o(1)\}n^3\|f'\|_2^2/12$, and the extra factor 4 with respect to~\eqref{eq:fisher_intro} comes from the choice of the origin of time; cf.~\cite{Golubev:1988}, \cite[Section~2 and Remark~3]{Castillo:2007} (after the change from period to frequency), \cite{Hall:2000} for the same phenomenon in random design regression, and Proposition~\ref{prop:information} below.
With the origin in the middle of the window the score for the frequency is asymptotically orthogonal to all directions in the shape, and the information with known shape is already the efficient information; in this sense, not knowing the shape costs nothing beyond not knowing its position in time.

Bayesian methods are widely used for frequency estimation in practice, from the Bayesian spectrum analysis of Bretthorst~\cite{Bretthorst:1988} and reversible jump samplers for sinusoids in noise~\cite{Andrieu:1999}, to Gaussian process models for the period of light curves~\cite{Wang:2012} and for detecting periodicities~\cite{Durrande:2016}, and to models for non-stationary periodic time series~\cite{HadjAmar:2020}.
The estimator of~\cite{Castillo:2007} is close to a Bayesian procedure, but it is a point estimator.
As far as we are aware, the frequentist properties of the posterior distribution of the frequency have not been studied when the shape is non-parametric.

The Bayesian approach to frequency estimation has some advantages.
Frequency and shape are estimated jointly, and samples from the posterior, obtained for instance with Markov chain Monte Carlo algorithms, quantify the uncertainty about both.
The frequency need not be restricted to a grid or to a known range, and conditions on the estimator of the frequency, needed for plug-in estimation, are replaced by conditions on the priors.
The role of the constraints above, which keep the frequency away from multiples and fractions of the true one, is played by the choice of the prior alone, and the prior of Section~\ref{sec:examples} needs no such knowledge.
The posterior inherits, however, the multimodality of the likelihood near multiples and fractions of the true frequency.

The theory of non-parametric Bayesian procedures has developed considerably since the articles of Ghosal, Ghosh and van der Vaart~\cite{Ghosal:2000} and Shen and Wasserman~\cite{Shen:2001}; cf.~\cite{Ghosal:2017} for an overview, and~\cite{Ghosal:2007} for non-i.i.d.\ models such as the Gaussian white noise model.
For appropriately chosen priors, the posterior gives rate-optimal estimates, and adapts to unknown smoothness when, for instance, the prior is a random series with a random number of terms~\cite{Arbel:2013,Shen:2015}.

If one is only interested in a low-dimensional functional of the model, such as the frequency, semi-parametric versions of the Bernstein-von Mises theorem describe the asymptotic behaviour of its marginal posterior~\cite{Shen:2002,Bickel:2012,Castillo:2012,Rivoirard:2012,Castillo:2015}; see~\cite{Ray:2020,Walker:2023,Franssen:2024,Giordano:2025,Duverger:2026} for recent developments.
Under such a theorem the marginal posterior is asymptotically Gaussian, centred at an efficient estimator, with variance equal to the inverse of the efficient information, so that credible sets are asymptotically efficient confidence sets.
These results concern functionals that are estimable at rate $n^{-1/2}$, in models where the parameter of interest and the nuisance parameter carry information of the same order in $n$.
In parametric models, local asymptotic normality (LAN), the approximation of the log-likelihood ratio by that of a Gaussian shift experiment, is classical with other normalising rates, as are the corresponding results for Bayes estimators~\cite{Ibragimov:1981}; for instance, the maximum likelihood and Bayes estimators of the frequency of a periodic drift of known shape are asymptotically normal at rate $n^{-3/2}$~\cite{Hopfner:2012}.
Posterior limits at rates other than $n^{-1/2}$ have also been obtained in non-regular models, where the limits need not be Gaussian~\cite{Ghosal:1995,Bochkina:2014}, including semi-parametric ones~\cite{Kleijn:2012}.

In this paper we study the posterior distribution of the frequency and of the shape in model~\eqref{eq:gwn_intro}, when frequency and shape are endowed with independent priors and the prior on the shape is a random series in the trigonometric basis.
We denote the posterior by $\Pi(\cdot\mid X^{(n)})$.
Our results are of three kinds.
First, we give conditions on the priors under which the posterior for the signal contracts around $\theta_0$, in the $L^2$ norm on $[0,n]$ normalised by $n$, at a rate $\epsilon_n$.
For the priors of Section~\ref{sec:examples}, this is the minimax rate $n^{-\alpha/(2\alpha+1)}$ for shapes of Sobolev smoothness $\alpha$, up to a logarithmic factor, attained adaptively when the number of terms in the series is random.

Second, we show that the marginal posterior for the frequency contracts around $\tau_0$ at rate $\epsilon_n/n$, and that the marginal posterior for the shape contracts around $f_0$ at rate $\epsilon_n$ in $L^2\big([0,1]\big)$; see Theorem~\ref{theo:checking_conditions_2}.
This requires an identifiability condition.
The pairs $(\tau_0,f_0)$ and $(\tau_0/p,f_0(p\,\cdot))$ give the same signal $\theta_0$ for every integer $p\ge2$, so shapes with a shorter period must be excluded and, for the posterior, also shapes that are arbitrarily close to having one; we impose this through a condition on the prior.
The localisation of the frequency is deterministic and holds for every frequency ratio $\tau/\tau_0$: a signal close to $\theta_0$ is close to its own shifts by multiples of $1/\tau_0$, which forces $\tau/\tau_0$ to be close to a rational number $p/q$; the condition on the shapes excludes $q\ge2$, and the fact that $f_0$ has smallest period 1 excludes $p\ge2$ (Proposition~\ref{prop:smaller_sieve}).
The rate $\epsilon_n/n$ is slower than $n^{-3/2}$, since a frequency error $\delta$ shifts the signal by about $n\delta$ periods at the end of the window, so closeness of signals in $L^2$ controls $n\delta$ but not more; it suffices, however, to localise the frequency for the third result.

Third, we show that the marginal posterior for the frequency satisfies a semi-parametric Bernstein-von Mises theorem at rate $n^{-3/2}$: for shapes of Sobolev smoothness $\alpha>3/2$ and suitable stationary priors (Theorem~\ref{theo:BvM} and Corollary~\ref{cor:BvM_example}),
\begin{equation}\label{eq:BvM_intro}
\sup_{B}\bigg|\Pi\Big(n^{3/2}(\tau-\hat\tau_n)\in B\;\Big|\; X^{(n)}\Big)-N\Big(0,\frac{12}{\|f_0'\|_2^2}\Big)(B)\bigg|\to0,
\end{equation}
in $P_{\theta_0}^{(n)}$-probability, as $n\to\infty$, where the supremum is over all Borel sets $B$ and $\hat\tau_n$ is an efficient estimator of $\tau_0$, i.e., any estimator such that
\[
n^{3/2}(\hat\tau_n-\tau_0)=\frac{12}{n^{3/2}\|f_0'\|_2^2}\int_0^n\Big(t-\frac n2\Big)f_0'(\tau_0t)\,dW_t+o_P(1).
\]
In particular, credible intervals for the frequency are asymptotically efficient confidence intervals at the scale $n^{-3/2}$.
The rate $n^{-3/2}$ alone does not make this result new, since semi-parametric theorems such as~\cite[Theorem~1]{Castillo:2012} hold with a similar proof at other normalising rates (Lemma~\ref{lemma:BvM_general}).
What is new is that the frequency, a fixed parameter that does not depend on $n$, and the shape carry information of different orders, $n^3$ and $n$, and that the lack of orthogonality between the two comes from the choice of the origin of time, which allows us to remove it, up to terms of lower order, by a reparametrisation that does not change the posterior, as explained below.
We are not aware of other semi-parametric Bernstein-von Mises theorems with a Gaussian limit in which a fixed parameter of interest and the nuisance parameter carry information of different orders.

The proof of~\eqref{eq:BvM_intro} rests on the version of~\cite[Theorem~1]{Castillo:2012} at rate $n^{-3/2}$ given by Lemma~\ref{lemma:BvM_general}.
What has to be shown is that the log-likelihoods at $(\tau,f)$ and at $(\tau_0,f)$ differ by a quadratic function of $s=n^{3/2}(\tau-\tau_0)$, uniformly over the frequencies and shapes on which the posterior concentrates, up to an error that is small relative to $1+s^2$.
Since we localise the posterior on a product of a set of frequencies and a set of shapes, we do not need the concentration of the posterior with known frequency assumed in~\cite{Castillo:2012}; in exchange, the comparison has to hold on the whole product.
It meets two difficulties.
The first is specific to the frequency: a change of frequency deforms the signal more and more along the observation window, so the remainder terms in the expansion of the log-likelihood involve the derivatives of the shapes, at frequencies up to distance $\epsilon_n/n$ from $\tau_0$.
We control them using the smoothness $\alpha>3/2$ of the shape and a bound on stochastic integrals that is uniform over the set on which the posterior concentrates; the approach of Bickel and Kleijn~\cite{Bickel:2012} would instead require the marginal posterior for the frequency to be localised at rate $n^{-3/2}$ beforehand, which our first two results do not give.

The second difficulty is the lack of orthogonality.
With the origin of time at the start of the observation window, as in~\eqref{eq:gwn_intro}, the least favourable direction for the shape, the direction that reduces the information about the frequency the most, is $g=-sf_0'/2$ in the local parameters
\[
\tau=\tau_0+sn^{-3/2},\qquad f=f_0+gn^{-1/2}.
\]
Semi-parametric Bernstein-von Mises theorems allow for this only under conditions on the prior in that direction.
Castillo~\cite[Theorem~2]{Castillo:2012} covers Gaussian process priors, through a change of variables that moves the shape along an approximation of the least favourable direction, and~\cite{Castillo:2015} requires the prior to be approximately unchanged by a perturbation along the representer of the functional.
Here a perturbation $f_0+cf_0'$ of the shape along $f_0'$, with $c\in\mathbb{R}$, is, up to terms of order $c^2$, a time shift of $f_0$, but not of a generic $f$, so even a prior on the shape that is unchanged by time shifts does not satisfy these conditions directly.

We use instead the following.
If the frequency and the shape are independent under the prior, and the prior on the shape is \emph{stationary}, in the sense that $f(\cdot+c)$ has the same prior distribution as $f$ for every $c\in\mathbb{R}$ (see Section~\ref{sec:estimation} for the precise definition), then the prior on $(\tau,\theta)$, and therefore the marginal posterior for the frequency, is the same in every parametrisation
\[
\theta(t)=f\big(\tau(t-c)\big),\qquad t\in[0,n],
\]
with the origin of time moved to $c\in\mathbb{R}$.
We may therefore place the origin in the middle of the observation window, where the score for the frequency is asymptotically orthogonal to the shape, and prove the Bernstein-von Mises theorem in that parametrisation.
The price is that the true shape in the centred parametrisation, $f_0(\cdot+\tau_0n/2)$, depends on $n$, so we localise the posterior in the original parametrisation and transfer the localisation to the centred one.
A reparametrisation that removes the loss of information so that independent priors can be used appears in~\cite{Walker:2023}, for the partially linear model, where it changes the prior; here, under stationarity, the two parametrisations induce the same prior, so the result concerns the posterior in the natural parametrisation.
Centring the observation window at the origin, as in~\cite{Golubev:1988,Castillo:2007}, is not free for a Bayesian procedure, since moving the origin of time changes the prior in general; the invariance of the posterior under this change for stationary priors (Proposition~\ref{prop:origin}), and its use to remove the loss of information, are to our knowledge new.

Gaussian random series priors with independent centred coefficients and equal variances for the two coefficients of each harmonic, as in~\cite[Section~2]{Castillo:2007}, are stationary, but, when $f_0$ has a non-zero harmonic of index at least 2, no Gaussian random series prior with independent centred coefficients satisfies our identifiability condition together with the other conditions of our second result (Proposition~\ref{prop:gaussian_fails} in Appendix~\ref{sec:technical:examples}).
In Section~\ref{sec:examples} we condition a Gaussian random series prior on a set of shapes that depends only on the amplitudes of the harmonics, which keeps it stationary, and show that it satisfies the conditions of our first two results for every true shape with smallest period 1 and Sobolev smoothness $\alpha\ge1$, adaptively when the number of terms is random, and those of the third when, moreover, $\alpha>3/2$.
We work with the Gaussian white noise model, where the computations are the most transparent; other models with periodic features, such as Poisson processes, are discussed in Section~\ref{sec:discussion}.

Section~\ref{sec:preliminaries} introduces the notation and Section~\ref{sec:statement_of_problem} the model and the information about the frequency; Sections~\ref{sec:estimation} to~\ref{sec:BVM} contain the results above, and Section~\ref{sec:discussion} is a discussion.
The proofs are in Appendices~\ref{sec:technical:statement} to~\ref{sec:technical:BVM}.

\section{Preliminaries}\label{sec:preliminaries}

In this section we introduce some notation and basic definitions that will be used throughout this paper.

In what follows, $\mathbb{N}=\{1,2,\dots\}$, $I\subseteq\mathbb{R}$ is an interval, and $M>0$ is a constant.
For $x\in\mathbb{R}$, $\lfloor x\rfloor$ is the largest integer smaller than or equal to $x$, and $\lceil x\rceil$ the smallest integer larger than or equal to $x$.
For positive sequences $a_n$ and $b_n$, $a_n\lesssim b_n$ and $b_n\gtrsim a_n$ mean that $a_n\le Cb_n$ for all $n$, for a constant $C$ that does not depend on $n$.
We write $a_n\asymp b_n$ if $a_n\lesssim b_n$ and $b_n\lesssim a_n$, and $a_n=o(b_n)$ if $a_n/b_n\to0$, as $n\to\infty$.
We write $N(\mu,\sigma^2)$ for the normal distribution with mean $\mu$ and variance $\sigma^2$, $\Phi$ for the standard normal distribution function, and $\chi_d^2$ for the chi-square distribution with $d$ degrees of freedom.
We write $\mathbb{P}$ and $\mathbb{E}$ for the probability and the expectation of events and variables that involve only auxiliary random variables, such as a standard normal variable.
The notation for the model ($X^{(n)}$, $P_\theta^{(n)}$, $\theta_0$, $\tau_0$, $f_0$), used informally in Section~\ref{sec:introduction}, is made precise in Section~\ref{sec:statement_of_problem}.
We denote by $L^2\big(I\big)$ the space of square integrable functions on $I$, and by $\ell^2$ the space of square summable doubly infinite sequences:
\[
L^2\big(I\big) = \Big\{ f:I\to\mathbb{R}: \int_I \big\{f(t)\big\}^2\,dt < \infty \Big\}, \;\;
\ell^2 = \Big\{ \bm{x}=(x_i, \, i\in\mathbb{Z}): \sum_{i\in\mathbb{Z}} x_i^2 < \infty \Big\}.
\]
We also consider some subsets of these collections, namely the Sobolev space of $\alpha$-smooth functions, the periodic Sobolev space of $\alpha$-smooth functions, and the Sobolev ellipsoid of $\alpha$-regular sequences, respectively,
\begin{align*}
\mathcal{S}_\alpha\big(M\big) &= \Big\{ f:[0,1]\to\mathbb{R}: f^{(\alpha-1)}\in AC\big([0,1]\big), \; \int_0^1 \big\{f^{(\alpha)}(t)\big\}^2\,dt \le M^2 \Big\}, \;\; \alpha\in\mathbb{N},\\
\mathcal{W}_\alpha\big(M\big) &= \Big\{ f\in\mathcal{S}_\alpha(M): f^{(j)}(0)=f^{(j)}(1),\; j=0,\dots,\alpha-1 \Big\},\;\; \alpha\in\mathbb{N},\\
w_\alpha\big(M\big) &= \Big\{ \bm{x}\in\ell^2: \sum_{i\in\mathbb{Z}^*} v_i^{2\alpha} x_i^2 \le M^2 \Big\}, \;\; \alpha>0,
\end{align*}
where $AC(I)$ is the space of absolutely continuous functions on $I$, and $f^{(m)}$ is the $m$-th derivative of $f$, with $f^{(0)}=f$ and $f'=f^{(1)}$.
In the definition of $w_\alpha(M)$, $\mathbb{Z}^*=\mathbb{Z}\setminus\{0\}$ and $v_i=2\pi|i|$, and no condition is imposed on $x_0\in\mathbb{R}$.
On $L^2\big([0,n]\big)$, $L^2\big([0,1]\big)$ and $\ell^2$ we consider the norms
\[
\|f\|_{2,n}^2 = \frac1n\int_0^n \big\{f(t)\big\}^2\,dt, \;\;
\|f\|_{2}^2 = \|f\|_{2,1}^2, \;\;
\|\bm{x}\|_2^2=\sum_{i\in\mathbb{Z}} x_i^2, \;\; \text{and} \;\;
\|\bm{x}\|_\infty=\sup_{i\in\mathbb{Z}} |x_i|.
\]

The spaces $\mathcal{W}_\alpha\big(M\big)$ and $w_\alpha\big(M\big)$ are related in the following way.
The trigonometric orthonormal basis of $L^2\big([0,1]\big)$, extended to $\mathbb{R}$ by periodicity, is defined as
\begin{equation*}
\varphi_0(t) = 1, \quad \varphi_{i}(t) = \sqrt2\cos(2\pi i t), \quad \varphi_{-i}(t) = \sqrt2\sin(2\pi i t),
\quad i\in\mathbb{N}, \, t\in\mathbb{R}.
\end{equation*}
We abbreviate the basis as $\bm{\varphi}=\big(\varphi_i, \, i\in\mathbb{Z}\big)$.
A function $f\in L^2\big([0,1]\big)$ admits a representation\footnote{If $\bm{x}=(x_{-m},\dots,x_m)\in\mathbb{R}^{2m+1}$ for some $m\in\mathbb{N}$, and we write $\bm{x}^T\bm{\varphi}(t)$, then we see $\bm{x}$ as a sequence in $\ell^2$ with $x_{\pm(m+1)},x_{\pm(m+2)},\dots$ equal to zero.}
\[
f(t) = \sum_{i\in\mathbb{Z}} x_i \varphi_i(t) = \bm{x}^T\bm{\varphi}(t), \quad t\in[0,1],
\]
in $L^2\big([0,1]\big)$, with $x_i=\int_0^1 f(t)\varphi_i(t)\,dt$, so that $\bm{x}\in\ell^2$ and $\|f\|_2=\|\bm{x}\|_2$.
The relation between the two spaces is then as follows: a function $f$ as above belongs to $\mathcal{W}_\alpha(M)$ if, and only if, the corresponding sequence $\bm{x}=(x_i, \, i\in\mathbb{Z})$ belongs to $w_\alpha(M)$; cf.~\cite[Proposition~1.14]{Tsybakov:2009}.
We call $\bm x$ the coordinates of $f$.
The coordinate $x_0=\int_0^1f(t)\,dt$ is the mean of $f$, which we also denote by $\bar f$, so that $\|f-\bar f\|_2^2=\sum_{i\in\mathbb{Z}^*}x_i^2$.
We write $\mathcal{W}_\alpha=\bigcup_{M>0}\mathcal{W}_\alpha(M)$ and $w_\alpha=\bigcup_{M>0}w_\alpha(M)$.

Because of this relation we identify a function $f\in\mathcal{W}_\alpha$ with its coordinates $\bm x\in w_\alpha$ and we may also write $\bm{x}\in\mathcal{W}_\alpha$ or $f\in w_\alpha$.
Since $w_\alpha$ is defined for any $\alpha>0$, we say that $f\in\mathcal{W}_\alpha$, $\alpha>0$, if the coordinates of $f$ are in $w_\alpha$.

For a function $f(t)$ on $\mathbb{R}$, if we write $f\in\mathcal{W}_\alpha$, we mean that $f$ is 1-periodic, and the restriction of $f$ to the interval $[0,1]$ belongs to $\mathcal{W}_\alpha$.
Conversely, the restriction to $[0,1)$ of a function $f\in\mathcal{W}_\alpha$ can be extended to a 1-periodic function on $\mathbb{R}$.

For $k\in\mathbb{N}$ we call $a_k=(x_k^2+x_{-k}^2)^{1/2}$ the \emph{amplitude of the $k$-th harmonic} of $f$.
If $f$ has coordinates $\bm x$, then the coordinates of the time shift $f(\cdot+c)$, $c\in\mathbb{R}$, are obtained by rotating each pair $(x_k,x_{-k})$, $k\in\mathbb{N}$, through the angle $-2\pi kc$, and $x_0$ is unchanged.
In particular, the amplitudes of the harmonics, the norm $\|f\|_2$ and membership of $\mathcal{W}_\alpha(M)$ are unchanged by time shifts, and, since a rotation through an angle $\phi$ moves a vector of length $a$ by the distance $2a|\sin(\phi/2)|$,
\begin{equation}\label{eq:shift}
\|f(\cdot+c)-f\|_2^2=\sum_{k=1}^\infty a_k^2\,4\sin^2(\pi kc),\quad c\in\mathbb{R}.
\end{equation}
For $p\in\mathbb{N}$ we have $\varphi_i(p\,\cdot)=\varphi_{pi}$, so the coordinate of $f(p\,\cdot)$ with index $pi$ is $x_i$, $i\in\mathbb{Z}$, and its coordinates with indices that are not multiples of $p$ vanish.

Let $f$ be a function on $\mathbb{R}$ that is square integrable on bounded intervals.
Functions that are equal almost everywhere are identified.
We call $c\in\mathbb{R}$ a \emph{period} of $f$ if $f(t+c)=f(t)$ for almost every $t\in\mathbb{R}$, and we call $f$ constant if it is almost everywhere equal to a constant.
If $f$ is not constant and has a positive period, then it has a smallest positive period (Lemma~\ref{lemma:periods}(ii) below), and we call the reciprocal of this period the \emph{frequency} of $f$.
A \emph{shape} is a 1-periodic function on $\mathbb{R}$ that is square integrable on $[0,1]$, which we identify with its restriction to $[0,1]$ and with its coordinates.

For $k\in\mathbb{N}$ and a shape $f$ with coordinates $\bm x$, we write
\[
f_{[k]}=\sum_{i=-k}^kx_i\varphi_i
\]
for the truncation of $f$ after the $k$-th harmonic, and we call $f$ a \emph{trigonometric polynomial of degree at most $k$} if $f=f_{[k]}$, i.e., if $x_i=0$ for $|i|>k$.
If $f\in w_\alpha(M)$, then, since $v_i\ge2\pi(k+1)$ for $|i|>k$,
\begin{equation}\label{eq:truncation}
\|f-f_{[k]}\|_2^2=\sum_{|i|>k}x_i^2\le\frac{1}{\{2\pi(k+1)\}^{2\alpha}}\sum_{|i|>k}v_i^{2\alpha}x_i^2\le\frac{M^2}{\{2\pi(k+1)\}^{2\alpha}}.
\end{equation}

\section{Statement of the problem}\label{sec:statement_of_problem}

Consider a sequence of statistical experiments $\big(\mathcal{X}^{(n)},\mathcal{U}^{(n)},\mathcal{P}^{(n)}\big)$, $n\in\mathbb{N}$.
The non-parametric model $\mathcal{P}^{(n)}=\{P_\theta^{(n)}:\theta\in\Theta\}$ is assumed to have periodic features, in the sense that the signals $\theta$ form a set $\Theta$ of periodic functions, made precise in Section~\ref{sec:statement_of_problem:parametrisation}.
We observe $X^{(n)}$, distributed according to $P_{\theta_0}^{(n)}$, for some $\theta_0\in\Theta$.
Our main goal is to make semi-parametric inference on the frequency $\tau_0$ of $\theta_0$.

Here we use the Gaussian white noise model as a prototype for $\mathcal{P}^{(n)}$.
We observe $X^{(n)}=\big(Y_t:t\in[0,n]\big)$, where $Y_0=0$ and
\begin{equation}\label{eq:def_gwn}
dY_t=\theta(t)\,dt+dW_t,\quad t\in[0,n],
\end{equation}
for a standard Brownian motion $W$ and a signal $\theta$ that is square integrable on $[0,n]$; model~\eqref{eq:gwn_intro} is the case $\theta=f(\tau\,\cdot)$.
The sample space $\mathcal{X}^{(n)}$ is the space $C\big([0,n]\big)$ of continuous functions on $[0,n]$ with its Borel $\sigma$-algebra $\mathcal{U}^{(n)}$, and we denote by $\nu^{(n)}$ the law of $W$ on this space, which is the law of $X^{(n)}$ when $\theta\equiv0$.
The measures $P_\theta^{(n)}$ are equivalent to $\nu^{(n)}$, with densities
\begin{equation}\label{eq:likelihood}
p_\theta^{(n)}\big(X^{(n)}\big)=\frac{dP_\theta^{(n)}}{d\nu^{(n)}}\big(X^{(n)}\big)=\exp\bigg\{\int_0^n\theta(t)\,dY_t-\frac12\int_0^n\big\{\theta(t)\big\}^2\,dt\bigg\};
\end{equation}
cf.~\cite[Appendix~II, Theorem~1]{Ibragimov:1981}.
Since $dY_t=\theta_0(t)\,dt+dW_t$ under $P_{\theta_0}^{(n)}$, the log-likelihood ratio is
\begin{equation}\label{eq:log_lik_ratio}
\log\frac{p_\theta^{(n)}}{p_{\theta_0}^{(n)}}\big(X^{(n)}\big)=\int_0^n\big\{\theta(t)-\theta_0(t)\big\}\,dW_t-\frac n2\|\theta-\theta_0\|_{2,n}^2,
\end{equation}
where $\|\cdot\|_{2,n}$ is applied to the restriction to $[0,n]$.
Under $P_{\theta_0}^{(n)}$ the log-likelihood ratio is therefore Gaussian, with mean and variance respectively
\[
-\frac n2\|\theta-\theta_0\|_{2,n}^2\qquad\text{and}\qquad n\|\theta-\theta_0\|_{2,n}^2,
\]
and the Kullback-Leibler divergence of $P_\theta^{(n)}$ from $P_{\theta_0}^{(n)}$ is minus the mean.

\subsection{Frequency and shape}\label{sec:statement_of_problem:parametrisation}

For the semi-parametric problem of estimating the frequency $\tau_0$ of the signal $\theta_0$, it is convenient to parametrise the model in terms of a frequency and a shape.
For $\tau>0$ and a shape $f$, write
\begin{equation}\label{eq:def_psi}
\psi(\tau,f)(t)=f(\tau t),\quad t\in\mathbb{R}.
\end{equation}
The map $\psi$ is not one-to-one, since, for every integer $p\ge2$, $f(p\,\cdot)$ is again a shape and
\[
\psi\big(\tau/p,f(p\,\cdot)\big)=\psi(\tau,f).
\]
Also, the frequency of $\psi(\tau,f)$ itself need not be $\tau$.
The following lemma specifies how a shape $f$ and frequency $\tau$ determine the frequency of $\psi(\tau,f)$.
The proofs of this and of the other results of this section are in Appendix~\ref{sec:technical:statement}.

\begin{lemma}\label{lemma:periods}
The following statements hold:
\begin{enumerate}[label=(\roman*)]
\item Let $f$ be a shape that is not constant, with coordinates $\bm x$, and let $d$ be the greatest common divisor of the set $\{k\in\mathbb{N}: x_k^2+x_{-k}^2>0\}$ of the indices of the non-zero harmonics of $f$.
Then the set of periods of $f$ is $\{j/d: j\in\mathbb{Z}\}$, and, for every $\tau>0$, the frequency of $\psi(\tau,f)$ is $d\tau$.
\item Let $\theta$ be a function on $\mathbb{R}$ that is square integrable on bounded intervals and not constant.
If $\theta$ has a positive period, then it has a smallest positive period.
\item If two functions as in (ii), with frequencies $\tau$ and $\tau'$, are equal almost everywhere on $[0,1/\tau+1/\tau']$, then they are equal.
\end{enumerate}
\end{lemma}

Denote by $\mathcal{F}$ the set of shapes that are not constant and have smallest positive period 1, i.e., with $d=1$ in Lemma~\ref{lemma:periods}(i), and by $\Theta$ the set of functions on $\mathbb{R}$ that are square integrable on bounded intervals, are not constant and have a positive period.
By Lemma~\ref{lemma:periods}, the map $\psi:(0,\infty)\times\mathcal{F}\to\Theta$ is one-to-one and onto: a signal $\theta\in\Theta$ with frequency $\tau$ equals $\psi(\tau,f)$ for $f=\theta(\cdot/\tau)$, and for no other pair in $(0,\infty)\times\mathcal{F}$.
To ease the notation, we identify $\Theta$ with $(0,\infty)\times\mathcal{F}$ and simply write $\theta=(\tau,f)$, which should be understood as $\theta=\psi(\tau,f)$.
By Lemma~\ref{lemma:periods}(iii), two different signals in $\Theta$ with frequencies $\tau$ and $\tau'$ differ on a subset of $[0,n]$ of positive measure, and therefore give different laws $P_\theta^{(n)}$, as soon as $n\ge1/\tau+1/\tau'$.
The model is therefore identifiable in the parametrisation $(0,\infty)\times\mathcal{F}$: every two parameters give different laws for all $n\ge1/\tau+1/\tau'$.

Throughout, the data are generated by a true signal $\theta_0\in\Theta$, with frequency $\tau_0$ and shape $f_0=\theta_0(\cdot/\tau_0)\in\mathcal{F}$, so that $\theta_0(t)=f_0(\tau_0t)$.
We denote the coordinates of $f_0$ by $x_{0,i}$, $i\in\mathbb{Z}$, and write $\sigma_0=\|f_0-\bar f_0\|_2$, which is positive since $f_0$ is not constant.
Further assumptions on $f_0$ are stated where they are needed.

Identifiability is not enough as shapes in $\mathcal{F}$ can be arbitrarily close to shapes with a shorter period.
For instance, for $\eta>0$ consider the shape
\[
f_\eta=f_0(2\,\cdot)+\eta\varphi_1.
\]
Its first harmonic is non-zero, so $f_\eta\in\mathcal{F}$ by Lemma~\ref{lemma:periods}(i), and the signal $(\tau_0/2,f_\eta)$ has frequency $\tau_0/2$.
Yet, since $|\varphi_1|\le\sqrt2$,
\begin{equation}\label{eq:f_eta_close}
\psi(\tau_0/2,f_\eta)=\theta_0+\eta\varphi_1(\tau_0\,\cdot/2)
\quad\text{and}\quad
\|\psi(\tau_0/2,f_\eta)-\theta_0\|_{2,n}\le\sqrt2\,\eta,
\quad n>0.
\end{equation}
For small $\eta$ the signal $(\tau_0/2,f_\eta)$ is therefore close to $\theta_0$, while its frequency is half of $\tau_0$: a bound on $\|\theta-\theta_0\|_{2,n}$ alone does not localise the frequency.

To exclude such shapes we use the sets
\begin{equation}\label{eq:xi_set_definition}
\Xi(h)=\Big\{\bm{x}\in\ell^2: \sum_{i\in\mathbb{Z}^*}x_i^2>0,\;\sum_{i\in\mathbb{Z}^*}x_{pi}^2\le(1-h)\sum_{i\in\mathbb{Z}^*}x_i^2,\;p=2,3,\dots\Big\},\quad h\in(0,1],
\end{equation}
which are adapted from those of Golubev~\cite{Golubev:1988}.
The sum $\sum_{i\in\mathbb{Z}^*}x_{pi}^2$ is the sum of the squared amplitudes of the harmonics with indices $p,2p,3p,\dots$, so the sets $\Xi(h)$ depend on $\bm x$ only through the amplitudes of the harmonics, and, by Section~\ref{sec:preliminaries}, they are unchanged by time shifts.

The condition in~\eqref{eq:xi_set_definition} has a simple meaning in terms of the shape $f$ with coordinates $\bm x$.
The sum $\sum_{i\in\mathbb{Z}^*}x_i^2$ is $\|f-\bar f\|_2^2$, the squared distance from $f$ to its mean $\bar f= x_0$ (Section~\ref{sec:preliminaries}), i.e., to the constant shapes.
For $p\in\mathbb{N}$, let
\[
A_pf(u)=\frac1p\sum_{j=0}^{p-1}f\Big(u+\frac jp\Big),\quad u\in\mathbb{R},
\]
be the average of $f$ over its shifts by multiples of $1/p$.
By Section~\ref{sec:preliminaries}, the shift by $j/p$ rotates each pair $(x_k,x_{-k})$ through the angle $-2\pi kj/p$, and the average over $j=0,\dots,p-1$ of these rotations is the identity if $p$ divides $k$, and zero otherwise.
Hence $A_pf$ has coordinates $x_i$ for the indices $i$ divisible by $p$, and zero for the other indices.
By Lemma~\ref{lemma:periods}(i), the shapes with period $1/p$ are those whose coordinates vanish outside the multiples of $p$ (constants included), so $A_pf$ is the orthogonal projection in $L^2\big([0,1]\big)$ of $f$ onto these shapes, and
\[
\sum_{i\in\mathbb{Z}^*}x_{pi}^2=\|A_pf-\bar f\|_2^2=\|f-\bar f\|_2^2-\|f-A_pf\|_2^2.
\]
A shape $f$ is therefore in $\Xi(h)$ if, and only if, it is not constant and
\begin{equation}\label{eq:xi_set_functions}
\|f-A_pf\|_2^2\ge h\|f-\bar f\|_2^2,\quad p=2,3,\dots.
\end{equation}
In words: for every $p\ge2$, the squared distance from $f$ to the shapes with period $1/p$ is at least a fraction $h$ of the squared distance from $f$ to the constants.
For the shape $f_\eta$ above, $A_2f_\eta=f_0(2\,\cdot)$ and $\bar f_\eta=\bar f_0$, so that
\[
\|f_\eta-A_2f_\eta\|_2^2=\eta^2,\qquad\|f_\eta-\bar f_\eta\|_2^2=\sigma_0^2+\eta^2,
\]
and $f_\eta\in\Xi(h)$ only if $h\le\eta^2/(\sigma_0^2+\eta^2)$.

\begin{lemma}\label{lemma:xi_sets}
The following statements hold:
\begin{enumerate}[label=(\roman*)]
\item If $0<h\le h'\le1$, then $\Xi(h')\subseteq\Xi(h)$.
\item A shape belongs to $\mathcal{F}$ if, and only if, its coordinates belong to $\Xi(h)$ for some $h\in(0,1]$.
\item If the amplitude $a$ of the first harmonic of $\bm x\in\ell^2$ is positive, then $\bm x\in\Xi(h)$ for $h=a^2/\sum_{i\in\mathbb{Z}^*}x_i^2$.
\item If $\bm x\in\Xi(h)$ and $\bm y$ are the coordinates of a shape that is constant or has period $1/p$ for some integer $p\ge2$, then $\|\bm x-\bm y\|_2^2\ge h\sum_{i\in\mathbb{Z}^*}x_i^2$.
In particular, if $f$ is a shape with coordinates in $\Xi(h)$, then every shape $g$ with $\|f-g\|_2^2<h\|f-\bar f\|_2^2$ is in $\mathcal{F}$, and $\mathcal{F}$ is an open subset of $L^2\big([0,1]\big)$.
\end{enumerate}
\end{lemma}

We use two consequences of Lemma~\ref{lemma:xi_sets} for the true shape.
First, since $f_0\in\mathcal{F}$, Lemma~\ref{lemma:xi_sets}(ii) gives $f_0\in\Xi(H)$ for some $H\in(0,1]$, and then, by Lemma~\ref{lemma:xi_sets}(iv), every shape $f$ with $\|f-f_0\|_2^2<H\sigma_0^2$ is in $\mathcal{F}$.

Second, the sets $\Xi(h)$ keep the signals with frequency exactly $\tau_0/p$, $p\ge2$, away from $\theta_0$.
Let $p\ge2$ be an integer such that $n\tau_0\ge2p$, and let $f\in\Xi(h)$.
Then
\begin{equation}\label{eq:xi_separation}
\|\psi(\tau_0/p,f)-\theta_0\|_{2,n}^2\ge\Big(1-\frac p{n\tau_0}\Big)h\|f-\bar f\|_2^2.
\end{equation}
Indeed, $\theta_0=\psi\big(\tau_0/p,f_0(p\,\cdot)\big)$, so $\psi(\tau_0/p,f)-\theta_0=\psi(\tau_0/p,g)$ for the shape $g=f-f_0(p\,\cdot)$.
The first bound in~\eqref{eq:information_1} of Proposition~\ref{prop:information} below, with $\tau_0/p$ in place of $\tau$, and Lemma~\ref{lemma:xi_sets}(iv), since $f_0(p\,\cdot)$ has period $1/p$, give
\[
\|\psi(\tau_0/p,g)\|_{2,n}^2\ge\Big(1-\frac p{n\tau_0}\Big)\|g\|_2^2\ge\Big(1-\frac p{n\tau_0}\Big)h\|f-\bar f\|_2^2.
\]
For $f=f_\eta$, with the largest $h$ for which $f_\eta\in\Xi(h)$, the right side of~\eqref{eq:xi_separation} is $\{1-2/(n\tau_0)\}\eta^2$, which matches the upper bound $2\eta^2$ from~\eqref{eq:f_eta_close} up to a constant factor.
Frequencies near, but not equal to, $\tau_0/p$, and more generally frequencies near rational multiples of $\tau_0$, are treated in Proposition~\ref{prop:smaller_sieve}.

To exclude shapes such as $f_\eta$ with $\eta$ small, we work with shapes in $\Xi(h)$ for some $h>0$.
Some restriction of this kind seems unavoidable, at least for results that are uniform in the truth: the Kullback-Leibler divergence of $P_{(\tau_0/2,f_\eta)}^{(n)}$ from $P_{\theta_0}^{(n)}$ is at most $n\eta^2$, so for $\eta=n^{-1/2}$ the two laws cannot be told apart consistently, while the frequencies are $\tau_0/2$ and $\tau_0$.
In Section~\ref{sec:estimation} we impose the restriction through a condition on the prior, with $h=h_n$ allowed to decrease to zero with $n$, so that $H$ need not be known; the prior of Section~\ref{sec:examples} is conditioned on $\Xi(h_n)$.
By Lemma~\ref{lemma:xi_sets}(i) and~(iii), a shape $f$ whose first harmonic has amplitude $a>0$ is in $\Xi(h)$ for every $h\le a^2/\|f-\bar f\|_2^2$; cf.~\cite[Section~3.1]{Castillo:2007}, where a lower bound on the first Fourier coefficient and a Sobolev bound give a condition of the type~\eqref{eq:xi_set_definition}.
A non-zero first harmonic is sufficient but not necessary: for instance, the shape $\cos(4\pi u)+\cos(6\pi u)$, $u\in\mathbb{R}$, is in $\Xi(1/2)$ and has no first harmonic.

The estimation of a shape $f\in\mathcal{F}$ from the data $X^{(n)}$ when the frequency is known is rather standard and reduces to the estimation of a non-parametric curve on a bounded interval, at a rate that depends on the smoothness of $f$; cf.~\cite{Tsybakov:2009}.
The estimation of the frequency is less standard, since, as we now show, the information about it grows like $n^3$.

\subsection{Information about the frequency}\label{sec:statement_of_problem:information}

Let $f\in\mathcal{F}\cap\mathcal{W}_1$ and $\tau>0$.
For fixed $f$, the map $\tau\mapsto\psi(\tau,f)$ is differentiable in $L^2\big([0,n]\big)$, with derivative $tf'(\tau t)$ (see Appendix~\ref{sec:technical:statement}), and the Fisher information about $\tau$ when $f$ is known is
\[
\mathcal{I}_n(\tau,f)=\int_0^n t^2\big\{f'(\tau t)\big\}^2\,dt;
\]
cf.~\cite[p.~209]{Ibragimov:1981}.
A frequency error $\delta$ moves the signal at time $t$ by about $t\delta f'(\tau t)$, so late observations carry more information about the frequency, and the weight $t^2$ makes the information grow like $n^3$.
When the shape is unknown, it can be perturbed in the direction of any shape $g$, since $\mathcal{F}$ is open by Lemma~\ref{lemma:xi_sets}(iv); this perturbs the signal in the direction $g(\tau\,\cdot)$.
We call the squared distance in $L^2\big([0,n]\big)$ from the score to these directions,
\[
\tilde{\mathcal{I}}_n(\tau,f)=\inf_{g}\int_0^n\big\{tf'(\tau t)-g(\tau t)\big\}^2\,dt,
\]
where the infimum is over all shapes $g$, the efficient information about $\tau$; cf.~\cite[Section~25.4]{vanderVaart:1998} for the definition in i.i.d.\ models.
Its interpretation as a bound for estimators of the frequency comes from the LAN expansion of Section~\ref{sec:BVM}.
Recall that shapes are 1-periodic, so the directions $g(\tau\,\cdot)$ have period $1/\tau$, while the score $tf'(\tau t)$ has an amplitude that grows linearly in $t$.
A single periodic direction cannot follow this growth: the best it can do is to match, at each phase, the average of $tf'(\tau t)$ over the periods in $[0,n]$, which is about $(n/2)f'(\tau t)$.
What remains, about $(t-n/2)f'(\tau t)$, has squared norm of order $n^3$, and this is the efficient information (Proposition~\ref{prop:information} below and its proof).

The Fisher information depends on where the origin of time is placed, but the efficient information does not.
For $c\in\mathbb{R}$, write
\[
\psi_c(\tau,f)(t)=f\big(\tau(t-c)\big),\quad t\in\mathbb{R},
\]
so that $\psi=\psi_0$ and $\psi_c\big(\tau,f(\cdot+\tau c)\big)=\psi(\tau,f)$.
A time shift does not change the periods of a shape, so $(\tau,f)\mapsto\big(\tau,f(\cdot+\tau c)\big)$ maps $(0,\infty)\times\mathcal{F}$ one-to-one onto itself, and $\psi_c$ is another parametrisation of $\Theta$, with the origin of time moved to $c$.
The two parametrisations differ in what is held fixed when $\tau$ varies: in $\psi_c$ it is the phase of the signal at time $c$.
In $\psi_c$, the derivative of the signal with respect to $\tau$ for fixed shape, at the point $\big(\tau,f(\cdot+\tau c)\big)$ that corresponds to $(\tau,f)$, is $(t-c)f'(\tau t)$, and the Fisher information at this point is
\[
\mathcal{I}_n^{(c)}(\tau,f)=\int_0^n(t-c)^2\big\{f'(\tau t)\big\}^2\,dt.
\]
A perturbation of the shape in the direction $g$ perturbs the signal in the direction $g\big(\tau(\cdot-c)\big)$, which is $\tilde g(\tau\,\cdot)$ for $\tilde g=g(\cdot-\tau c)$, so the directions in the shape are again the functions $g(\tau\,\cdot)$, for shapes $g$.
Since $(t-c)f'(\tau t)$ differs from $tf'(\tau t)$ by the direction $-cf'(\tau t)$, the efficient information at this point is $\tilde{\mathcal{I}}_n(\tau,f)$ for every $c$.
In particular, $\tilde{\mathcal{I}}_n(\tau,f)\le\mathcal{I}_n^{(c)}(\tau,f)$ for every $c\in\mathbb{R}$.

\begin{proposition}\label{prop:information}
Let $f\in\mathcal{W}_1$ be not constant, and let $\tau>0$ and $n>0$ be such that $n\tau\ge2$.
Then, for every shape $g$ with $\|g\|_2>0$,
\begin{align}
\bigg|\frac{\|g(\tau\,\cdot)\|_{2,n}^2}{\|g\|_2^2}-1\bigg|&\le\frac1{n\tau}, &
\bigg|\frac{\mathcal{I}_n(\tau,f)}{n^3\|f'\|_2^2/3}-1\bigg|&\le\frac{6}{n\tau},\label{eq:information_1}\\
\bigg|\frac{\mathcal{I}_n^{(n/2)}(\tau,f)}{n^3\|f'\|_2^2/12}-1\bigg|&\le\frac{9}{n\tau}, &
\bigg|\frac{\tilde{\mathcal{I}}_n(\tau,f)}{n^3\|f'\|_2^2/12}-1\bigg|&\le\frac{4}{n\tau},\label{eq:information_2}
\end{align}
and
\begin{equation}\label{eq:information_3}
\bigg|\int_0^n\Big(t-\frac n2\Big)f'(\tau t)g(\tau t)\,dt\bigg|\le\frac{n}{2\tau}\|f'\|_2\|g\|_2.
\end{equation}
\end{proposition}

The first bound in~\eqref{eq:information_1} does not involve $f$, and its proof shows that it holds for every shape $g$ with $\|g\|_2>0$ and all $n,\tau>0$.

Three consequences of Proposition~\ref{prop:information} are used in the rest of the paper.

First, the frequency and the shape carry information of different orders.
By~\eqref{eq:information_1} and~\eqref{eq:information_2}, for $n\tau\ge2$ and every shape $g$ with $\|g\|_2>0$,
\begin{equation}\label{eq:information_orders}
\tilde{\mathcal{I}}_n(\tau,f)=\frac{n^3}{12}\|f'\|_2^2\Big\{1+O\Big(\frac1{n\tau}\Big)\Big\},
\qquad
\int_0^n\big\{g(\tau t)\big\}^2\,dt=n\|g\|_2^2\Big\{1+O\Big(\frac1{n\tau}\Big)\Big\},
\end{equation}
where the second integral is the information about the shape in the direction $g$, and the constants in the $O$-terms are absolute.
This gives the rate $n^{-3/2}$ for the frequency and the local parameters $\tau=\tau_0+sn^{-3/2}$ and $f=f_0+gn^{-1/2}$ of Section~\ref{sec:introduction}; estimators of the frequency that attain the rate $n^{-3/2}$ when the shape is unknown are given in~\cite{Golubev:1988}.

Second, the frequency and the shape are asymptotically orthogonal in the parametrisation $\psi_{n/2}$, with the origin of time in the middle of the observation window.
For $n\tau\ge18$ and every shape $g$ with $\|g\|_2>0$,
\begin{equation}\label{eq:orthogonality}
\frac{\big|\int_0^n(t-n/2)f'(\tau t)g(\tau t)\,dt\big|}{\big\{\mathcal{I}_n^{(n/2)}(\tau,f)\big\}^{1/2}\big\{\int_0^n\{g(\tau t)\}^2\,dt\big\}^{1/2}}\le\frac3{n\tau},
\end{equation}
i.e., the cosine of the angle in $L^2\big([0,n]\big)$ between the score $(t-n/2)f'(\tau t)$ in $\psi_{n/2}$ and any direction $g(\tau\,\cdot)$ in the shape is at most $3/(n\tau)$.
Indeed, by the first bounds in~\eqref{eq:information_1} and~\eqref{eq:information_2}, the denominator is at least
\[
n^2\|f'\|_2\|g\|_2\Big\{\frac{(1-9/18)(1-1/18)}{12}\Big\}^{1/2}\ge\frac{n^2\|f'\|_2\|g\|_2}{5.05},
\]
and by~\eqref{eq:information_3} the numerator is at most $n\|f'\|_2\|g\|_2/(2\tau)$.
Accordingly, $\mathcal{I}_n^{(n/2)}(\tau,f)$ and $\tilde{\mathcal{I}}_n(\tau,f)$ agree to first order, by~\eqref{eq:information_2}.
With the origin of time at the start of the observation window, by contrast, the efficient information is asymptotically a quarter of the Fisher information $\mathcal{I}_n(\tau,f)$.

Third, the bounds of Proposition~\ref{prop:information} and~\eqref{eq:orthogonality} are uniform: they depend on $f$ and $g$ only through $\|f'\|_2$ and $\|g\|_2$, which are unchanged by time shifts, and their constants do not depend on $f$, $g$ or $\tau$.
This is relevant in Section~\ref{sec:BVM}, where we use the parametrisation $\psi_{n/2}$, in which the true shape $f_0(\cdot+\tau_0n/2)$ depends on $n$.

\begin{remark}[Least favourable direction]\label{rem:least_favourable}
The following justifies the description of the least favourable direction in Section~\ref{sec:introduction}.
If $n\tau$ is an integer, then, by the proof of Proposition~\ref{prop:information}, the infimum in the definition of $\tilde{\mathcal{I}}_n(\tau,f)$ is attained at the shape $g^*$ given by
\[
g^*(u)=\frac n2f'(u)+\frac{u-1/2}\tau f'(u),\quad u\in[0,1).
\]
This $g^*$ equals $(n/2)f'$ up to a term of norm at most $\|f'\|_2/(2\tau)$, and $f'$ is the derivative at $c=0$ of the time shift $c\mapsto f(\cdot+c)$: to first order, the shape compensates a change of frequency by a time shift.
In the local parameters $\tau=\tau_0+sn^{-3/2}$ and $f=f_0+gn^{-1/2}$, the submodel
\[
\tau\mapsto\Big(\tau,f_0-(\tau-\tau_0)\frac n2f_0'\Big)
\]
corresponds to $g=-sf_0'/2$, the least favourable direction.
\end{remark}

\section{Contraction of the posterior}\label{sec:estimation}

We adopt a Bayesian non-parametric approach to estimate the unknown parameter $\theta_0=(\tau_0,f_0)$, and later $\tau_0$, from data $X^{(n)}$ sampled from $P_{\theta_0}^{(n)}$.
Our prior is induced by the parametrisation $\psi$: we endow $\tau$ and $f$ with independent priors $\Pi_\tau$ and $\Pi_f$, and the prior $\Pi$ is the distribution of the induced periodic function $\psi(\tau,f)$; for simplicity we write $\Pi=\Pi_\tau\times\Pi_f$.
The prior $\Pi_\tau$ is a probability measure on $(0,\infty)$.

The prior $\Pi_f$ must be a prior on shapes, and is specified in one of the following ways.
\begin{enumerate}[label=\alph*)]
\item Sample $K\in\mathbb{N}$ from a certain distribution and then, given $K=k$, sample a vector $\bm{Z}=(Z_{-k},\dots,Z_k)\in\mathbb{R}^{2k+1}$.
\item Set $K=k_n$ for some deterministic, non-decreasing sequence $k_n\in\mathbb{N}$ and sample $\bm{Z}=(Z_{-k_n},\dots,Z_{k_n})\in\mathbb{R}^{2k_n+1}$.
\end{enumerate}
Either way, the prior $\Pi_f$ is the distribution of the trigonometric polynomial
\begin{equation}\label{eq:prior_structure}
f=\sum_{i=-K}^KZ_i\varphi_i=\bm{Z}^T\bm{\varphi}.
\end{equation}
We denote by $\Pi_K$ the prior on the dimension $K$, which in case b) is the point mass at $k_n$, and by $\Pi_{\bm{Z}}(\cdot\mid K=k)$ the conditional prior on the coordinates $\bm{Z}$, given $K=k$, so that $\Pi_f$ is the mixture of the conditional priors $\Pi_{\bm Z}(\cdot\mid K=k)$ with weights $\Pi_K(K=k)$.
If we independently sample $\tau$ from $\Pi_\tau$, then $\Pi$ is the law of the random periodic function
\begin{equation}\label{eq:prior_structure_2}
\theta(t)=\sum_{i=-K}^KZ_i\varphi_i(\tau t)=\bm{Z}^T\bm{\varphi}(\tau t),\quad t\in[0,n].
\end{equation}
The prior is thus a prior on the parameters $(\tau,K,\bm Z)$.
We write sets of parameters and integrals over them in terms of $(\tau,f)$, with $f=\bm Z^T\bm\varphi$, and we still write $\theta=(\tau,f)$ for the signal $\psi(\tau,f)$, also when $f\notin\mathcal{F}$.
The prior need not be concentrated on $(0,\infty)\times\mathcal{F}$, since it may put mass on, or arbitrarily close to, shapes that are constant or have a shorter period, and for such shapes $\tau$ is not the frequency of $\psi(\tau,f)$, by Lemma~\ref{lemma:periods}(i).
This does not affect the posterior, which depends on $(\tau,f)$ only through the signal $\psi(\tau,f)$, defined for every shape.
In Section~\ref{sec:estimation:semiparametric_consistency} we show that the posterior concentrates on shapes in $\mathcal{F}$, for which $\tau$ is the frequency.

For $c\in\mathbb{R}$, write $\bm Z^{(c)}$ for the coordinates of the time shift $f(\cdot+c)$, obtained from $\bm Z$ by the rotations of Section~\ref{sec:preliminaries}.
We call the prior $\Pi_f$ \emph{stationary} if, for every $c\in\mathbb{R}$, $(K,\bm Z^{(c)})$ has the same distribution as $(K,\bm Z)$; then $f(\cdot+c)$ has the same distribution as $f$.
This holds, for instance, if, given $K$, the coordinate $Z_0$ and the pairs $(Z_k,Z_{-k})$, $k=1,\dots,K$, are independent and the distribution of each pair is unchanged by rotations.
Stationary priors are used in Section~\ref{sec:BVM}.

Bayes' formula requires the likelihood to be jointly measurable in the parameters and the data, with the stochastic integral in~\eqref{eq:likelihood} defined for all signals at once.
The signals $\theta=\bm z^T\bm\varphi(\tau\,\cdot)$ in~\eqref{eq:prior_structure_2} are continuously differentiable, and for them we use
\[
\int_0^n\theta(t)\,dY_t=\theta(n)Y_n-\theta(0)Y_0-\int_0^n\theta'(t)Y_t\,dt.
\]
For each fixed $\theta$, the right side equals the stochastic integral almost surely, by the integration by parts formula for continuous semimartingales~\cite[Problem~3.3.12]{Karatzas:1998}, since $\theta$ is deterministic and of bounded variation on $[0,n]$.
The right side is a continuous function of $(\tau,\bm z,X^{(n)})\in(0,\infty)\times\mathbb{R}^{2k+1}\times C\big([0,n]\big)$, with the supremum norm on $C\big([0,n]\big)$, for every $k$, so the likelihood is jointly measurable.

For every measurable set $S$ of parameters, Bayes' formula delivers the posterior distribution
\begin{equation}\label{eq:bayess_formula}
\Pi\big(S\mid X^{(n)}\big)=\frac
{\int_S p_{\psi(\tau,f)}^{(n)}\big(X^{(n)}\big)\,d\Pi(\tau,f)}
{\int p_{\psi(\tau,f)}^{(n)}\big(X^{(n)}\big)\,d\Pi(\tau,f)};
\end{equation}
the denominator is positive, and it is finite almost surely, since, by Fubini's theorem, its integral with respect to $\nu^{(n)}$ is 1.
From the posterior~\eqref{eq:bayess_formula} one can make non-parametric inference on the signal $\theta$.
The posterior is later seen as a joint posterior on frequencies and shapes, and can be used to make semi-parametric inference on the frequency, for instance through point estimates and credible intervals for $\tau$.
In the next sections we study the asymptotic behaviour of the posterior~\eqref{eq:bayess_formula} under increasingly stronger assumptions on the components of the prior $\Pi$.

\subsection{Contraction of the non-parametric posterior}\label{sec:estimation:nonparametric_consistency}

Theorem~\ref{theo:checking_conditions} below states that under certain conditions on the priors on the components $K$, $\bm{Z}$ and $\tau$ of~\eqref{eq:prior_structure_2}, the posterior concentrates on the sieve
\begin{equation}\label{eq:sieve_adaptive}
\Theta_n=\Big\{(\tau,K,\bm{Z}):K\le\bar k_n,\;\bm{Z}\in w_1(m_n),\;|Z_0|\le m_n,\;\tau\in[l_n^{-1},u_n]\Big\},
\end{equation}
for sequences $\bar k_n$, $l_n$, $u_n$ and $m_n$ specified in the theorem, and contracts around the signal $\theta_0$ in $\|\cdot\|_{2,n}$.
The proof of Theorem~\ref{theo:checking_conditions} is in Appendix~\ref{sec:technical:estimation}.

\begin{theorem}[Non-parametric contraction]\label{theo:checking_conditions}
Assume that $f_0\in w_\alpha(M)$, for some $\alpha\ge1$ and $M>0$, that $\tau_0>0$, and that $\theta$ is endowed with a prior which is the law of~\eqref{eq:prior_structure_2}.
Assume that $\Pi_\tau$ admits a Lebesgue density which is bounded away from zero on a neighbourhood of $\tau_0$.
Let $k_n\le\bar k_n$ be integers, and let $l_n,u_n,m_n\ge1$ be bounded by some polynomial in $n$.
Let $\epsilon_n>0$ be such that $\epsilon_n\to0$ and $n\epsilon_n^2\to\infty$, as $n\to\infty$.
Assume that
\begin{equation}\label{eq:cond}
l_n=o(n),\qquad
\bar k_n\log n\lesssim n\epsilon_n^2,\qquad\text{and}\qquad
6C_\alpha\le(2k_n+1)^\alpha\epsilon_n,
\end{equation}
where $C_\alpha=M/\pi^\alpha$, and that
\begin{align}
\label{eq:cond_K_1}
\Pi_K(K=k_n)&\gtrsim\exp\big(-n\epsilon_n^2/4\big),\\
\label{eq:cond_K_2}
\Pi_K(K>\bar k_n)&\lesssim\exp\big(-5n\epsilon_n^2\big),\\
\label{eq:cond_Z_1}
\Pi_{\bm{Z}}\big(\|\bm{Z}-\bm{z}_n\|_2\le\epsilon_n/2\mid K=k_n\big)&\gtrsim\exp\big(-n\epsilon_n^2/4\big),\\
\label{eq:cond_Z_0}
\Pi_{\bm{Z}}\big(|Z_0|>m_n\mid K=k\big)&\lesssim\exp\big(-5n\epsilon_n^2\big),\quad k=1,\dots,\bar k_n,\\
\label{eq:cond_Z_3}
\Pi_{\bm{Z}}\big(\bm{Z}\notin w_1(m_n)\mid K=k\big)&\lesssim\exp\big(-5n\epsilon_n^2\big),\quad k=1,\dots,\bar k_n,\\
\label{eq:cond_tau_2}
\Pi_\tau\big(\tau\notin[l_n^{-1},u_n]\big)&\lesssim\exp\big(-5n\epsilon_n^2\big),
\end{align}
where $\bm{z}_n=(x_{0,-k_n},\dots,x_{0,k_n})$ are the coordinates of the truncation $f_{0,[k_n]}$ of $f_0$, and the constants in~\eqref{eq:cond_Z_0} and~\eqref{eq:cond_Z_3} do not depend on $k$.
Then, for all large enough $L$,
\begin{equation}\label{eq:contraction_statement}
\Pi\big(\|\theta-\theta_0\|_{2,n}>L\epsilon_n\mid X^{(n)}\big)\to0
\quad\text{and}\quad
\Pi\big(\Theta_n\mid X^{(n)}\big)\to1,
\end{equation}
both in $P_{\theta_0}^{(n)}$-probability, as $n\to\infty$.
\end{theorem}

The roles of the sequences are as follows.
The integer $k_n$ is the dimension at which $f_0$ is approximated, and the last condition in~\eqref{eq:cond} gives $\|f_0-f_{0,[k_n]}\|_2\le\epsilon_n/6$, by~\eqref{eq:truncation}.
The sets $\Theta_n$ have dimension at most $2\bar k_n+2$, and the second condition in~\eqref{eq:cond} bounds their entropy.
The bound $m_n$ makes the signals Lipschitz in $\tau$ (Lemma~\ref{lemma:upper_bound_metrics} in Appendix~\ref{sec:technical:estimation}), and $\tau\ge l_n^{-1}$ with $l_n=o(n)$ ensures that the signals in $\Theta_n$ have at least $n/l_n\to\infty$ periods in $[0,n]$.

The conclusion of Theorem~\ref{theo:checking_conditions} is in a sense the best one could hope for in terms of non-parametric contraction for periodic signals.
For instance, if $k_n\asymp\bar k_n\asymp(n/\log n)^{1/(2\alpha+1)}$, then~\eqref{eq:cond} allows
\[
\epsilon_n\asymp\Big(\frac n{\log n}\Big)^{-\alpha/(2\alpha+1)},
\]
the minimax rate for $\alpha$-smooth shapes up to a logarithmic factor (see Section~\ref{sec:examples}).
For sequences satisfying~\eqref{eq:cond}, conditions~\eqref{eq:cond_K_1} to~\eqref{eq:cond_Z_3} are the usual sufficient conditions for adaptive non-parametric contraction with random series priors (cf.~\cite{Arbel:2013,Shen:2015}), so the added uncertainty introduced by not knowing the frequency can be dealt with without extra assumptions on the prior on the shape.
In the proof, the ball of radius $\epsilon_n$ around $\theta_0$ contains the signals $(\tau,f)$ with $f$ close to $f_{0,[k_n]}$ and $|\tau-\tau_0|$ of order $\epsilon_n/(nk_n^{1/2})$; the prior mass of these frequencies is at least a negative power of $n$, and its logarithm is absorbed by $n\epsilon_n^2$.

The sequences in Theorem~\ref{theo:checking_conditions} depend on $\alpha$, but the prior need not.
If $K$ is random, conditions~\eqref{eq:cond_K_1} to~\eqref{eq:cond_tau_2} can hold for the sequences of every $\alpha\ge1$ at once, and the posterior then contracts at the rate $\epsilon_n$ that corresponds to the smoothness of $f_0$, without knowledge of $\alpha$, so the result is adaptive.
The same holds for Theorem~\ref{theo:checking_conditions_2} below when $h_n$ does not depend on $\alpha$, as for the prior of Section~\ref{sec:examples}.
If $K=k_n$, then $\bar k_n\ge k_n$, and the last two conditions in~\eqref{eq:cond} give
\[
\epsilon_n^{-1/\alpha}\lesssim k_n\lesssim\frac{n\epsilon_n^2}{\log n},
\]
so the rate in the previous paragraph requires $k_n\asymp(n/\log n)^{1/(2\alpha+1)}$, which depends on $\alpha$, and the result is not adaptive.
Priors of both kinds are given in Section~\ref{sec:examples}.

The requirements on the prior on the frequency are mild.
Condition~\eqref{eq:cond_tau_2} holds, for instance, if $\Pi_\tau$ has a density bounded by a multiple of $\exp(-b/\tau-b\tau)$ for some $b>0$, with $l_n=u_n=5n\epsilon_n^2/b$, since then
\[
\Pi_\tau\big(\tau<l^{-1}\big)\lesssim e^{-bl},\qquad\Pi_\tau(\tau>u)\lesssim e^{-bu};
\]
moreover, $l_n=o(n)$, since $\epsilon_n\to0$.

Theorem~\ref{theo:checking_conditions} is, however, not adequate to study the marginal posterior for the frequency: the frequency of a signal sampled from the posterior is not necessarily close to $\tau_0$, even if the signal is close to $\theta_0$ in $\|\cdot\|_{2,n}$.
Indeed, if $f_0$ is a trigonometric polynomial, $l_n,u_n,m_n\to\infty$ and $\eta\le\epsilon_n/2$, then the signal $(\tau_0/2,f_\eta)$ of Section~\ref{sec:statement_of_problem:parametrisation}, whose frequency is $\tau_0/2$, has parameters in $\Theta_n$ for large $n$, and lies in the ball of radius $\epsilon_n$ around $\theta_0$, by~\eqref{eq:f_eta_close}.
To address this, we consider a smaller sieve and impose stronger conditions on the prior.

\subsection{Contraction of the marginal posteriors}\label{sec:estimation:semiparametric_consistency}

We work with the smaller sieve
\begin{equation}\label{eq:sieve_good_one}
\Xi_n=\Xi_n(\bar k_n,l_n,u_n,m_n,h_n)=\big\{(\tau,K,\bm Z)\in\Theta_n:\bm Z\in\Xi(h_n)\big\},
\end{equation}
for a sequence $h_n\in(0,1]$.
By Lemma~\ref{lemma:xi_sets}(ii), shapes in $\Xi(h_n)$ belong to $\mathcal{F}$, so $\tau$ is the frequency of the signal $(\tau,f)$ for every parameter in $\Xi_n$.
Proposition~\ref{prop:smaller_sieve} below states that signals in $\Xi_n$ that approximate $\theta_0$ well must have frequency close to $\tau_0$ and shape close to $f_0$.
It is a deterministic statement, which uses only that $\tau\ge l_n^{-1}$, that $f$ is a trigonometric polynomial of degree at most $\bar k_n$, and that $f\in\Xi(h_n)$; it imposes no condition on the ratio $\tau/\tau_0$.
Its constants depend on $f_0$ through $\sigma_0$, through the number $H$ of Section~\ref{sec:statement_of_problem:parametrisation}, and through the smallest $J\in\mathbb{N}$ such that $\sum_{|i|>J}x_{0,i}^2\le H\sigma_0^2/2$.
The proof is in Appendix~\ref{sec:technical:estimation}.

\begin{proposition}[Localisation]\label{prop:smaller_sieve}
Let $\tau_0>0$ and $f_0\in\mathcal{F}$, and let $\sigma_0$, $H$ and $J$ be as above.
Let $n,\bar k\in\mathbb{N}$, $l>0$, $h\in(0,1]$ and $\epsilon>0$ be such that
\begin{equation}\label{eq:new_conditions}
n\tau_0\ge16\bar k,\qquad
n\ge4l,\qquad
2\bar k\ge J,\qquad\text{and}\qquad
400\,\epsilon^2\le\min\{h,H\}\,\sigma_0^2.
\end{equation}
Let $(\tau,f)$ be such that $\tau\ge l^{-1}$, $f$ is a trigonometric polynomial of degree at most $\bar k$ with coordinates in $\Xi(h)$, and $\|\psi(\tau,f)-\theta_0\|_{2,n}\le\epsilon$.
Then
\begin{equation}\label{eq:localisation_tau}
|\tau-\tau_0|\le\frac{56\,\epsilon}{\sigma_0n},
\end{equation}
and, for an absolute constant $C$,
\begin{equation}\label{eq:localisation_f}
\|f-f_0\|_2^2\le C\bigg\{\epsilon^2+\|f_0-f_{0,[\bar k]}\|_2^2+\frac{(1+\log\bar k)\big(\|f_0\|_2^2+\epsilon^2\big)}{n\tau_0}\bigg\}.
\end{equation}
\end{proposition}

By Proposition~\ref{prop:smaller_sieve}, marginal contraction follows if the posterior also concentrates on the sieve $\Xi_n$, and the following theorem shows that it indeed does.

\begin{theorem}[Marginal contraction]\label{theo:checking_conditions_2}
Assume the conditions of Theorem~\ref{theo:checking_conditions}, and let $h_n\in(0,1]$ be such that
\begin{align}
\label{eq:cond_h}
\epsilon_n^2&=o(h_n),\\
\label{eq:cond_Z_2}
\sum_{k=1}^{\bar k_n}\Pi_K(K=k)\,\Pi_{\bm{Z}}\big(\bm{Z}\notin\Xi(h_n)\mid K=k\big)&\lesssim\exp\big(-5n\epsilon_n^2\big).
\end{align}
Then $\Pi\big(\Xi_n\mid X^{(n)}\big)\to1$ and, for all large enough $L$,
\[
\Pi\big(|\tau-\tau_0|>L\epsilon_n/n\mid X^{(n)}\big)\to0
\quad\text{and}\quad
\Pi\big(\|f-f_0\|_2>L\epsilon_n\mid X^{(n)}\big)\to0,
\]
all in $P_{\theta_0}^{(n)}$-probability, as $n\to\infty$.
\end{theorem}

In particular, the marginal posterior for the frequency is consistent.
For a prior that does not depend on $h_n$, condition~\eqref{eq:cond_Z_2} can be satisfied by letting $h_n$ converge to zero fast enough, but~\eqref{eq:cond_h} imposes a lower bound on $h_n$.
Together they require the prior to put all but an exponentially small mass on shapes whose distance to the shapes with a shorter period, which is at least $h_n^{1/2}\|f-\bar f\|_2$ by Lemma~\ref{lemma:xi_sets}(iv), is of larger order than $\epsilon_n$; priors of this type are given in Section~\ref{sec:examples}.
The conditions do not involve $H$, $\sigma_0$ or $J$.
Since~\eqref{eq:cond_Z_2} puts almost all prior mass in $\Xi(h_n)$, condition~\eqref{eq:cond_Z_1} requires enough prior mass near $f_{0,[k_n]}$ inside $\Xi(h_n)$.

The rate $\epsilon_n/n$ for the frequency is not optimal, since it is slower than $n^{-3/2}$ by the factor $\sqrt n\,\epsilon_n\to\infty$.
It is, however, the best rate that can be deduced from contraction in $\|\cdot\|_{2,n}$ at rate $\epsilon_n$ alone.
Indeed, let $\delta>0$, $n\tau_0\ge2$ and $f_0\in\mathcal{W}_1$.
Writing $\psi(\tau_0+\delta,f_0)-\theta_0$ as the integral of its derivative with respect to the frequency, as in the proof of Proposition~\ref{prop:information}, Jensen's inequality and~\eqref{eq:information_1} give
\[
\|\psi(\tau_0+\delta,f_0)-\theta_0\|_{2,n}^2\le\frac{\delta^2}{n}\int_0^1\mathcal{I}_n(\tau_0+r\delta,f_0)\,dr\le\frac{n^2\delta^2}{3}\|f_0'\|_2^2\Big(1+\frac6{n\tau_0}\Big).
\]
Hence the ball of radius $\epsilon_n$ around $\theta_0$ contains all signals $\psi(\tau,f_0)$ with $|\tau-\tau_0|\le\sqrt3\,\epsilon_n/(2n\|f_0'\|_2)$, since $n\tau_0\ge2$.
Contraction at rate $\epsilon_n$ is therefore compatible with a posterior that puts all its mass on frequencies at distance of order $\epsilon_n/n$ from $\tau_0$, so no rate faster than $\epsilon_n/n$ can be deduced from it, and Proposition~\ref{prop:smaller_sieve} is sharp in this sense.

The marginal posterior for the frequency may still concentrate faster.
In Section~\ref{sec:BVM} we take a different approach and study the asymptotic distribution of the marginal posterior for $\tau$, using Theorem~\ref{theo:checking_conditions_2} as a first localisation.
Before that, we give in Section~\ref{sec:examples} examples of priors and sequences for which the conditions of our results so far are satisfied.

\section{Examples of priors}\label{sec:examples}

In this section we give priors for which the conditions of Theorems~\ref{theo:checking_conditions} and~\ref{theo:checking_conditions_2} hold, adaptively to the smoothness of the shape.
Throughout, the prior $\Pi_\tau$ on the frequency has a continuous and positive Lebesgue density on $(0,\infty)$ which is bounded by a multiple of $\exp(-b/\tau-b\tau)$ for some $b>0$, for instance the density proportional to $\exp(-b/\tau-b\tau)$.
By the discussion after Theorem~\ref{theo:checking_conditions}, $\Pi_\tau$ then satisfies the conditions of that theorem with $l_n=u_n=5n\epsilon_n^2/b$.
The proofs of the results of this section are in Appendix~\ref{sec:technical:examples}.

Let $h_n\in(0,1/4]$.
Given $K=k$, let $\bm Z=(Z_{-k},\dots,Z_k)$ be distributed as a vector of independent standard normal variables conditioned on the event $\{\bm Z\in\Xi(h_n)\}$.
This event has probability at least an absolute constant $c_*>0$, for every $k$ (Lemma~\ref{lemma:xi_mass} in Appendix~\ref{sec:technical:examples}).
Since the event does not involve $Z_0$, the coordinate $Z_0$ is standard normal and independent of the other coordinates.

The unconditioned prior is stationary, by the example after the definition in Section~\ref{sec:estimation}, and $\Xi(h_n)$ is unchanged by time shifts (Section~\ref{sec:statement_of_problem:parametrisation}), so the conditioned prior is again stationary.

For the dimension $K$ we consider the two cases of Section~\ref{sec:estimation}:
\begin{enumerate}[label=\alph*)]
\item $\Pi_K$ satisfies, for some constants $b_1,b_2>0$ and all large enough $k$,
\begin{equation}\label{eq:prior_K_example}
\exp(-b_1k\log k)\le\Pi_K(K=k)\le\exp(-b_2k\log k);
\end{equation}
this holds, for instance, if $K-1$ has a Poisson distribution with any mean, by the bounds $(k/e)^k\le k!\le k^k$;
\item $K=k_n$ for a deterministic sequence $k_n$.
\end{enumerate}

\begin{proposition}\label{prop:prior_example}
Let $\tau_0>0$ and $f_0\in\mathcal{F}\cap w_\alpha(M)$, for some $\alpha\ge1$ and $M>0$.
Let
\[
\tilde k_n=\Big(\frac n{\log n}\Big)^{1/(2\alpha+1)},\qquad k_n=\lceil\tilde k_n\rceil,\qquad\epsilon_n=C_0\tilde k_n^{-\alpha},
\]
for a constant $C_0>0$.
Let the prior be as above, with $K$ as in case a) or b), and with $h_n\in(0,1/4]$ such that $h_n\to0$ and $\epsilon_n^2=o(h_n)$.
Then, for all large enough $C_0$, the conditions of Theorems~\ref{theo:checking_conditions} and~\ref{theo:checking_conditions_2} hold, with
\[
l_n=u_n=\frac{5n\epsilon_n^2}b,\qquad m_n=10\pi\bar k_n\big(n\epsilon_n^2\big)^{1/2},
\]
$\bar k_n=\lceil D\tilde k_n\rceil$ for a large enough constant $D$ in case a), and $\bar k_n=k_n$ in case b).
In particular, for all large enough $L$,
\begin{align*}
&\Pi\big(\|\theta-\theta_0\|_{2,n}>L\epsilon_n\mid X^{(n)}\big)\to0,\qquad
\Pi\big(|\tau-\tau_0|>L\epsilon_n/n\mid X^{(n)}\big)\to0,\\
&\Pi\big(\|f-f_0\|_2>L\epsilon_n\mid X^{(n)}\big)\to0,
\end{align*}
in $P_{\theta_0}^{(n)}$-probability, as $n\to\infty$.
\end{proposition}

Condition~\eqref{eq:cond_h} holds for $h_n=1/\max\{4,\log n\}$ and every $\alpha\ge1$, since $\epsilon_n^2$ is a negative power of $n$ times a power of $\log n$.
In case a), the prior with this $h_n$ does not depend on $\alpha$, $M$, $\tau_0$ or $f_0$, and Proposition~\ref{prop:prior_example} gives contraction at the rates stated there for every $\alpha\ge1$, $M>0$ and $f_0\in\mathcal{F}\cap w_\alpha(M)$.
The posterior for the signal and the marginal posteriors for the frequency and the shape are therefore adaptive, also to the number $H$ of Section~\ref{sec:statement_of_problem:parametrisation}.
The logarithmic factor comes from the entropy condition in~\eqref{eq:cond} and from~\eqref{eq:cond_Z_1}.

The condition $h_n\to0$ is used only for~\eqref{eq:cond_Z_1}, through $\limsup_nh_n<H$: the shapes close to $f_{0,[k_n]}$ are then in $\Xi(h_n)$, since $f_0\in\Xi(H)$.
The proof shows that a fixed $h_n=h\in(0,1/4]$ can be used when $f_0\in\Xi(H)$ for some $H>h$.
The verification uses the normal distribution only through a lower bound on its density, tail bounds for $Z_0$ and for sums of squares, and the lower bound on the probability of $\Xi(1/4)$ in Lemma~\ref{lemma:xi_mass}, so it carries over to other distributions with these properties.
Distributions of the pairs $(Z_k,Z_{-k})$ that are unchanged by rotations keep the prior stationary, by Section~\ref{sec:estimation}.

Up to the logarithmic factor, the rates for the signal and for the shape are the minimax rates over $w_\alpha(M)$ when the frequency is known, and a fortiori when it is not (Lemma~\ref{lemma:minimax} in Appendix~\ref{sec:technical:examples}).
The rate for the frequency is not, as discussed in Section~\ref{sec:estimation:semiparametric_consistency}.

Without the conditioning, independent standard normal coordinates satisfy the conditions of Theorem~\ref{theo:checking_conditions}, by the proof of Proposition~\ref{prop:prior_example}.
When $f_0$ has a non-zero harmonic of index at least 2, however, Gaussian random series priors with independent centred coordinates fail condition~\eqref{eq:cond_Z_2} of Theorem~\ref{theo:checking_conditions_2}, for any choice of the variances (Proposition~\ref{prop:gaussian_fails} in Appendix~\ref{sec:technical:examples}); this motivates the conditioning above.

\section{Asymptotic distribution of the posterior}\label{sec:BVM}

By Theorem~\ref{theo:checking_conditions_2}, the marginal posterior for the frequency concentrates on frequencies at distance at most $L\epsilon_n/n$ from $\tau_0$.
In this section we study it at the scale $n^{-3/2}$, through the local parameter $s=n^{3/2}(\tau-\tau_0)$.
We first show that, for stationary priors on the shape, the origin of time can be moved to the middle of the observation window without changing the posterior.
In the corresponding parametrisation we then give a LAN expansion, which identifies the efficient information and the efficient estimators of the frequency, and the Bernstein-von Mises theorem.
The proofs are in Appendix~\ref{sec:technical:BVM}.

\subsection{Stationary priors and the origin of time}\label{sec:BVM:stationary}

Recall from Section~\ref{sec:statement_of_problem:information} the parametrisation $\psi_c$, with the origin of time at $c\in\mathbb{R}$, in which the parameter $(\tau,f)$ of $\psi$ becomes $\big(\tau,f(\cdot+\tau c)\big)$.
Recall also the definition of a stationary prior from Section~\ref{sec:estimation}; the prior of Section~\ref{sec:examples} is an example.
For a parameter $(\tau,K,\bm Z)$, we write $\big(\tau,f(\cdot+\tau c)\big)$ for the parameter $(\tau,K,\bm Z^{(\tau c)})$, whose shape is $f(\cdot+\tau c)$.

\begin{proposition}[Origin of time]\label{prop:origin}
Let $\Pi=\Pi_\tau\times\Pi_f$, with $\Pi_f$ stationary.
Then, for every $c\in\mathbb{R}$, the map $(\tau,f)\mapsto\big(\tau,f(\cdot+\tau c)\big)$ leaves $\Pi$ invariant, and, for every measurable set $S$ of parameters,
\[
\Pi\Big(\big(\tau,f(\cdot+\tau c)\big)\in S\;\Big|\;X^{(n)}\Big)=\frac
{\int_S p_{\psi_c(\tau,f)}^{(n)}\big(X^{(n)}\big)\,d\Pi(\tau,f)}
{\int p_{\psi_c(\tau,f)}^{(n)}\big(X^{(n)}\big)\,d\Pi(\tau,f)}.
\]
In particular, the marginal posterior for the frequency is the same in the parametrisations $\psi$ and $\psi_c$, with the same prior $\Pi$.
\end{proposition}

The right side of the display in Proposition~\ref{prop:origin} is the posterior in the model parametrised by $\psi_c$, with prior $\Pi$.
We use $c=n/2$, for which the frequency and the shape are asymptotically orthogonal, by~\eqref{eq:orthogonality}.
In the parametrisation $\psi_{n/2}$ the shape of the parameter $(\tau,f)$ is $f(\cdot+\tau n/2)$, and the true shape is
\[
\tilde f_{0,n}=f_0(\cdot+\tau_0n/2),
\]
which depends on $n$; its amplitudes, its norm and the norm of its derivative are those of $f_0$, by Section~\ref{sec:preliminaries}.
Under the conditions of Theorem~\ref{theo:checking_conditions_2}, the marginal posterior for the shape in $\psi_{n/2}$ contracts around $\tilde f_{0,n}$ at the rate $\epsilon_n$ (see the proof of Theorem~\ref{theo:BvM}).

\subsection{Local asymptotic normality and efficiency}\label{sec:BVM:LAN}

Let $f_0\in\mathcal{W}_1$, and write
\begin{align*}
\mathcal{I}_0&=\frac{\|f_0'\|_2^2}{12},
\qquad
\Delta_n=\frac1{n^{3/2}}\int_0^n\Big(t-\frac n2\Big)f_0'(\tau_0t)\,dW_t,\\
\mathbb{W}_n(g)&=\frac1{n^{1/2}}\int_0^ng\Big(\tau_0\Big(t-\frac n2\Big)\Big)\,dW_t,
\end{align*}
for shapes $g$.
The variable $\Delta_n$ is the score for the frequency in the parametrisation $\psi_{n/2}$, scaled by $n^{-3/2}$, since
\[
\tilde f_{0,n}'\Big(\tau_0\Big(t-\frac n2\Big)\Big)=f_0'(\tau_0t),\quad t\in[0,n].
\]
It is normally distributed with mean zero and variance $\mathcal{I}_n^{(n/2)}(\tau_0,f_0)/n^3$, which converges to $\mathcal{I}_0$, by~\eqref{eq:information_2}.

\begin{proposition}[LAN]\label{prop:LAN}
Let $\tau_0>0$ and $f_0\in\mathcal{F}\cap\mathcal{W}_1$.
For $s\in\mathbb{R}$ and a shape $g$, let $\tau_n=\tau_0+sn^{-3/2}$ and $f_n=\tilde f_{0,n}+n^{-1/2}g$.
Then, under $P_{\theta_0}^{(n)}$,
\[
\log\frac{p_{\psi_{n/2}(\tau_n,f_n)}^{(n)}}{p_{\theta_0}^{(n)}}\big(X^{(n)}\big)=s\Delta_n+\mathbb{W}_n(g)-\frac12\big(s^2\mathcal{I}_0+\|g\|_2^2\big)+o_P(1),
\]
and $\big(\Delta_n,\mathbb{W}_n(g)\big)$ converges in distribution to a pair of independent normal variables with mean zero and variances $\mathcal{I}_0$ and $\|g\|_2^2$, as $n\to\infty$.
\end{proposition}

With the origin of time at $n/2$, the model is that of~\cite{Golubev:1988} and~\cite{Castillo:2007}, observed on a window centred at the origin, after the change from period to frequency; the efficient information for the frequency is computed in~\cite[proof of Theorem~2]{Castillo:2007}, and the case of a known shape, for twice continuously differentiable shapes, is in~\cite[Section~III.5]{Ibragimov:1981}.
Proposition~\ref{prop:LAN} gives the joint expansion in the frequency and the shape for $f_0\in\mathcal{W}_1$.
The expansion has no cross term between $s$ and $g$.
Hence, in the parametrisation $\psi_{n/2}$, not knowing the shape does not reduce the information about $s$, which is $\mathcal{I}_0$, the limit of $\tilde{\mathcal{I}}_n(\tau_0,f_0)/n^3$, by~\eqref{eq:information_2}.
With $g\equiv0$, Proposition~\ref{prop:LAN} states that the model in which the shape in $\psi_{n/2}$ is fixed at $\tilde f_{0,n}$, and only the frequency varies, is LAN at $\tau_0$ with normalising factor $(n^{3/2}\mathcal{I}_0^{1/2})^{-1}$, in the sense of~\cite[Definition~II.2.1]{Ibragimov:1981}.
Write $E_\tau^*$ for the expectation in this model, at the frequency $\tau$.
By~\cite[Theorem~II.12.1]{Ibragimov:1981}, for every sequence of estimators $T_n$ based on $X^{(n)}$, every loss function $w$ in the class $\mathbf{W}_{e,2}$ of~\cite[Section~I.2]{Ibragimov:1981}, which contains the functions $w(x)=|x|^m$, $m>0$, and every $\delta>0$,
\[
\liminf_{n\to\infty}\sup_{|\tau-\tau_0|<\delta}E_\tau^*\Big[w\Big(\mathcal{I}_0^{1/2}n^{3/2}(T_n-\tau)\Big)\Big]\ge\mathbb{E}\big[w(\xi)\big],
\]
where $\xi$ is a standard normal variable.
In the parametrisation $\psi$ the shapes of this model are the time shifts $f_0(\cdot+(\tau_0-\tau)n/2)$ of $f_0$, so the same bound holds for every model that contains them.
By the second part of~\cite[Theorem~II.12.1]{Ibragimov:1981}, for loss functions in this class whose sets $\{w<y\}$ are bounded for small $y>0$, an estimator $\hat\tau_n$ attains the bound, in the limit as $\delta\to0$, if and only if
\begin{equation}\label{eq:efficient_estimator}
n^{3/2}(\hat\tau_n-\tau_0)=\frac{\Delta_n}{\mathcal{I}_0}+o_P(1),
\end{equation}
under $P_{\theta_0}^{(n)}$.
This justifies the definition of efficient estimators in Section~\ref{sec:introduction}, which is~\eqref{eq:efficient_estimator} in the present notation; cf.~\cite[Section~1.2]{Castillo:2012}.

\subsection{The Bernstein-von Mises theorem}\label{sec:BVM:BvM}

\begin{theorem}[Bernstein-von Mises]\label{theo:BvM}
Assume the conditions of Theorems~\ref{theo:checking_conditions} and~\ref{theo:checking_conditions_2}, that the density of $\Pi_\tau$ is continuous and positive at $\tau_0$, that $\Pi_f$ is stationary, and that
\begin{equation}\label{eq:cond_BvM}
n^{1/2}\bar k_n\epsilon_n^2\to0.
\end{equation}
Then
\begin{equation}\label{eq:BvM_statement}
\sup_B\bigg|\Pi\Big(n^{3/2}(\tau-\tau_0)\in B\;\Big|\;X^{(n)}\Big)-N\Big(\frac{\Delta_n}{\mathcal{I}_0},\frac1{\mathcal{I}_0}\Big)(B)\bigg|\to0,
\end{equation}
in $P_{\theta_0}^{(n)}$-probability, as $n\to\infty$, where the supremum is over all Borel sets $B$.
Moreover, for $a\in(0,1)$, let $Q_n(a)$ be the $a$-quantile of the marginal posterior for the frequency.
Then, as $n\to\infty$,
\[
P_{\theta_0}^{(n)}\big(Q_n(a/2)\le\tau_0\le Q_n(1-a/2)\big)\to1-a.
\]
\end{theorem}

If $\hat\tau_n$ is efficient, then, by~\eqref{eq:efficient_estimator}, the normal distributions with variance $1/\mathcal{I}_0$ centred at $\Delta_n/\mathcal{I}_0$ and at $n^{3/2}(\hat\tau_n-\tau_0)$ are asymptotically equal in total variation, so~\eqref{eq:BvM_statement} is equivalent to~\eqref{eq:BvM_intro}.
By the proof, for every $a\in(0,1)$,
\[
n^{3/2}\big(Q_n(a)-\tau_0\big)=\frac{\Delta_n}{\mathcal{I}_0}+\frac{z_a}{\mathcal{I}_0^{1/2}}+o_P(1),
\]
where $z_a$ is the $a$-quantile of the standard normal distribution.
In particular, under the conditions of Theorem~\ref{theo:BvM}, the posterior median is an efficient estimator.
Since $\Delta_n$ is bounded in probability, \eqref{eq:BvM_statement} also shows that the marginal posterior for the frequency contracts at the rate $n^{-3/2}$, faster than the rate $\epsilon_n/n$ of Theorem~\ref{theo:checking_conditions_2}.
The last assertion of Theorem~\ref{theo:BvM} states that equal-tailed credible intervals for the frequency are asymptotically confidence intervals with the same level.
By the previous display, their length is $n^{-3/2}$ times $\mathcal{I}_0^{-1/2}(z_{1-a/2}-z_{a/2})+o_P(1)$, the length of the confidence interval based on an efficient estimator.

The proof of Theorem~\ref{theo:BvM} rests on Lemma~\ref{lemma:BvM_general} in Appendix~\ref{sec:technical:BVM}, a version of~\cite[Theorem~1]{Castillo:2012} for an arbitrary normalising rate.
It holds for every product prior whose marginal for the frequency has a density that is continuous and positive at $\tau_0$, and it requires the log-likelihood ratios at $(\tau,f)$ and $(\tau_0,f)$ to differ by $s\Delta_n-s^2\mathcal{I}_0/2$, up to an error that is small relative to $1+s^2$, uniformly over a product of a set of frequencies and a set of shapes on which the posterior concentrates.
This is condition~(N) of~\cite{Castillo:2012}, required there on a ball; the localisation on the product replaces condition~(C) there, whose part on the posterior with known parameter of interest is not needed.
The part of the proof that is specific to the model is the verification of this expansion in the parametrisation $\psi_{n/2}$, and it is there that the choice of the origin of time and condition~\eqref{eq:cond_BvM} enter.

Together with~\eqref{eq:cond}, condition~\eqref{eq:cond_BvM} can hold only if $\alpha>3/2$.
Indeed, by~\eqref{eq:cond}, $\bar k_n\ge k_n\gtrsim\epsilon_n^{-1/\alpha}$ and $n\epsilon_n^2\gtrsim\bar k_n$, so that $\epsilon_n\gtrsim n^{-\alpha/(2\alpha+1)}$ and
\[
n^{1/2}\bar k_n\epsilon_n^2\gtrsim n^{1/2}\epsilon_n^{2-1/\alpha}\gtrsim n^{(3-2\alpha)/(4\alpha+2)},
\]
which does not tend to zero if $\alpha\le3/2$.

\begin{corollary}\label{cor:BvM_example}
Let the prior and $f_0$ be as in Proposition~\ref{prop:prior_example}, with $\alpha>3/2$.
Then~\eqref{eq:BvM_statement} and the last assertion of Theorem~\ref{theo:BvM} hold.
\end{corollary}

In case a) of Section~\ref{sec:examples}, with $h_n=1/\max\{4,\log n\}$, the prior does not depend on $\alpha$, so one prior gives the Bernstein-von Mises theorem for every $\alpha>3/2$, $M>0$ and $f_0\in\mathcal{F}\cap w_\alpha(M)$.

\section{Discussion}\label{sec:discussion}

Our results are pointwise in the true parameter.
The constants in Proposition~\ref{prop:smaller_sieve} depend on $\theta_0$ through $\tau_0$, $\sigma_0$, $H$, $J$ and $\|f_0\|_2$, and the coverage of the credible intervals in Theorem~\ref{theo:BvM} holds for each fixed $\theta_0$.
By the example in Section~\ref{sec:statement_of_problem:parametrisation}, some restriction on the class of true shapes is unavoidable for uniform statements.
A natural class consists of the true parameters with $\tau_0$ in a compact subset of $(0,\infty)$ and shapes in $w_\alpha(M)\cap\Xi(H)$, for a fixed $H\in(0,1]$, with bounded mean and with $\sigma_0$ bounded below, on which the constants of Proposition~\ref{prop:smaller_sieve} are uniform.
We have not pursued uniform versions of the other results.

Theorem~\ref{theo:BvM} requires the smoothness $\alpha>3/2$, through condition~\eqref{eq:cond_BvM}, while the efficient information is finite as soon as $f_0\in\mathcal{W}_1$.
In the proof, the condition controls the remainder terms of the expansion, which involve the derivatives of the shapes, uniformly over the set on which the posterior concentrates.
The same threshold appears in~\cite[Section~2.1, Theorem~3]{Castillo:2012}, where it comes from the verification of the analogous expansion in the translation model.
We do not know whether a Bernstein-von Mises theorem for the frequency holds when $1\le\alpha\le3/2$.

Stationarity of the prior on the shape is used only in Proposition~\ref{prop:origin}, to move the origin of time to the middle of the observation window without changing the posterior.
For a prior that is not stationary, for instance one under which the phases of the harmonics are not uniformly distributed, the prior induced on the parameters of $\psi_{n/2}$ is in general not a product prior and depends on $n$, so Lemma~\ref{lemma:BvM_general} does not apply in $\psi_{n/2}$.
In the original parametrisation the prior would have to satisfy conditions in the least favourable direction, as discussed in Section~\ref{sec:introduction}, and we do not know whether the marginal posterior for the frequency is then still centred at an efficient estimator.

When the frequency is known, the minimax rate over $w_\alpha(M)$ has no logarithmic factor, and it is attained adaptively over $\alpha$ and $M$~\cite[Theorem~3.6]{Tsybakov:2009}.
For the Bernstein-von Mises theorem the logarithmic factor in the rates of Section~\ref{sec:examples} is immaterial, since~\eqref{eq:cond_BvM} holds for the priors of that section for every $\alpha>3/2$.

The contraction results use the Gaussian likelihood mainly through the existence of tests and bounds on Kullback-Leibler divergences, and the Bernstein-von Mises theorem through the expansion in Lemma~\ref{lemma:BvM_general}, whose verification uses the quadratic form of the log-likelihood ratio in the signal and a maximal inequality for Gaussian processes (Lemma~\ref{lemma:maximal}).
Our proof extends~\cite[Theorem~1]{Castillo:2012} in three ways: to an arbitrary normalising rate, to localisation on a product of a set of parameters of interest and a set of nuisance parameters, in place of the concentration of the posterior with known parameter of interest, and to models in which the loss of information disappears under a reparametrisation that leaves the prior unchanged.
The last is a general principle: Proposition~\ref{prop:origin} uses the model only through the fact that the likelihood depends on the frequency and the shape through the signal, so it holds, with the same proof, in every dominated model with a jointly measurable likelihood of this kind.
Similar arguments should therefore apply to other models with periodic features, such as non-parametric regression, or Poisson processes on $[0,n]$ with intensity $f(\tau t)$, for which non-parametric Bayesian estimation of the intensity is studied in~\cite{Belitser:2015,Kirichenko:2015}.
In the Poisson model, for a known shape $f$ that is smooth and bounded away from zero, the Fisher information about $\tau$ is
\[
\int_0^nt^2\frac{\big\{f'(\tau t)\big\}^2}{f(\tau t)}\,dt=\int_0^nt^2\Big\{\big(2\sqrt f\big)'(\tau t)\Big\}^2\,dt,
\]
which is~\eqref{eq:fisher_intro} with $2\sqrt f$ in place of $f$, so we expect the results to carry over with $2\sqrt f$ in the role of the shape, for stationary priors on $\sqrt f$ or on another transformation of the shape.

Model~\eqref{eq:gwn_intro} does not reduce to a fixed sequence model, since the natural basis $\varphi_i(\tau\,\cdot)$ for the signal depends on the unknown frequency.
In non-parametric regression with observations on a grid of spacing $\Delta$, the frequencies $\tau$ and $\tau+j/\Delta$, $j\in\mathbb{Z}$, give the same values on the grid, and so do $\tau$ and $1/\Delta-\tau$ with the shape reflected in time.
The frequency is then identifiable only in a known range below the Nyquist frequency $1/(2\Delta)$, or with irregular or random designs~\cite{Hall:2000,LevyLeduc:2008}.
With many observations per period we expect our results to carry over, although the equivalence of~\cite{Brown:1996} concerns a fixed observation window, while here the window grows with $n$.

\appendix

\section{Proofs for Section~\ref{sec:statement_of_problem}}\label{sec:technical:statement}

\paragraph{Proof of Lemma~\ref{lemma:periods}:}
For (i), let $A=\{k\in\mathbb{N}:x_k^2+x_{-k}^2>0\}$, which is non-empty since $f$ is not constant.
Since $f$ and its time shifts are 1-periodic, $c\in\mathbb{R}$ is a period of $f$ if, and only if, $f(\cdot+c)=f$ in $L^2\big([0,1]\big)$, i.e., if, and only if, $f(\cdot+c)$ and $f$ have the same coordinates.
By~\eqref{eq:shift}, this holds if, and only if, $\sin(\pi kc)=0$ for every $k\in A$, so $c$ is a period of $f$ if, and only if, $kc\in\mathbb{Z}$ for every $k\in A$.
If $c=j/d$ with $j\in\mathbb{Z}$, then $kc=jk/d\in\mathbb{Z}$ for every $k\in A$, since $d$ divides $k$.
Conversely, for $r\in\mathbb{N}$ let $d_r$ be the greatest common divisor of the $r$ smallest elements of $A$ (of all elements of $A$, if $A$ has fewer than $r$ elements).
The sequence $(d_r)$ is a non-increasing sequence of positive integers, so it is eventually equal to some $d'$, which divides every element of $A$; hence $d'$ divides $d$, and since $d$ divides every $d_r$, $d=d'$.
So $d$ is the greatest common divisor of finitely many elements $k_1,\dots,k_r$ of $A$, and, by B\'ezout's identity, $d=\sum_{l=1}^r m_lk_l$ for some integers $m_1,\dots,m_r$.
If $kc\in\mathbb{Z}$ for every $k\in A$, then $dc=\sum_{l=1}^r m_l(k_lc)\in\mathbb{Z}$.
This proves the first assertion of (i).
For the second, the variable change $u=\tau t$ shows that $c$ is a period of $\psi(\tau,f)$ if, and only if, $\tau c$ is a period of $f$, i.e., $\tau c\in\{j/d:j\in\mathbb{Z}\}$.
The smallest positive period of $\psi(\tau,f)$ is therefore $1/(d\tau)$, and its frequency is $d\tau$.

For (ii), let $c>0$ be a period of such a function $\theta$.
Then $\theta(c\,\cdot)$ is a shape that is not constant, so by (i) its periods are $\{j/d:j\in\mathbb{Z}\}$ for some $d\in\mathbb{N}$, and the periods of $\theta$ are $\{jc/d:j\in\mathbb{Z}\}$.
The smallest positive period of $\theta$ is $c/d$.

For (iii), let $\theta$ and $\theta'$ have frequencies $\tau$ and $\tau'$ and be equal almost everywhere on $[0,1/\tau+1/\tau']$.
For almost every $s\in[0,1/\tau']$ we have $\theta'(s+1/\tau)=\theta(s+1/\tau)=\theta(s)=\theta'(s)$.
The function $s\mapsto\theta'(s+1/\tau)-\theta'(s)$ has period $1/\tau'$ and vanishes almost everywhere on $[0,1/\tau']$, so it vanishes almost everywhere on $\mathbb{R}$, and $1/\tau$ is a period of $\theta'$.
Then $\theta$ and $\theta'$ both have period $1/\tau$ and are equal almost everywhere on $[0,1/\tau]$, so they are equal.
\qed

\paragraph{Proof of Lemma~\ref{lemma:xi_sets}:}
Assertion (i) is clear.
For $\bm x\in\ell^2$ and $p\in\mathbb{N}$, write $S_p=\sum_{i\in\mathbb{Z}^*}x_{pi}^2$, the sum of the squared amplitudes of the harmonics with indices in $\{p,2p,3p,\dots\}$.
For (ii), let $f$ be a shape with coordinates $\bm x\in\Xi(h)$, for some $h\in(0,1]$.
Then $S_1>0$, so $f$ is not constant; let $d$ be as in Lemma~\ref{lemma:periods}(i).
If $d\ge2$, then $x_i=0$ unless $d$ divides $i$, so that $S_d=S_1>(1-h)S_1$, which contradicts $\bm x\in\Xi(h)$.
Hence $d=1$ and $f\in\mathcal{F}$.
Conversely, let $f\in\mathcal{F}$ have coordinates $\bm x$.
Then $S_1>0$, and, since $d=1$, for every $p\ge2$ some non-zero harmonic of $f$ has an index $k$ that is not a multiple of $p$, so that $S_p\le S_1-(x_k^2+x_{-k}^2)<S_1$.
Moreover, $S_p\le\sum_{|i|\ge p}x_i^2\to0$ as $p\to\infty$, so there is an integer $p_0\ge2$ such that $S_p\le S_1/2$ for all $p\ge p_0$.
Then $r=\max\{1/2,S_2/S_1,\dots,S_{p_0-1}/S_1\}<1$, and $\bm x\in\Xi(1-r)$.
For (iii), note that $S_1\ge a^2>0$ and that $\{p,2p,3p,\dots\}$ does not contain 1 if $p\ge2$, so that, for such $p$, $S_p\le S_1-a^2=(1-a^2/S_1)S_1$, with $a^2/S_1\in(0,1]$.
For (iv), if $\bm y$ are the coordinates of a constant shape, then $\|\bm x-\bm y\|_2^2\ge S_1\ge hS_1$.
If $\bm y$ are the coordinates of a shape that is not constant and has period $1/p$, then $1/p=j/d$ for some $j\in\mathbb{Z}$, with $d$ as in Lemma~\ref{lemma:periods}(i) for this shape, so $p$ divides $d$ and $y_i=0$ unless $p$ divides $i$.
Hence $\|\bm x-\bm y\|_2^2\ge\sum_{i\in\mathbb{Z}^*,\,p\nmid i}x_i^2=S_1-S_p\ge hS_1$, where $p\nmid i$ means that $p$ does not divide $i$.
Finally, let $f\in\mathcal{F}$, with coordinates $\bm x\in\Xi(h)$ by (ii), and let $g$ be a shape with coordinates $\bm y$ such that $\|\bm x-\bm y\|_2^2<hS_1$.
By the first part of (iv), $g$ is not constant and has no period $1/p$ with $p\ge2$, so by Lemma~\ref{lemma:periods}(i) its smallest positive period is 1, and $g\in\mathcal{F}$.
\qed

\paragraph{Proof of Proposition~\ref{prop:information}:}
We first prove the differentiability claimed in Section~\ref{sec:statement_of_problem:information}.
For $\delta\neq0$, the absolute continuity of $f$ gives $f\big((\tau+\delta)t\big)-f(\tau t)=\delta t\int_0^1f'\big((\tau+r\delta)t\big)\,dr$, so that, by Minkowski's integral inequality,
\[
\sqrt n\,\bigg\|\frac{\psi(\tau+\delta,f)-\psi(\tau,f)}{\delta}-tf'(\tau t)\bigg\|_{2,n}\le n^{3/2}\sup_{|\sigma-\tau|\le|\delta|}\big\|f'(\sigma\,\cdot)-f'(\tau\,\cdot)\big\|_{2,n}.
\]
The right side tends to zero as $\delta\to0$, since $\sigma\mapsto f'(\sigma\,\cdot)$ is continuous in $L^2\big([0,n]\big)$, as follows by approximating $f'$ in $L^2\big([0,1]\big)$ by continuous 1-periodic functions.

Write $\kappa=n\tau$.
Let $R$ be a 1-periodic function, integrable on $[0,1]$, with mean $\mu=\int_0^1R(u)\,du$, and let $G(u)=\int_0^u\{R(v)-\mu\}\,dv$, $u\ge0$.
The function $G$ is absolutely continuous and 1-periodic, with $G(0)=0$, and since also $G(u)=-\int_u^1\{R(v)-\mu\}\,dv$ for $u\in[0,1]$,
\[
|G(u)|\le\frac12\int_0^1|R(v)-\mu|\,dv\le\int_0^1|R(v)|\,dv,\quad u\ge0.
\]
For an absolutely continuous function $\omega$ on $[0,\kappa]$, integration by parts gives
\[
\int_0^\kappa\omega(u)R(u)\,du=\mu\int_0^\kappa\omega(u)\,du+\omega(\kappa)G(\kappa)-\int_0^\kappa\omega'(u)G(u)\,du,
\]
so that
\begin{equation}\label{eq:proof_information_ibp}
\bigg|\int_0^\kappa\omega(u)R(u)\,du-\mu\int_0^\kappa\omega(u)\,du\bigg|\le\int_0^1|R(v)|\,dv\,\bigg\{|\omega(\kappa)|+\int_0^\kappa|\omega'(u)|\,du\bigg\}.
\end{equation}
For the first bound in~\eqref{eq:information_1}, take $R=g^2$ and $\omega\equiv1$.
By the variable change $u=\tau t$, $\|g(\tau\,\cdot)\|_{2,n}^2$ is the average of $R$ over $[0,\kappa]$, and~\eqref{eq:proof_information_ibp} gives
\[
\Big|\|g(\tau\,\cdot)\|_{2,n}^2-\|g\|_2^2\Big|\le\frac{\|g\|_2^2}\kappa.
\]

For the next two bounds take $R=(f')^2$, so that $\mu=\int_0^1|R(v)|\,dv=\|f'\|_2^2$, which is positive since $f$ is absolutely continuous and not constant.
By the variable change $u=\tau t$,
\[
\mathcal{I}_n(\tau,f)=\frac1{\tau^3}\int_0^\kappa u^2R(u)\,du,
\qquad
\mathcal{I}_n^{(n/2)}(\tau,f)=\frac1{\tau^3}\int_0^\kappa\Big(u-\frac\kappa2\Big)^2R(u)\,du.
\]
With $\omega(u)=u^2$ we have $\int_0^\kappa\omega(u)\,du=\kappa^3/3$ and $|\omega(\kappa)|+\int_0^\kappa|\omega'(u)|\,du=2\kappa^2$, and with $\omega(u)=(u-\kappa/2)^2$ these quantities are $\kappa^3/12$ and $3\kappa^2/4$.
Hence~\eqref{eq:proof_information_ibp} gives
\[
\bigg|\mathcal{I}_n(\tau,f)-\frac{\|f'\|_2^2\kappa^3}{3\tau^3}\bigg|\le\frac{2\|f'\|_2^2\kappa^2}{\tau^3},
\qquad
\bigg|\mathcal{I}_n^{(n/2)}(\tau,f)-\frac{\|f'\|_2^2\kappa^3}{12\tau^3}\bigg|\le\frac{3\|f'\|_2^2\kappa^2}{4\tau^3}.
\]
Since $\kappa^3/\tau^3=n^3$, these are the second bound in~\eqref{eq:information_1} and the first bound in~\eqref{eq:information_2}.
For~\eqref{eq:information_3}, take $R=f'g$, so that $\int_0^1|R(v)|\,dv\le\|f'\|_2\|g\|_2$ by the Cauchy-Schwarz inequality, and $\omega(u)=u-\kappa/2$, so that $\int_0^\kappa\omega(u)\,du=0$ and $|\omega(\kappa)|+\int_0^\kappa|\omega'(u)|\,du=3\kappa/2$.
Since the integral in~\eqref{eq:information_3} is $\tau^{-2}\int_0^\kappa\omega(u)R(u)\,du$, by the variable change $u=\tau t$, this gives~\eqref{eq:information_3} with the constant $3/2$ in place of $1/2$.
For the constant $1/2$, write $\kappa=m+r$ with $m=\lfloor\kappa\rfloor$ and $r\in[0,1)$.
Splitting $[0,\kappa]$ into the intervals $[j,j+1)$, $j=0,\dots,m-1$, and $[m,\kappa]$, and using the periodicity of $R$,
\[
\int_0^\kappa\Big(u-\frac\kappa2\Big)R(u)\,du=\int_0^1w(v)R(v)\,dv,
\qquad
w(v)=\sum_{j=0}^{m-1}\Big(j+v-\frac\kappa2\Big)+\Big(m+v-\frac\kappa2\Big)\chi(v),
\]
where $\chi(v)=1$ if $v<r$ and $\chi(v)=0$ otherwise.
The sum equals $m\{v-(1+r)/2\}$, and $m-\kappa/2=(m-r)/2$, so that $w(v)=(m+1)(v-r/2)$ for $v\in[0,r)$, and $w(v)=m\{v-(1+r)/2\}$ for $v\in[r,1)$.
Hence $\sup_{v\in[0,1)}|w(v)|\le\max\{(m+1)r,m(1-r)\}/2\le\kappa/2$, since $(m+1)r=mr+r\le m+r$ and $m(1-r)\le m$.
It follows that $\big|\int_0^\kappa\omega(u)R(u)\,du\big|\le(\kappa/2)\int_0^1|R(v)|\,dv\le(\kappa/2)\|f'\|_2\|g\|_2$, and the variable change $u=\tau t$ gives~\eqref{eq:information_3}.

For the last bound in~\eqref{eq:information_2}, let $q\in\mathbb{N}$ and
\[
V(q)=\inf_{g}\int_0^{q/\tau}\big\{tf'(\tau t)-g(\tau t)\big\}^2\,dt.
\]
By the variable change $u=\tau t$ and the periodicity of $f'$ and $g$,
\[
\int_0^{q/\tau}\big\{tf'(\tau t)-g(\tau t)\big\}^2\,dt
=\frac1\tau\int_0^1\sum_{j=0}^{q-1}\Big\{\frac{j+u}\tau f'(u)-g(u)\Big\}^2\,du.
\]
For fixed $u$, the sum is minimised over the value $g(u)$ by the mean $f'(u)\{u+(q-1)/2\}/\tau$ of the numbers $(j+u)f'(u)/\tau$, $j=0,\dots,q-1$, and the minimum is
\[
\frac{\{f'(u)\}^2}{\tau^2}\sum_{j=0}^{q-1}\Big(j-\frac{q-1}2\Big)^2=\frac{\{f'(u)\}^2}{\tau^2}\,\frac{q(q^2-1)}{12}.
\]
The minimising function $g(u)=f'(u)\{u+(q-1)/2\}/\tau$, $u\in[0,1)$, extended to $\mathbb{R}$ by periodicity, is a shape, so $V(q)=\|f'\|_2^2q(q^2-1)/(12\tau^3)$; for $q=\kappa$ this $g$ is the function $g^*$ of Remark~\ref{rem:least_favourable}.
Since the integrand in the definition of $\tilde{\mathcal{I}}_n(\tau,f)$ is non-negative and $\lfloor\kappa\rfloor/\tau\le n\le\lceil\kappa\rceil/\tau$, we have $V(\lfloor\kappa\rfloor)\le\tilde{\mathcal{I}}_n(\tau,f)\le V(\lceil\kappa\rceil)$.
The function $x\mapsto x(x^2-1)$ is increasing on $[1,\infty)$, and $\kappa-1\le\lfloor\kappa\rfloor\le\lceil\kappa\rceil\le\kappa+1$ with $\kappa-1\ge1$, so
\[
\frac{(\kappa-1)(\kappa-2)}{\kappa^2}\le\frac{\tilde{\mathcal{I}}_n(\tau,f)}{\|f'\|_2^2\kappa^3/(12\tau^3)}\le\frac{(\kappa+1)(\kappa+2)}{\kappa^2}.
\]
The left side is at least $1-3/\kappa$, and, since $\kappa\ge2$, the right side is at most $1+4/\kappa$.
\qed

\section{Proofs for Section~\ref{sec:estimation}}\label{sec:technical:estimation}

In this appendix, for $\epsilon>0$ and a set $S$ of parameters, $N(\epsilon,S,\|\cdot\|_{2,n})$ is the minimal number of balls of radius $\epsilon$ in $L^2\big([0,n]\big)$, for the norm $\|\cdot\|_{2,n}$, needed to cover the set of signals of the parameters in $S$, restricted to $[0,n]$.
For $x\in\mathbb{R}$ we write $\|x\|_{\mathbb{Z}}=\min_{j\in\mathbb{Z}}|x-j|$ for the distance from $x$ to the nearest integer.

\begin{lemma}[Posterior contraction]\label{lemma:contraction_general}
Let $\epsilon_n$ be a positive sequence such that $\epsilon_n\to0$ and $n\epsilon_n^2\to\infty$.
Assume that there exist measurable sets $S_n$ of parameters $(\tau,K,\bm Z)$, and constants $C_1,\dots,C_4>0$, such that, for all large enough $n$,
\begin{align}
\log N\big(C_1\epsilon_n,S_n,\|\cdot\|_{2,n}\big)&\le C_2n\epsilon_n^2,\label{eq:condition_entropy}\\
\Pi\big(S_n^c\big)&\lesssim\exp\big\{-n\epsilon_n^2(C_3^2+4C_4^2)\big\},\label{eq:condition_exceptional_set}\\
\Pi\big(\|\theta-\theta_0\|_{2,n}\le C_4\epsilon_n\big)&\ge\exp\big(-n\epsilon_n^2C_3^2\big),\label{eq:condition_prior_mass}
\end{align}
where $S_n^c$ is the complement of $S_n$.
Then, for all large enough $L$,
\[
\Pi\big(\|\theta-\theta_0\|_{2,n}>L\epsilon_n\mid X^{(n)}\big)\to0
\quad\text{and}\quad
\Pi\big(S_n\mid X^{(n)}\big)\to1,
\]
in $P_{\theta_0}^{(n)}$-probability, as $n\to\infty$.
\end{lemma}
\begin{proof}
The lemma follows from~\cite[Theorem~1 and Lemma~1]{Ghosal:2007}, with $k=2$ there, and we only indicate the reduction.
We apply these results on the parameter space of $(\tau,K,\bm Z)$, with the semimetric $d_n=e_n$ given by the distance $\|\cdot\|_{2,n}$ between the signals, which their framework allows.
The process $\tilde Y_u=Y_{nu}/n$, $u\in[0,1]$, satisfies $d\tilde Y_u=\theta(nu)\,du+n^{-1/2}\,d\tilde W_u$ for the standard Brownian motion $\tilde W_u=W_{nu}/n^{1/2}$, and $\|\theta(n\,\cdot)\|_2=\|\theta\|_{2,n}$.
This is the white noise model of~\cite[Section~5]{Ghosal:2007}, with the signal $\theta(n\,\cdot)$ in place of $\theta$.
By~\cite[Lemma~5]{Ghosal:2007} and the bound $1-\Phi(x)\le e^{-x^2/2}$, $x\ge0$, tests as in~\cite[(2.2)]{Ghosal:2007} exist, with the constants denoted there by $\xi$ and $K$ equal to $1/4$ and $1/32$.
By~\cite[Lemma~6]{Ghosal:2007}, the Kullback-Leibler neighbourhoods $B_n(\theta_0,\epsilon;2)$ there are the balls $\{\|\theta-\theta_0\|_{2,n}\le\epsilon\}$.

For the second assertion, apply~\cite[Lemma~1]{Ghosal:2007} with $C_4\epsilon_n$ in place of $\epsilon_n$: by~\eqref{eq:condition_exceptional_set} and~\eqref{eq:condition_prior_mass}, $\Pi(S_n^c)/\Pi(\|\theta-\theta_0\|_{2,n}\le C_4\epsilon_n)\lesssim e^{-4C_4^2n\epsilon_n^2}$, which is $o(e^{-2C_4^2n\epsilon_n^2})$.

For the first assertion, apply~\cite[Theorem~1]{Ghosal:2007} to the sets $S_n$, with $\bar\epsilon_n=C_5\epsilon_n$ in place of $\epsilon_n$ and $C_5^2\ge\max\{64C_1^2,C_2,C_4^2,64C_3^2\}$.
Condition~(2.4) there holds, since $\epsilon/8\ge C_1\epsilon_n$ for $\epsilon>\bar\epsilon_n$, so that, by~\eqref{eq:condition_entropy}, the local entropy is at most $C_2n\epsilon_n^2\le n\bar\epsilon_n^2$.
Condition~(2.5) there holds for every $j\in\mathbb{N}$, since its numerator is at most 1 and, by~\eqref{eq:condition_prior_mass} and since $\bar\epsilon_n\ge C_4\epsilon_n$, its denominator is at least $e^{-C_3^2n\epsilon_n^2}\ge e^{-n\bar\epsilon_n^2j^2/64}$.
Since $n\bar\epsilon_n^2\to\infty$, the bound~(8.6) in the proof of~\cite[Theorem~1]{Ghosal:2007} gives the conclusion~(2.6) there for a fixed, large enough multiple of $\bar\epsilon_n$, instead of $M_n\bar\epsilon_n$ with $M_n\to\infty$.
Hence $\Pi\big(S_n\cap\{\|\theta-\theta_0\|_{2,n}>L\epsilon_n\}\mid X^{(n)}\big)\to0$ for all large enough $L$, and, with the second assertion, the first follows.
\qed\end{proof}

\begin{lemma}[Distance between signals]\label{lemma:upper_bound_metrics}
Let $k\in\mathbb{N}$, $\bm{x},\bm{y}\in\mathbb{R}^{2k+1}$, $M>0$ with $\bm{y}\in w_1(M)$, and $\tau,\tau_1>0$ with $n\tau\ge2$.
Then
\[
\|\bm{x}^T\bm{\varphi}(\tau\,\cdot)-\bm{y}^T\bm{\varphi}(\tau_1\,\cdot)\|_{2,n}\le\Big(1+\frac1{n\tau}\Big)^{1/2}\|\bm{x}-\bm{y}\|_2+\Big(\frac{4k}3\Big)^{1/2}nM|\tau-\tau_1|.
\]
\end{lemma}
\begin{proof}
By the triangle inequality, the left side is at most $\|(\bm x-\bm y)^T\bm\varphi(\tau\,\cdot)\|_{2,n}+\|\bm{y}^T\bm{\varphi}(\tau\,\cdot)-\bm{y}^T\bm{\varphi}(\tau_1\,\cdot)\|_{2,n}$.
The first bound in~\eqref{eq:information_1} bounds the first term by $(1+1/(n\tau))^{1/2}\|\bm x-\bm y\|_2$.
By the mean value theorem, $|\varphi_i(\tau t)-\varphi_i(\tau_1t)|\le\sqrt2\,v_it|\tau-\tau_1|$, $i\in\mathbb{Z}^*$, so that, by the Cauchy-Schwarz inequality and since $\bm y\in w_1(M)$,
\begin{align*}
\big|\bm{y}^T\bm{\varphi}(\tau t)-\bm{y}^T\bm{\varphi}(\tau_1t)\big|&\le\sqrt2\,t|\tau-\tau_1|\sum_{0<|i|\le k}v_i|y_i|\\
&\le\sqrt2\,t|\tau-\tau_1|(2k)^{1/2}\bigg(\sum_{0<|i|\le k}v_i^2y_i^2\bigg)^{1/2}\le2k^{1/2}Mt|\tau-\tau_1|.
\end{align*}
Squaring and integrating over $[0,n]$ gives $\|\bm{y}^T\bm{\varphi}(\tau\,\cdot)-\bm{y}^T\bm{\varphi}(\tau_1\,\cdot)\|_{2,n}^2\le4kM^2(\tau-\tau_1)^2n^2/3$, and the statement of the lemma follows.
\qed\end{proof}

\paragraph{Proof of Theorem~\ref{theo:checking_conditions}:}
The result follows if we show that the conditions of Lemma~\ref{lemma:contraction_general} hold for the sets $S_n=\Theta_n$ in~\eqref{eq:sieve_adaptive}, with $C_1=C_3=C_4=1$ and some $C_2>0$.
Throughout, $n$ is large enough that $n\ge2l_n$.
By~\eqref{eq:cond}, $\bar k_n\lesssim n\epsilon_n^2/\log n\le n$, so $\bar k_n$, $m_n$ and $u_n$ are bounded by a polynomial in $n$.

\emph{Checking condition~\eqref{eq:condition_entropy}}:
Let $(\tau_1,K,\bm y)\in\Theta_n$, with $K=k\le\bar k_n$, and let $\tau\in[l_n^{-1},u_n]$ and $\bm x\in\mathbb{R}^{2k+1}$.
Since $n\tau\ge n/l_n\ge2$ and $\bm y\in w_1(m_n)$, Lemma~\ref{lemma:upper_bound_metrics} gives
\[
\|\psi(\tau,\bm{x}^T\bm{\varphi})-\psi(\tau_1,\bm{y}^T\bm{\varphi})\|_{2,n}\le\sqrt2\,\|\bm{x}-\bm{y}\|_2+\big(4\bar k_n/3\big)^{1/2}nm_n|\tau-\tau_1|.
\]
If $|z_0|\le m_n$ and $\bm{z}\in w_1(m_n)$, then $\|\bm{z}\|_\infty\le m_n$, since $v_i\ge1$, $i\in\mathbb{Z}^*$.
Let $r>0$.
Cover $[l_n^{-1},u_n]$ by at most $1+2\bar k_n^{1/2}nm_nu_n/r$ intervals of length $r/(2\bar k_n^{1/2}nm_n)$ with centres in $[l_n^{-1},u_n]$, and, for each $k\le\bar k_n$, the cube $\{\bm z\in\mathbb{R}^{2k+1}:\|\bm z\|_\infty\le m_n\}$ by at most $\{1+4m_n(2k+1)^{1/2}/r\}^{2k+1}$ cubes of side $r/\{2(2k+1)^{1/2}\}$.
If $\tau$ is the centre of an interval that contains $\tau_1$ and $\bm x$ the centre of a cube that contains $\bm y$, then $|\tau-\tau_1|\le r/(4\bar k_n^{1/2}nm_n)$ and $\|\bm x-\bm y\|_2\le r/4$, and the previous display gives a distance at most $(\sqrt2+2/\sqrt3)\,r/4\le r$.
Hence
\begin{align*}
\log N\big(r,\Theta_n,\|\cdot\|_{2,n}\big)&\le\log\Big(1+\frac{2\bar k_n^{1/2}nm_nu_n}{r}\Big)+\log\bar k_n\\
&\quad+(2\bar k_n+1)\log\Big(1+\frac{4m_n(2\bar k_n+1)^{1/2}}{r}\Big).
\end{align*}
For $r=\epsilon_n$, we have $\epsilon_n\ge n^{-1/2}$ for large $n$, since $n\epsilon_n^2\to\infty$, so every logarithm in the previous display is $\lesssim\log n$.
Hence $\log N(\epsilon_n,\Theta_n,\|\cdot\|_{2,n})\lesssim\bar k_n\log n\lesssim n\epsilon_n^2$, by~\eqref{eq:cond}, and the condition holds for some $C_2$.

\emph{Checking condition~\eqref{eq:condition_exceptional_set}}:
By conditioning appropriately,
\begin{align*}
\Pi\big(\Theta_n^c\big)&\le\Pi_K(K>\bar k_n)+\Pi_\tau\big(\tau\notin[l_n^{-1},u_n]\big)\\
&\quad+\sum_{k=1}^{\bar k_n}\Big\{\Pi_{\bm{Z}}\big(|Z_0|>m_n\mid K=k\big)+\Pi_{\bm{Z}}\big(\bm{Z}\notin w_1(m_n)\mid K=k\big)\Big\}\Pi_K(K=k).
\end{align*}
By~\eqref{eq:cond_K_2}, \eqref{eq:cond_Z_0}, \eqref{eq:cond_Z_3} and~\eqref{eq:cond_tau_2}, the previous display is $\lesssim\exp(-5n\epsilon_n^2)$, so the condition holds with $C_3^2+4C_4^2=5$.

\emph{Checking condition~\eqref{eq:condition_prior_mass}}:
Since $\alpha\ge1$ and $v_i\ge1$, $\sum_{0<|i|\le k_n}v_i^2x_{0,i}^2\le\sum_{i\in\mathbb{Z}^*}v_i^{2\alpha}x_{0,i}^2\le M^2$, i.e., $\bm z_n\in w_1(M)$.
By~\eqref{eq:truncation} and since $2\pi(k_n+1)\ge\pi(2k_n+1)$,
\[
\|f_0-f_{0,[k_n]}\|_2\le\frac{M}{\pi^\alpha(2k_n+1)^\alpha}=\frac{C_\alpha}{(2k_n+1)^\alpha}\le\frac{\epsilon_n}6,
\]
by~\eqref{eq:cond}.
Let $\bm{z}\in\mathbb{R}^{2k_n+1}$ with $\|\bm{z}-\bm{z}_n\|_2\le\epsilon_n/2$, and $\tau$ with $|\tau-\tau_0|\le\epsilon_n/(4\sqrt3\,k_n^{1/2}nM)$, so that $\tau\ge\tau_0/2$ for large $n$.
Then $(4k_n/3)^{1/2}nM|\tau-\tau_0|\le\epsilon_n/6$.
By the triangle inequality, Lemma~\ref{lemma:upper_bound_metrics} with $\bm y=\bm z_n\in w_1(M)$, and the first bound in~\eqref{eq:information_1} applied to the shape $f_0-f_{0,[k_n]}$,
\begin{align*}
\|\psi(\tau,\bm{z}^T\bm{\varphi})-\theta_0\|_{2,n}&\le\|\bm{z}^T\bm{\varphi}(\tau\,\cdot)-\bm{z}_n^T\bm{\varphi}(\tau_0\,\cdot)\|_{2,n}+\|(f_0-f_{0,[k_n]})(\tau_0\,\cdot)\|_{2,n}\\
&\le\Big(1+\frac2{n\tau_0}\Big)^{1/2}\frac{\epsilon_n}2+\frac{\epsilon_n}6+\Big(1+\frac1{n\tau_0}\Big)^{1/2}\frac{\epsilon_n}6,
\end{align*}
which is at most $\epsilon_n$ for large $n$.
Since the density of $\Pi_\tau$ is bounded away from zero near $\tau_0$, it follows from~\eqref{eq:cond_K_1} and~\eqref{eq:cond_Z_1} that
\begin{align*}
\Pi\big(\|\theta-\theta_0\|_{2,n}\le\epsilon_n\big)&\ge
\Pi_K(K=k_n)\,
\Pi_{\bm{Z}}\big(\|\bm{Z}-\bm{z}_n\|_2\le\epsilon_n/2\mid K=k_n\big)\\
&\quad\times\Pi_\tau\Big(|\tau-\tau_0|\le\frac{\epsilon_n}{4\sqrt3\,k_n^{1/2}nM}\Big)\\
&\gtrsim\exp\Big(-\frac{n\epsilon_n^2}2-\log\frac{k_n^{1/2}nM}{\epsilon_n}\Big),
\end{align*}
which is at least $\exp(-n\epsilon_n^2)$ for large $n$, since the logarithm is $\lesssim\log n$, and $\log n=o(n\epsilon_n^2)$ by~\eqref{eq:cond}, as $\bar k_n\ge k_n\to\infty$ by the last condition in~\eqref{eq:cond}, since $\epsilon_n\to0$.
The condition holds with $C_3=C_4=1$.
\qed

The proof of Proposition~\ref{prop:smaller_sieve} rests on two observations, made precise in Lemmas~\ref{lemma:shift_average} and~\ref{lemma:near_diagonal}.
First, $\theta_0$ has period $1/\tau_0$, so a signal $(\tau,f)$ close to $\theta_0$ is close to its own shifts by multiples of $1/\tau_0$, and the shape $f$ is then close to its own shifts by multiples of $\tau/\tau_0$.
The condition $n\tau_0\ge16\bar k$ gives at least $4\bar k$ such shifts, enough to detect that $k\tau/\tau_0$ is at distance more than $1/(4\bar k)$ from the integers.
This puts most of the squared amplitudes of $f$ on harmonics $k$ for which $k\tau/\tau_0$ is close to an integer, so that $\tau/\tau_0$ is close to a rational number $p/q$.
For $f\in\Xi(h)$ only $q=1$ remains, and $p=0$ is excluded, since $\tau\ge l^{-1}$ and $n\ge4l$.
Symmetrically, $\theta_0$ is then close to its own shifts by multiples of $1/\tau$, which excludes $p\ge2$, since $f_0\in\Xi(H)$.
Second, let $\tau=\tau_0(1+\zeta)$ with $\zeta$ small.
At time $t=u/\tau_0$ the error in the frequency shifts the shape by $\zeta u$ periods, and Lemma~\ref{lemma:near_diagonal} shows that $\|\psi(\tau,f)-\psi(\tau_0,f_{0,[\bar k]})\|_{2,n}^2$ is the average over $u\in[0,n\tau_0]$ of $\|f(\cdot+\zeta u)-f_{0,[\bar k]}\|_2^2$, up to a term of order $(1+\log\bar k)(\|f\|_2^2+\|f_0\|_2^2)/(n\tau_0)$.
Similar reductions to rational ratios appear in~\cite[proof of Lemma~2]{Golubev:1988}, \cite[Lemma~1]{Castillo:2007} and~\cite[Lemmas~5 and~6]{LevyLeduc:2008}, where they locate the maximum of a contrast, evaluated at the true signal, near the rational multiples of the true period or frequency.
For a fixed and known shape, a lower bound of the type obtained in Case~5 of the proof is~\cite[Lemma~III.5.3]{Ibragimov:1981}.

\begin{lemma}[Shift averaging]\label{lemma:shift_average}
The following statements hold:
\begin{enumerate}[label=(\roman*)]
\item For every integer $Q\ge2$ and every $x\in\mathbb{R}$, with $\|x\|_{\mathbb{Z}}$ the distance from $x$ to the nearest integer (see the start of this appendix),
\[
\frac1Q\sum_{j=1}^Q4\sin^2(\pi jx)\ge\frac16\min\big\{1,(Q\|x\|_{\mathbb{Z}})^2\big\}.
\]
\item Let $f$ and $g$ be shapes, and let $\tau,\tau_1>0$ and $n>0$ be such that $n\tau\ge4$ and $n\tau_1\ge4$.
Let $a_k$ be the amplitude of the $k$-th harmonic of $f$, $k\in\mathbb{N}$, let $Q=\lfloor n\tau_1/2\rfloor$ and $\rho=\tau/\tau_1$.
Then
\[
\sum_{k=1}^\infty a_k^2\min\big\{1,(Q\|k\rho\|_{\mathbb{Z}})^2\big\}\le96\,\|\psi(\tau,f)-\psi(\tau_1,g)\|_{2,n}^2.
\]
\end{enumerate}
\end{lemma}
\begin{proof}
For (i), suppose first that $Q\|x\|_{\mathbb{Z}}\ge1$.
Since $|\sin(\pi x)|=\sin(\pi\|x\|_{\mathbb{Z}})\ge2\|x\|_{\mathbb{Z}}$,
\[
\bigg|\sum_{j=1}^Q\cos(2\pi jx)\bigg|\le\bigg|\sum_{j=1}^Qe^{2\pi\mathrm{i}jx}\bigg|=\frac{|\sin(\pi Qx)|}{|\sin(\pi x)|}\le\frac1{2\|x\|_{\mathbb{Z}}},
\]
so that $Q^{-1}\sum_{j=1}^Q4\sin^2(\pi jx)=2-2Q^{-1}\sum_{j=1}^Q\cos(2\pi jx)\ge2-1/(Q\|x\|_{\mathbb{Z}})\ge1$.
Suppose now that $Q\|x\|_{\mathbb{Z}}<1$, and let $m=\lfloor Q/2\rfloor$, so that $m\ge Q/3$.
For $j\le m$ we have $j\|x\|_{\mathbb{Z}}<1/2$, hence $\|jx\|_{\mathbb{Z}}=j\|x\|_{\mathbb{Z}}$ and $\sin^2(\pi jx)\ge4j^2\|x\|_{\mathbb{Z}}^2$.
Therefore
\[
\frac1Q\sum_{j=1}^Q4\sin^2(\pi jx)\ge\frac{16\|x\|_{\mathbb{Z}}^2}Q\sum_{j=1}^mj^2\ge\frac{16\|x\|_{\mathbb{Z}}^2m^3}{3Q}\ge\frac{16}{81}(Q\|x\|_{\mathbb{Z}})^2,
\]
and $16/81\ge1/6$; we use the simpler constant $1/6$ in the statement.

For (ii), write $\theta=\psi(\tau,f)$, $\theta_1=\psi(\tau_1,g)$ and $\epsilon=\|\theta-\theta_1\|_{2,n}$; the shape $g$ only enters through $\epsilon$.
Let $j\in\{1,\dots,Q\}$ and $s=j/\tau_1$, so that $s\le n/2$, and $s$ is a period of $\theta_1$.
For $t\in[0,n-s]$, $\theta(t+s)-\theta(t)=\{\theta(t+s)-\theta_1(t+s)\}-\{\theta(t)-\theta_1(t)\}$, so that
\[
\int_0^{n-s}\big\{\theta(t+s)-\theta(t)\big\}^2\,dt\le4n\epsilon^2.
\]
On the other hand, $\theta(t+s)-\theta(t)=g_j(\tau t)$ for the shape $g_j=f(\cdot+j\rho)-f$, since $\tau s=j\rho$.
By the first bound in~\eqref{eq:information_1}, with $n-s$ in place of $n$, and since $(n-s)\tau\ge n\tau/2\ge2$,
\[
\int_0^{n-s}\big\{g_j(\tau t)\big\}^2\,dt\ge(n-s)\Big(1-\frac1{(n-s)\tau}\Big)\|g_j\|_2^2\ge\frac n4\|g_j\|_2^2.
\]
Hence $\|f(\cdot+j\rho)-f\|_2^2\le16\epsilon^2$.
By~\eqref{eq:shift}, $\|f(\cdot+j\rho)-f\|_2^2=\sum_{k=1}^\infty a_k^2\,4\sin^2(\pi kj\rho)$.
Averaging over $j=1,\dots,Q$, with $Q\ge2$ since $n\tau_1\ge4$, and using (i) with $x=k\rho$ gives the assertion, with the constant $96=6\times16$ from the factor 16 above and the factor 6 from (i); neither factor is optimised.
\qed\end{proof}

\begin{lemma}[Near-diagonal expansion]\label{lemma:near_diagonal}
Let $\bar k\in\mathbb{N}$, $\kappa>0$ and $\zeta\in\mathbb{R}$ with $\bar k|\zeta|\le1/4$, and let $\rho=1+\zeta$.
Let $g$ and $\tilde g$ be trigonometric polynomials of degree at most $\bar k$.
Then
\[
\bigg|\frac1\kappa\int_0^\kappa\big\{g(\rho u)-\tilde g(u)\big\}^2\,du-\frac1\kappa\int_0^\kappa\|g(\cdot+\zeta u)-\tilde g\|_2^2\,du\bigg|\le\frac{14+8\log\bar k}\kappa\big(\|g\|_2^2+\|\tilde g\|_2^2\big).
\]
\end{lemma}
\begin{proof}
For $k=1,\dots,\bar k$, let $a_k$ and $b_k$ be the amplitudes of the $k$-th harmonics of $g$ and $\tilde g$, and write their pairs of coordinates as $a_k(\cos\gamma_k,\sin\gamma_k)$ and $b_k(\cos\tilde\gamma_k,\sin\tilde\gamma_k)$, so that the $k$-th harmonic of $g$ is $\sqrt2\,a_k\cos(2\pi ku-\gamma_k)$, and similarly for $\tilde g$.
Let $d_0=\bar g-\bar{\tilde g}$ be the difference of the means of $g$ and $\tilde g$, and, for $u,v\in\mathbb{R}$,
\[
T_k(v,u)=\sqrt2\,a_k\cos\big(2\pi k(v+\zeta u)-\gamma_k\big)-\sqrt2\,b_k\cos(2\pi kv-\tilde\gamma_k),
\]
so that $g(v+\zeta u)-\tilde g(v)=d_0+\sum_{k=1}^{\bar k}T_k(v,u)$, and $g(\rho u)-\tilde g(u)$ is this function at $v=u$.
By Parseval's identity in $v$,
\[
\|g(\cdot+\zeta u)-\tilde g\|_2^2=d_0^2+\sum_{k=1}^{\bar k}\Big\{a_k^2+b_k^2-2a_kb_k\cos(2\pi k\zeta u-\gamma_k+\tilde\gamma_k)\Big\}.
\]
We use repeatedly that, for $\lambda\neq0$ and $\xi\in\mathbb{R}$,
\begin{equation}\label{eq:proof_near_diagonal_cos}
\bigg|\frac1\kappa\int_0^\kappa\cos(2\pi\lambda u+\xi)\,du\bigg|\le\frac1{\pi\kappa|\lambda|}.
\end{equation}
Write $S=d_0^2+\sum_{k=1}^{\bar k}(a_k^2+b_k^2)$.
We expand the square of $d_0+\sum_kT_k(u,u)$ and compare the averages over $[0,\kappa]$ of its terms one by one with the previous display.
The terms $d_0^2$ agree.

By the product formula for cosines,
\begin{align*}
T_k(u,u)^2&=a_k^2+b_k^2-2a_kb_k\cos(2\pi k\zeta u-\gamma_k+\tilde\gamma_k)\\
&\quad+a_k^2\cos(4\pi k\rho u-2\gamma_k)+b_k^2\cos(4\pi ku-2\tilde\gamma_k)\\
&\quad-2a_kb_k\cos\big(2\pi k(1+\rho)u-\gamma_k-\tilde\gamma_k\big).
\end{align*}
The first line is the $k$-th term of the previous display.
The frequencies $2k\rho$, $2k$ and $k(1+\rho)$ in the second line are at least $3k/2$, since $\rho\ge3/4$, so, by~\eqref{eq:proof_near_diagonal_cos}, the average of the second line is at most $(a_k+b_k)^2\,2/(3\pi\kappa k)\le(a_k^2+b_k^2)/\kappa$ in absolute value.
The sum over $k$ of these bounds is at most $S/\kappa$.

The function $T_k(u,u)$ is a sum of two cosines with frequencies $k\rho\ge3k/4$ and $k$, so the average of $2d_0T_k(u,u)$ is at most $2|d_0|\sqrt2(a_k+b_k)\,4/(3\pi\kappa k)$ in absolute value.
By the Cauchy-Schwarz inequality and $\sum_kk^{-2}\le\pi^2/6$, the sum over $k$ of these bounds is at most
\[
\frac{8\sqrt2}{3\pi\kappa}|d_0|\frac\pi{\sqrt6}\bigg\{\sum_{k=1}^{\bar k}(a_k+b_k)^2\bigg\}^{1/2}\le\frac{16}{3\sqrt6\,\kappa}|d_0|\bigg\{\sum_{k=1}^{\bar k}(a_k^2+b_k^2)\bigg\}^{1/2}\le\frac{2S}{\kappa}.
\]

For $k\neq k'$, the product $T_k(u,u)T_{k'}(u,u)$ is, by the product formula, a sum of eight cosines with amplitudes whose total is at most $2(a_k+b_k)(a_{k'}+b_{k'})$, and with frequencies of the form $\lambda\pm\lambda'$, with $\lambda\in\{k\rho,k\}$ and $\lambda'\in\{k'\rho,k'\}$.
Since $\bar k|\zeta|\le1/4$, every such frequency is at least $|k-k'|-1/4\ge3|k-k'|/4$ in absolute value, so the average of $T_k(u,u)T_{k'}(u,u)$ is at most $8(a_k+b_k)(a_{k'}+b_{k'})/(3\pi\kappa|k-k'|)$ in absolute value.
With $c_k=a_k+b_k$, the Cauchy-Schwarz inequality gives
\[
\sum_{k\neq k'}\frac{c_kc_{k'}}{|k-k'|}=\sum_{j=1}^{\bar k-1}\frac2j\sum_{k=1}^{\bar k-j}c_kc_{k+j}\le2(1+\log\bar k)\sum_{k=1}^{\bar k}c_k^2\le4(1+\log\bar k)S,
\]
so that the sum over $k\neq k'$ of the averages is at most $32(1+\log\bar k)S/(3\pi\kappa)\le4(1+\log\bar k)S/\kappa$.

Collecting the bounds, the difference in the statement is at most $(7+4\log\bar k)S/\kappa$.
Finally, $S=\|g\|_2^2+\|\tilde g\|_2^2-2\bar g\,\bar{\tilde g}\le2(\|g\|_2^2+\|\tilde g\|_2^2)$.
\qed\end{proof}

\paragraph{Proof of Proposition~\ref{prop:smaller_sieve}:}
Let $a_k$ and $b_k$ be the amplitudes of the $k$-th harmonics of $f$ and $f_0$, $k\in\mathbb{N}$, so that $a_k=0$ for $k>\bar k$, and $\sum_{k\ge1}a_k^2=\|f-\bar f\|_2^2$.
Write $\theta=\psi(\tau,f)$, $\rho=\tau/\tau_0$, $Q=\lfloor n\tau_0/2\rfloor$ and $Q'=\lfloor n\tau/2\rfloor$.
By~\eqref{eq:new_conditions}, $n\tau_0\ge16$ and $n\tau\ge n/l\ge4$, so that $Q\ge n\tau_0/4\ge4\bar k$ and $Q'\ge n\tau/4$.
Also $400\epsilon^2\le\sigma_0^2$.
The constant 400 in~\eqref{eq:new_conditions} is a convenient value rather than an optimal one: it gives $\epsilon\le\sigma_0/20$, which is used in the next step, and Cases~1 to~5 below need $192\epsilon^2<\min\{h,H\}\sigma_0^2$, where $192=2\times96$ is the constant of Lemma~\ref{lemma:shift_average}(ii) times a factor 2 lost in each case.

We first bound $\|f-\bar f\|_2$ from below.
By the first bound in~\eqref{eq:information_1}, applied to the shapes $f_0-\bar f$ and $f-\bar f$, and since $\bar f_0$ minimises $c\mapsto\|f_0-c\|_2$,
\[
\Big(\frac54\Big)^{1/2}\|f-\bar f\|_2\ge\|\theta-\bar f\|_{2,n}\ge\|\theta_0-\bar f\|_{2,n}-\epsilon\ge\Big(\frac34\Big)^{1/2}\|f_0-\bar f\|_2-\epsilon\ge\Big\{\Big(\frac34\Big)^{1/2}-\frac1{20}\Big\}\sigma_0,
\]
so that $\|f-\bar f\|_2^2\ge\sigma_0^2/2$.

By Lemma~\ref{lemma:shift_average}(ii), applied to $(f,\tau)$ and $(f_0,\tau_0)$, and to $(f_0,\tau_0)$ and $(f,\tau)$,
\begin{equation}\label{eq:proof_smaller_sieve_shifts}
\sum_{k=1}^{\bar k}a_k^2\min\big\{1,(Q\|k\rho\|_{\mathbb{Z}})^2\big\}\le96\epsilon^2
\quad\text{and}\quad
\sum_{j=1}^\infty b_j^2\min\big\{1,(Q'\|j/\rho\|_{\mathbb{Z}})^2\big\}\le96\epsilon^2.
\end{equation}
Let $\mathcal{K}$ be the set of $k\in\{1,\dots,\bar k\}$ with $\|k\rho\|_{\mathbb{Z}}\le1/(4\bar k)$, and, for $k\in\mathcal{K}$, let $p_k$ be the integer nearest to $k\rho$, so that $p_k\ge0$.
For $k,k'\in\mathcal{K}$,
\[
|k'p_k-kp_{k'}|=\big|k'(p_k-k\rho)-k(p_{k'}-k'\rho)\big|\le\frac{2\bar k}{4\bar k}=\frac12,
\]
so that $k'p_k=kp_{k'}$.
Hence, if $\mathcal{K}$ is not empty, $p_k/k$ is the same rational number for all $k\in\mathcal{K}$, which we write as $p/q$ with $q\in\mathbb{N}$ and $p\ge0$ coprime integers (with $q=1$ if $p=0$).
Then $kp/q=p_k$ is an integer, so $q$ divides $k$ for every $k\in\mathcal{K}$, and $|\rho-p/q|=|k\rho-p_k|/k\le1/(4\bar kq)$ for every $k\in\mathcal{K}$.
Every $k\le\bar k$ not in $\mathcal{K}$ has $\|k\rho\|_{\mathbb{Z}}>1/(4\bar k)$, hence weight $\min\{1,(Q\|k\rho\|_{\mathbb{Z}})^2\}=1$ in~\eqref{eq:proof_smaller_sieve_shifts}, since $Q\ge4\bar k$.
We distinguish five cases, and show that the first four are impossible.

\emph{Case 1: $\mathcal{K}$ is empty.}
Then all weights in the first sum in~\eqref{eq:proof_smaller_sieve_shifts} are 1, so $\|f-\bar f\|_2^2\le96\epsilon^2<\sigma_0^2/2$, a contradiction.

\emph{Case 2: $q\ge2$.}
Then every $k\le\bar k$ that is not a multiple of $q$ is outside $\mathcal{K}$, so the first bound in~\eqref{eq:proof_smaller_sieve_shifts} gives $\sum_{k\le\bar k,\,q\nmid k}a_k^2\le96\epsilon^2$.
Since $f\in\Xi(h)$, the left side is $\|f-\bar f\|_2^2-\sum_{i\in\mathbb{Z}^*}x_{qi}^2\ge h\|f-\bar f\|_2^2\ge h\sigma_0^2/2$, where $\bm x$ are the coordinates of $f$.
This contradicts $400\epsilon^2\le h\sigma_0^2$.

\emph{Case 3: $q=1$ and $p=0$.}
Then $p_k=0$ for $k\in\mathcal{K}$, so that $\|k\rho\|_{\mathbb{Z}}=k\rho\ge\rho$ and $Q\|k\rho\|_{\mathbb{Z}}\ge Q\tau/\tau_0\ge n\tau/4\ge1$.
All weights are 1 again, and we conclude as in Case~1.

\emph{Case 4: $q=1$ and $p\ge2$.}
Then $|\rho-p|\le1/(4\bar k)\le1/4$, so that $\rho\ge p/2$ and $Q'\ge n\tau/4\ge np\tau_0/8$.
Let $j\in\{1,\dots,J\}$ not be a multiple of $p$.
Then $\|j/p\|_{\mathbb{Z}}\ge1/p$ and, since $J\le2\bar k$ by~\eqref{eq:new_conditions},
\[
\Big|\frac j\rho-\frac jp\Big|=\frac{j|p-\rho|}{\rho p}\le\frac{2J}{4\bar kp^2}\le\frac1{p^2}\le\frac1{2p},
\]
so that $\|j/\rho\|_{\mathbb{Z}}\ge1/(2p)$ and $Q'\|j/\rho\|_{\mathbb{Z}}\ge n\tau_0/16\ge1$.
The second bound in~\eqref{eq:proof_smaller_sieve_shifts} then gives $\sum_{j\le J,\,p\nmid j}b_j^2\le96\epsilon^2$.
Since $f_0\in\Xi(H)$, as in Case~2, $\sum_{j:\,p\nmid j}b_j^2\ge H\sigma_0^2$, and, by the definition of $J$, $\sum_{j>J}b_j^2\le H\sigma_0^2/2$.
Hence $H\sigma_0^2/2\le96\epsilon^2$, which contradicts $400\epsilon^2\le H\sigma_0^2$.

\emph{Case 5: $q=1$ and $p=1$.}
Then $\rho=1+\zeta$ with $|\zeta|\le1/(4\bar k)$, so that $k|\zeta|\le1/4$ and $\|k\rho\|_{\mathbb{Z}}=k|\zeta|$ for $k\le\bar k$.
By the first bound in~\eqref{eq:proof_smaller_sieve_shifts},
\[
\min\big\{1,(Q|\zeta|)^2\big\}\|f-\bar f\|_2^2\le\sum_{k=1}^{\bar k}a_k^2\min\big\{1,(Qk|\zeta|)^2\big\}\le96\epsilon^2.
\]
Since $96\epsilon^2<\|f-\bar f\|_2^2$, the minimum is $(Q|\zeta|)^2$, and $Q|\zeta|\le4\sqrt6\,\epsilon/\|f-\bar f\|_2\le8\sqrt3\,\epsilon/\sigma_0$.
Since $Q\ge n\tau_0/4$ and $32\sqrt3\le56$, this gives~\eqref{eq:localisation_tau}.

It remains to prove~\eqref{eq:localisation_f}, in Case~5.
Let $\kappa_0=n\tau_0$ and $\tilde f=f_{0,[\bar k]}$, and let $r=\|f_0-\tilde f\|_2$, so that $\|\tilde f\|_2\le\|f_0\|_2$.
By the variable change $u=\tau_0t$, the triangle inequality and the first bound in~\eqref{eq:information_1},
\[
\bigg\{\frac1{\kappa_0}\int_0^{\kappa_0}\big\{f(\rho u)-\tilde f(u)\big\}^2\,du\bigg\}^{1/2}=\|\theta-\psi(\tau_0,\tilde f)\|_{2,n}\le\epsilon+\Big(\frac54\Big)^{1/2}r.
\]
By the first bound in~\eqref{eq:information_1}, $\|f\|_2^2\le(4/3)\|\theta\|_{2,n}^2$ and $\|\theta_0\|_{2,n}^2\le(5/4)\|f_0\|_2^2$, so that
\[
\|f\|_2^2\le\frac43\big(\|\theta_0\|_{2,n}+\epsilon\big)^2\le\frac83\Big(\frac54\|f_0\|_2^2+\epsilon^2\Big)\le4\|f_0\|_2^2+4\epsilon^2.
\]
Since $\bar k|\zeta|\le1/4$, Lemma~\ref{lemma:near_diagonal}, applied to $g=f$, $\tilde g=\tilde f$ and $\kappa=\kappa_0$, gives
\[
\frac1{\kappa_0}\int_0^{\kappa_0}\|f(\cdot+\zeta u)-\tilde f\|_2^2\,du\le2\epsilon^2+3r^2+\frac{14+8\log\bar k}{\kappa_0}\big(5\|f_0\|_2^2+4\epsilon^2\big).
\]
Moreover, for every $u$, $\|f-\tilde f\|_2^2\le2\|f(\cdot+\zeta u)-\tilde f\|_2^2+2\|f(\cdot+\zeta u)-f\|_2^2$, and, by~\eqref{eq:shift}, $\|f(\cdot+\zeta u)-f\|_2^2=\sum_ka_k^2\,4\sin^2(\pi k\zeta u)$.
Since $4\sin^2(\pi k\zeta u)\le\min\{4,(2\pi k\zeta u)^2\}$ and $\kappa_0<2Q+2\le4Q$,
\begin{align*}
\frac1{\kappa_0}\int_0^{\kappa_0}\|f(\cdot+\zeta u)-f\|_2^2\,du&\le\sum_{k=1}^{\bar k}a_k^2\min\Big\{4,\frac{(2\pi k\zeta\kappa_0)^2}3\Big\}\\
&\le\frac{64\pi^2}{3}\sum_{k=1}^{\bar k}a_k^2\min\big\{1,(Qk|\zeta|)^2\big\}\le2048\pi^2\epsilon^2,
\end{align*}
by the first bound in~\eqref{eq:proof_smaller_sieve_shifts}.
Averaging the inequality for $\|f-\tilde f\|_2^2$ over $u\in[0,\kappa_0]$, and using $\|f-f_0\|_2^2\le2\|f-\tilde f\|_2^2+2r^2$, gives~\eqref{eq:localisation_f}.
\qed

\paragraph{Proof of Theorem~\ref{theo:checking_conditions_2}:}
The conditions of Lemma~\ref{lemma:contraction_general} hold for the sets $S_n=\Xi_n$, with the same constants as in the proof of Theorem~\ref{theo:checking_conditions}.
Indeed, since $\Xi_n\subseteq\Theta_n$, \eqref{eq:condition_entropy} holds; \eqref{eq:condition_prior_mass} does not involve the sets; and, by the bound on $\Pi(\Theta_n^c)$ in that proof and by~\eqref{eq:cond_Z_2},
\[
\Pi\big(\Xi_n^c\big)\le\Pi\big(\Theta_n^c\big)+\sum_{k=1}^{\bar k_n}\Pi_K(K=k)\,\Pi_{\bm{Z}}\big(\bm{Z}\notin\Xi(h_n)\mid K=k\big)\lesssim\exp\big(-5n\epsilon_n^2\big).
\]
Hence, for all large enough $L$, $\Pi\big(\|\theta-\theta_0\|_{2,n}>L\epsilon_n\mid X^{(n)}\big)\to0$ and $\Pi\big(\Xi_n\mid X^{(n)}\big)\to1$, in $P_{\theta_0}^{(n)}$-probability.

Fix such an $L$.
For large $n$, the conditions~\eqref{eq:new_conditions} hold with $\bar k=\bar k_n$, $l=l_n$, $h=h_n$ and $\epsilon=L\epsilon_n$: $\bar k_n=o(n)$ and $l_n=o(n)$ by~\eqref{eq:cond}, $\bar k_n\ge k_n\to\infty$ by the last condition in~\eqref{eq:cond}, $\epsilon_n^2=o(h_n)$ by~\eqref{eq:cond_h}, and $\epsilon_n\to0$.
For a parameter in $\Xi_n$, $f=\bm Z^T\bm\varphi$ is a trigonometric polynomial of degree at most $\bar k_n$ with coordinates in $\Xi(h_n)$, and $\tau\ge l_n^{-1}$.
By Proposition~\ref{prop:smaller_sieve}, every parameter in $\Xi_n$ with $\|\psi(\tau,f)-\theta_0\|_{2,n}\le L\epsilon_n$ then satisfies $|\tau-\tau_0|\le56L\epsilon_n/(\sigma_0n)$ and
\[
\|f-f_0\|_2^2\le C\bigg\{L^2\epsilon_n^2+\frac{\epsilon_n^2}{36}+\frac{(1+\log\bar k_n)\big(\|f_0\|_2^2+L^2\epsilon_n^2\big)}{n\tau_0}\bigg\}\lesssim\epsilon_n^2,
\]
since $\|f_0-f_{0,[\bar k_n]}\|_2\le\|f_0-f_{0,[k_n]}\|_2\le\epsilon_n/6$, as in the proof of Theorem~\ref{theo:checking_conditions}, and $\log\bar k_n\lesssim\log n=o(n\epsilon_n^2)$.
The posterior probability of the complement of this set of parameters is at most $\Pi\big(\Xi_n^c\mid X^{(n)}\big)+\Pi\big(\|\theta-\theta_0\|_{2,n}>L\epsilon_n\mid X^{(n)}\big)\to0$, which proves the theorem, with a larger constant in place of $L$.
\qed

\section{Proofs for Section~\ref{sec:examples}}\label{sec:technical:examples}

We use the following tail bound for the chi-square distribution~$\chi_d^2$ of Section~\ref{sec:preliminaries}: for $d\in\mathbb{N}$, $N\sim\chi_d^2$ and $t>0$,
\begin{equation}\label{eq:chi_square_tail}
\mathbb{P}\big(N>d+2(dt)^{1/2}+2t\big)\le e^{-t};
\end{equation}
cf.~\cite[Section~2.4]{Boucheron:2013}, where $N-d$ is shown to be sub-gamma on the right tail with variance factor $2d$ and scale parameter 2.

\begin{lemma}[Mass of the sets $\Xi(h)$]\label{lemma:xi_mass}
There is an absolute constant $c_*>0$ such that, for every $k\in\mathbb{N}$ and $h\in(0,1/4]$, a vector $\bm Z=(Z_{-k},\dots,Z_k)$ of independent standard normal variables satisfies $\mathbb{P}\big(\bm Z\in\Xi(h)\big)\ge c_*$.
\end{lemma}
\begin{proof}
By Lemma~\ref{lemma:xi_sets}(i), it suffices to consider $h=1/4$.
For $j=1,\dots,k$, let $a_j^2=Z_j^2+Z_{-j}^2$ be the squared amplitude of the $j$-th harmonic, let $S_1=\sum_{j=1}^ka_j^2$, which is positive almost surely, and, for $p\ge2$, let $S_p=\sum_{j\le k:\,p\mid j}a_j^2$, where $p\mid j$ means that $p$ divides $j$.
Then $\bm Z\in\Xi(1/4)$ if, and only if, $S_p\le3S_1/4$ for every $p=2,\dots,k$, since $S_p=0$ for $p>k$.
The variables $a_j^2/2$ are independent and exponentially distributed with mean 1.
Fix $2\le p\le k$, and let $m=\lfloor k/p\rfloor$, so that $1\le m\le k/2$.
Then $S_p/S_1$ is the ratio of a sum of $m$ of the variables $a_j^2/2$ to the sum of all $k$, so it has the beta distribution with parameters $m$ and $k-m$, which is the law of the $m$-th order statistic of $k-1$ independent uniform variables on $[0,1]$.
Hence $\mathbb{P}(S_p>3S_1/4)$ is the probability that at most $m-1$ of these uniform variables are at most $3/4$, and, by Hoeffding's inequality~\cite[Theorem~2.8]{Boucheron:2013},
\[
\mathbb{P}\big(S_p>3S_1/4\big)\le\mathbb{P}\Big(T\le\frac k2-1\Big)\le\exp\Big\{-\frac{(k+1)^2}{8(k-1)}\Big\}\le e^{-k/8},
\]
where $T$ has the binomial distribution with parameters $k-1$ and $3/4$, whose mean exceeds $k/2-1$ by $(k+1)/4$.
The union bound over $p$ gives $\mathbb{P}\big(\bm Z\notin\Xi(1/4)\big)\le(k-1)e^{-k/8}$, which is at most $1/2$ for $k\ge40$, a convenient round value (the bound holds from $k=34$ on); for small $k$ the union bound is useless, since $(k-1)e^{-k/8}$ exceeds 1, and we argue differently.
For $k<40$, the event $\big\{a_1^2>3\sum_{j=2}^ka_j^2\big\}$ has positive probability, and on it $S_p\le S_1-a_1^2<S_1/4$ for every $p\ge2$, so that $\mathbb{P}\big(\bm Z\in\Xi(1/4)\big)>0$.
The lemma follows with $c_*$ the minimum of $1/2$ and of $\mathbb{P}\big(\bm Z\in\Xi(1/4)\big)$ for $k=1,\dots,39$.
\qed\end{proof}

\paragraph{Proof of Proposition~\ref{prop:prior_example}:}
Since $\tilde k_n^{2\alpha+1}=n/\log n$, we have $n\epsilon_n^2=C_0^2n\tilde k_n^{-2\alpha}=C_0^2\tilde k_n\log n$, so that $\epsilon_n\to0$ and $n\epsilon_n^2\to\infty$.
For large $n$, $\tilde k_n\le k_n\le2\tilde k_n$, so that $k_n\log n\le2n\epsilon_n^2/C_0^2$.
In case a), let $D\ge\max\{1,10(2\alpha+1)C_0^2/b_2\}$ and $\bar k_n=\lceil D\tilde k_n\rceil$, and in case b), let $\bar k_n=k_n$ and $D=1$; in both cases $k_n\le\bar k_n$ and $\bar k_n\log n\le2Dn\epsilon_n^2/C_0^2$ for large $n$, so that $\bar k_n=o(n\epsilon_n^2)$.
The sequences $l_n$, $u_n$ and $m_n$ are at least 1 for large $n$, bounded by a polynomial in $n$, and $l_n=o(n)$, since $\epsilon_n\to0$.
We take $C_0^2\ge\max\{96,8b_1\}$ and $C_0\ge6M/(2\pi)^\alpha$, and check the conditions of Theorem~\ref{theo:checking_conditions} one by one; the lower bounds on $D$ and $C_0$ are those needed at the places where we refer to the choice of $D$ and of $C_0$, and are not optimised.

\emph{Condition~\eqref{eq:cond}}:
The first two conditions are shown above.
For the third, $(2k_n+1)^\alpha\epsilon_n\ge(2\tilde k_n)^\alpha C_0\tilde k_n^{-\alpha}=2^\alpha C_0\ge6M/\pi^\alpha=6C_\alpha$.

\emph{Conditions~\eqref{eq:cond_K_1} and~\eqref{eq:cond_K_2}}:
In case b), $\Pi_K(K=k_n)=1$ and $\Pi_K(K>\bar k_n)=0$.
In case a), by~\eqref{eq:prior_K_example} and since $k_n\le n$, for large $n$,
\[
\Pi_K(K=k_n)\ge\exp\big(-b_1k_n\log n\big)\ge\exp\big(-2b_1n\epsilon_n^2/C_0^2\big)\ge\exp\big(-n\epsilon_n^2/4\big).
\]
Also, $\log\bar k_n\ge\log\tilde k_n\ge\log n/\{2(2\alpha+1)\}$ and $\bar k_n^{-b_2}\le1/2$ for large $n$, so that, by~\eqref{eq:prior_K_example},
\begin{align*}
\Pi_K(K>\bar k_n)&\le\sum_{k>\bar k_n}e^{-b_2k\log\bar k_n}\le2e^{-b_2\bar k_n\log\bar k_n}\\
&\le2\exp\Big(-\frac{b_2D\tilde k_n\log n}{2(2\alpha+1)}\Big)\le2\exp\big(-5n\epsilon_n^2\big),
\end{align*}
by the choice of $D$.

\emph{Conditions~\eqref{eq:cond_Z_0} and~\eqref{eq:cond_Z_3}}:
Let $k\le\bar k_n$.
By Lemma~\ref{lemma:xi_mass}, conditioning on $\{\bm Z\in\Xi(h_n)\}$ multiplies the probability of every event by at most $1/c_*$, and, under the unconditioned law given $K=k$, the sum $\sum_{0<|i|\le k}Z_i^2$ has the $\chi_{2k}^2$ distribution, which is stochastically smaller than that of $N\sim\chi_{2\bar k_n}^2$.
For large $n$, $2\bar k_n\le5n\epsilon_n^2$, so that~\eqref{eq:chi_square_tail} with $d=2\bar k_n$ and $t=5n\epsilon_n^2$ gives
\begin{equation}\label{eq:proof_prior_example_chi}
\mathbb{P}\big(N>25n\epsilon_n^2\big)\le e^{-5n\epsilon_n^2},
\end{equation}
since $d+2(dt)^{1/2}+2t\le5t$.
Since $v_i\le2\pi\bar k_n$ for $|i|\le k$, the event $\bm Z\notin w_1(m_n)$ implies that $\sum_{0<|i|\le k}Z_i^2>m_n^2/(2\pi\bar k_n)^2=25n\epsilon_n^2$, so that $\Pi_{\bm Z}\big(\bm Z\notin w_1(m_n)\mid K=k\big)\le e^{-5n\epsilon_n^2}/c_*$.
Also, $Z_0$ is standard normal, so that $\Pi_{\bm Z}\big(|Z_0|>m_n\mid K=k\big)\le2e^{-m_n^2/2}\le2e^{-5n\epsilon_n^2}$, since $m_n^2\ge10n\epsilon_n^2$.

\emph{Condition~\eqref{eq:cond_Z_1}}:
Let $\bm z_n$ be the coordinates of $f_{0,[k_n]}$, as in Theorem~\ref{theo:checking_conditions}.
We first show that, for large $n$, the ball $\{\bm x\in\mathbb{R}^{2k_n+1}:\|\bm x-\bm z_n\|_2\le\epsilon_n/2\}$ is contained in $\Xi(h_n)$.
Let $g$ be the shape with coordinates $\bm x$ in this ball, and let $p\ge2$.
The maps $f\mapsto f-A_pf$ and $f\mapsto f-\bar f$ are orthogonal projections in $L^2\big([0,1]\big)$, and they commute with truncation, since the two maps and the truncation act on each coordinate separately.
By~\eqref{eq:truncation} and~\eqref{eq:cond}, $\|f_0-f_{0,[k_n]}\|_2\le\epsilon_n/6$, as in the proof of Theorem~\ref{theo:checking_conditions}, and, by~\eqref{eq:xi_set_functions} for $f_0\in\Xi(H)$,
\[
\|f_{0,[k_n]}-A_pf_{0,[k_n]}\|_2^2=\|f_0-A_pf_0\|_2^2-\big\|(f_0-f_{0,[k_n]})-A_p(f_0-f_{0,[k_n]})\big\|_2^2\ge H\sigma_0^2-\frac{\epsilon_n^2}{36}.
\]
Hence
\[
\|g-A_pg\|_2\ge\|f_{0,[k_n]}-A_pf_{0,[k_n]}\|_2-\frac{\epsilon_n}2\ge\Big(H\sigma_0^2-\frac{\epsilon_n^2}{36}\Big)^{1/2}-\frac{\epsilon_n}2.
\]
Also, $\|g-\bar g\|_2\le\|f_{0,[k_n]}-\bar f_0\|_2+\epsilon_n/2\le\sigma_0+\epsilon_n/2$.
Since $\epsilon_n\to0$ and $h_n\to0$, for large $n$ the right side of the previous display is positive and at least $h_n^{1/2}(\sigma_0+\epsilon_n/2)$, for every $p\ge2$, so that $g$ is not constant and $\bm x\in\Xi(h_n)$, by~\eqref{eq:xi_set_functions}.

Let $\delta_n=\epsilon_n/\{2(2k_n+1)^{1/2}\}$.
If $|Z_i-z_{n,i}|\le\delta_n$ for all $|i|\le k_n$, then $\|\bm Z-\bm z_n\|_2\le\epsilon_n/2$.
By the previous paragraph, the conditioned probability of the event $\{\|\bm Z-\bm z_n\|_2\le\epsilon_n/2\}$ is at least its unconditioned probability.
For a standard normal variable $\xi$, $z\in\mathbb{R}$ and $\delta>0$,
\[
\mathbb{P}\big(|\xi-z|\le\delta\big)\ge\frac{2\delta}{(2\pi)^{1/2}}e^{-(|z|+\delta)^2/2}\ge\frac{2\delta}{(2\pi)^{1/2}}e^{-z^2-\delta^2}.
\]
Hence, by independence under the unconditioned law,
\begin{align*}
\Pi_{\bm Z}\big(\|\bm Z-\bm z_n\|_2\le\epsilon_n/2\mid K=k_n\big)&\ge\prod_{|i|\le k_n}\frac{2\delta_n}{(2\pi)^{1/2}}e^{-z_{n,i}^2-\delta_n^2}\\
&=\exp\bigg\{-\|\bm z_n\|_2^2-\frac{\epsilon_n^2}4-(2k_n+1)\log\frac{(2\pi)^{1/2}}{2\delta_n}\bigg\}.
\end{align*}
Here $\|\bm z_n\|_2\le\|f_0\|_2$, $\epsilon_n\le1$ for large $n$, and $(2\pi)^{1/2}/(2\delta_n)=\{2\pi(2k_n+1)\}^{1/2}/\epsilon_n\le n^2$ for large $n$, since $k_n\le n$ and $\epsilon_n\ge n^{-1/2}$.
Since $2k_n+1\le3k_n$, the last term in the exponent is at most $6k_n\log n\le12n\epsilon_n^2/C_0^2\le n\epsilon_n^2/8$, by the choice of $C_0$.
Hence the left side is at least $\exp\big(-\|f_0\|_2^2-1-n\epsilon_n^2/8\big)\gtrsim\exp\big(-n\epsilon_n^2/4\big)$.

\emph{Conditions on $\Pi_\tau$}:
The density of $\Pi_\tau$ is continuous and positive, hence bounded away from zero on a neighbourhood of $\tau_0$, and~\eqref{eq:cond_tau_2} holds by the discussion after Theorem~\ref{theo:checking_conditions}.

This completes the check of the conditions of Theorem~\ref{theo:checking_conditions}.
For Theorem~\ref{theo:checking_conditions_2}, condition~\eqref{eq:cond_h} holds by assumption, and the left side of~\eqref{eq:cond_Z_2} is zero, since $\Pi_{\bm Z}\big(\bm Z\notin\Xi(h_n)\mid K=k\big)=0$ for every $k$.
The last assertion of the proposition follows from Theorems~\ref{theo:checking_conditions} and~\ref{theo:checking_conditions_2}.
\qed

\begin{lemma}[Minimax lower bound]\label{lemma:minimax}
Let $\alpha>0$, $M>0$ and $\tau_0>0$, and write $E_\theta$ for the expectation under $P_\theta^{(n)}$.
There is a constant $c>0$, depending only on $\alpha$ and $M$, such that, for all large enough $n$,
\begin{align*}
\inf_{\hat f_n}\sup_{f\in w_\alpha(M)}E_{\psi(\tau_0,f)}\|\hat f_n-f\|_2^2&\ge cn^{-2\alpha/(2\alpha+1)},\\
\inf_{\hat\theta_n}\sup_{f\in w_\alpha(M)}E_{\psi(\tau_0,f)}\|\hat\theta_n-\psi(\tau_0,f)\|_{2,n}^2&\ge cn^{-2\alpha/(2\alpha+1)},
\end{align*}
where the infima are over all estimators based on $X^{(n)}$ with values in $L^2\big([0,1]\big)$ and in $L^2\big([0,n]\big)$, respectively.
\end{lemma}
\begin{proof}
The first bound is, up to the constant, the lower bound of~\cite[Theorem~3]{Castillo:2006}, which is proved there at a fixed period, with Pinsker's constant; that theorem is stated for Sobolev smoothness at least 2, so we give the short argument, which covers every $\alpha>0$.
Since the estimators in the infima may depend on $\tau_0$, the bounds hold a fortiori for estimators that do not.
Let $n$ be large enough that $n\tau_0\ge2$, and let $n'=\lceil n\tau_0\rceil/\tau_0$, so that $n\le n'\le n+1/\tau_0\le2n$.
Observing $Y$ on $[0,n']$ instead of $[0,n]$ can only decrease the minimax risks, so it suffices to prove the first bound for estimators based on $(Y_t:t\in[0,n'])$.
On $[0,n']$, the process
\[
\bar Y(u)=\frac1{n'}\sum_{j=0}^{n'\tau_0-1}\big(Y_{(j+u)/\tau_0}-Y_{j/\tau_0}\big),\quad u\in[0,1],
\]
satisfies $d\bar Y(u)=f(u)\,du+(n')^{-1/2}d\bar W(u)$ for a standard Brownian motion $\bar W$.
It is sufficient for $f$, since, in~\eqref{eq:likelihood} with $n'$ in place of $n$, $\int_0^{n'}f(\tau_0t)\,dY_t=n'\int_0^1f(u)\,d\bar Y(u)$ and $\int_0^{n'}\{f(\tau_0t)\}^2\,dt=n'\|f\|_2^2$.
The squared error $\|\hat f-f\|_2^2$ is the sum of the squared errors of the coordinates of $\hat f$, and replacing each coordinate by its conditional expectation given $\bar Y$ does not increase its squared error, by the Rao-Blackwell theorem~\cite[Theorem~1.7.8]{Lehmann:1998}.
Hence the minimax risk for estimating $f\in w_\alpha(M)$ in $\|\cdot\|_2^2$ from $(Y_t:t\in[0,n'])$ is at least that of the Gaussian white noise model on $[0,1]$ with noise level $(n')^{-1/2}$.
By~\cite[Theorem~3.1]{Tsybakov:2009}, where the class $\tilde W(\alpha,M)$ is $w_\alpha(M)$, the latter is at least a constant multiple of $(n')^{-2\alpha/(2\alpha+1)}\ge(2n)^{-2\alpha/(2\alpha+1)}$ for large $n$, with a constant that depends only on $\alpha$ and $M$.
This gives the first bound.

For the second, let $V$ be the set of restrictions to $[0,n]$ of the signals $\psi(\tau_0,g)$, for shapes $g$.
By the first bound in~\eqref{eq:information_1}, $\|\psi(\tau_0,g)\|_{2,n}^2\ge\{1-1/(n\tau_0)\}\|g\|_2^2\ge\|g\|_2^2/2$ for every shape $g$, so $V$ is the image of $L^2\big([0,1]\big)$ under a linear map that is bounded and bounded below, and is therefore a closed subspace of $L^2\big([0,n]\big)$.
Let $\hat\theta_n$ be an estimator of the signal, and let $\psi(\tau_0,\hat g_n)$ be the orthogonal projection of $\hat\theta_n$ onto $V$, which defines an estimator $\hat g_n$ of the shape.
For $f\in w_\alpha(M)$, $\psi(\tau_0,f)\in V$, so that
\[
\|\hat\theta_n-\psi(\tau_0,f)\|_{2,n}^2\ge\|\psi(\tau_0,\hat g_n-f)\|_{2,n}^2\ge\frac12\|\hat g_n-f\|_2^2,
\]
and the second bound follows from the first, with $c/2$ in place of $c$, so both bounds hold with $c/2$.
\qed\end{proof}

\paragraph{Gaussian random series priors.}

\begin{proposition}\label{prop:gaussian_fails}
Let $f_0$ be a shape with a non-zero $p$-th harmonic for some integer $p\ge2$.
Let $k_n\to\infty$ and $\epsilon_n\to0$ with $n\epsilon_n^2\to\infty$, and let $h_n\in(0,1]$ with $\epsilon_n^2=o(h_n)$.
Assume that, given $K=k_n$, the coordinates $Z_i$, $|i|\le k_n$, are independent, centred and normally distributed, and that~\eqref{eq:cond_K_1} and~\eqref{eq:cond_Z_1} hold.
Then
\[
\Pi_K(K=k_n)\,\Pi_{\bm{Z}}\big(\bm{Z}\notin\Xi(h_n)\mid K=k_n\big)\gtrsim\exp\big(-n\epsilon_n^2/2\big),
\]
so that~\eqref{eq:cond_Z_2} does not hold, for any sequence $\bar k_n\ge k_n$.
\end{proposition}

Among the shapes in $\mathcal{F}$, the assumption on $f_0$ excludes only those whose single non-zero harmonic is the first.
The lower bound in Proposition~\ref{prop:gaussian_fails} is of the same order as the lower bound on the prior mass near $\theta_0$ given by~\eqref{eq:cond_K_1} and~\eqref{eq:cond_Z_1}, so changing the constant in~\eqref{eq:cond_Z_2} does not help.
The reason is that a Gaussian prior with independent centred coordinates puts at least as much mass near the shape $A_pf_{0,[k_n]}$, which has period $1/p$, as near $f_{0,[k_n]}$ (see the proof); conditioning on $\Xi(h_n)$, as in Section~\ref{sec:examples}, removes this mass.
Proposition~\ref{prop:gaussian_fails} concerns the conditions of Theorem~\ref{theo:checking_conditions_2}; we do not know whether the marginal posterior for the frequency is inconsistent under Gaussian priors.

\paragraph{Proof of Proposition~\ref{prop:gaussian_fails}:}
Let $a_p>0$ be the amplitude of the $p$-th harmonic of $f_0$, and let $n$ be large enough that $k_n\ge p$.
Given $K=k_n$, split $\bm Z$ into $\bm Z'=(Z_i:|i|\le k_n,\,p\nmid i)$ and $\bm Z''=(Z_i:|i|\le k_n,\,p\mid i)$, and $\bm z_n$, the coordinates of $f_{0,[k_n]}$, into $\bm z_n'$ and $\bm z_n''$ in the same way.
The vectors $\bm Z'$ and $\bm Z''$ are independent and centred normal.
By Anderson's lemma~\cite[Theorem~2.4.4]{Gine:2016}, applied to $\bm Z'$ and the ball of radius $\epsilon_n/2$ centred at the origin,
\[
\Pi_{\bm Z}\big(\|\bm Z'\|_2\le\epsilon_n/2\mid K=k_n\big)\ge\Pi_{\bm Z}\big(\|\bm Z'-\bm z_n'\|_2\le\epsilon_n/2\mid K=k_n\big).
\]
Let $E$ be the event that $\|\bm Z'\|_2\le\epsilon_n/2$ and $\|\bm Z''-\bm z_n''\|_2\le\epsilon_n/2$.
Since $\|\bm Z-\bm z_n\|_2^2=\|\bm Z'-\bm z_n'\|_2^2+\|\bm Z''-\bm z_n''\|_2^2$, independence, the previous display and~\eqref{eq:cond_Z_1} give
\begin{align*}
\Pi_{\bm Z}\big(E\mid K=k_n\big)&\ge\Pi_{\bm Z}\big(\|\bm Z'-\bm z_n'\|_2\le\epsilon_n/2,\;\|\bm Z''-\bm z_n''\|_2\le\epsilon_n/2\mid K=k_n\big)\\
&\ge\Pi_{\bm Z}\big(\|\bm Z-\bm z_n\|_2\le\epsilon_n/2\mid K=k_n\big)\gtrsim e^{-n\epsilon_n^2/4}.
\end{align*}
The coordinates of $\bm z_n''$ with indices $\pm p$ have squared sum $a_p^2$, so on $E$, by the triangle inequality, $\sum_{i\in\mathbb{Z}^*}Z_{pi}^2\ge(a_p-\epsilon_n/2)^2\ge a_p^2/4$ for large $n$, while $\sum_{p\nmid i}Z_i^2\le\epsilon_n^2/4$.
Since $\epsilon_n^2=o(h_n)$, for large $n$ we have $\epsilon_n^2/4<h_na_p^2/4$, so that, on $E$,
\[
\sum_{i\in\mathbb{Z}^*}Z_{pi}^2=\sum_{i\in\mathbb{Z}^*}Z_i^2-\sum_{p\nmid i}Z_i^2>(1-h_n)\sum_{i\in\mathbb{Z}^*}Z_i^2,
\]
i.e., $\bm Z\notin\Xi(h_n)$.
By~\eqref{eq:cond_K_1}, the left side in the statement is at least $\Pi_K(K=k_n)\,\Pi_{\bm Z}\big(E\mid K=k_n\big)\gtrsim e^{-n\epsilon_n^2/2}$.
Since $n\epsilon_n^2\to\infty$, this is not bounded by a multiple of $e^{-5n\epsilon_n^2}$, and it is one of the terms of the sum in~\eqref{eq:cond_Z_2} when $\bar k_n\ge k_n$.
\qed

\section{Proofs for Section~\ref{sec:BVM}}\label{sec:technical:BVM}

In this appendix, for a function $h$ that is square integrable on $[0,n]$, we write $|h|_n=n^{1/2}\|h\|_{2,n}=\big(\int_0^n\{h(t)\}^2\,dt\big)^{1/2}$, and $\langle h_1,h_2\rangle_n=\int_0^nh_1(t)h_2(t)\,dt$.
For $\tau\in\mathbb{R}$ and a shape $g$, we write $g_\tau$ for the function $t\mapsto g\big(\tau(t-n/2)\big)$ on $[0,n]$, so that $\psi_{n/2}(\tau,g)=g_\tau$ on $[0,n]$, and $(g')_\tau$ for the corresponding function of the derivative $g'$.
Since $g_\tau=\tilde g(\tau\,\cdot)$ for the shape $\tilde g=g(\cdot-\tau n/2)$, which has the same norm as $g$, the first bound in~\eqref{eq:information_1}, which holds for all $n,\tau>0$, gives $|g_\tau|_n^2\le2n\|g\|_2^2$ if $n\tau\ge1$.
We use three consequences of this bound.
Let $\tau,\tau'>0$ with $n\min\{\tau,\tau'\}\ge1$, and let $g$ be a shape.
If $g\in\mathcal{W}_1$, then
\[
g_\tau(t)-g_{\tau'}(t)=(\tau-\tau')\Big(t-\frac n2\Big)\int_0^1(g')_{\tau'+r(\tau-\tau')}(t)\,dr,
\]
so that, by Jensen's inequality and since $|t-n/2|\le n/2$,
\begin{equation}\label{eq:proof_BvM_lipschitz}
|g_\tau-g_{\tau'}|_n\le\frac{|\tau-\tau'|n^{3/2}}{\sqrt2}\|g'\|_2.
\end{equation}
If $g\in\mathcal{W}_2$, then, in the same way,
\[
g_\tau(t)-g_{\tau'}(t)-(\tau-\tau')\Big(t-\frac n2\Big)(g')_{\tau'}(t)=(\tau-\tau')^2\Big(t-\frac n2\Big)^2\int_0^1\!\!\int_0^1r\,(g'')_{\tau'+rr'(\tau-\tau')}(t)\,dr'\,dr,
\]
and, since $\int_0^1\int_0^1r\,dr'\,dr=1/2$,
\begin{equation}\label{eq:proof_BvM_second_order}
\Big|g_\tau-g_{\tau'}-(\tau-\tau')\Big(\cdot-\frac n2\Big)(g')_{\tau'}\Big|_n\le\frac{(\tau-\tau')^2n^{5/2}}{\sqrt{32}}\|g''\|_2.
\end{equation}
Finally, $|(\cdot-n/2)g_{\tau}|_n\le(n/2)|g_\tau|_n\le n^{3/2}\|g\|_2/\sqrt2$.

The stochastic integrals below have deterministic integrands.
When the integrand is continuously differentiable in $t$ and depends continuously on a parameter, we use, as in Section~\ref{sec:estimation}, the version given by the integration by parts formula, which is continuous in the parameter; identities between such integrals that hold almost surely for each value of the parameter then hold almost surely for all values simultaneously.

\paragraph{Proof of Proposition~\ref{prop:origin}:}
For each $K$, the map $(\tau,\bm Z)\mapsto\bm Z^{(\tau c)}$ rotates each pair $(Z_k,Z_{-k})$ through the angle $-2\pi k\tau c$ (Section~\ref{sec:preliminaries}), so it is continuous, and the map $(\tau,K,\bm Z)\mapsto(\tau,K,\bm Z^{(\tau c)})$ is measurable.
Let $S_1\subseteq(0,\infty)$ and $S_2$ be measurable sets of frequencies and of parameters $(K,\bm Z)$, respectively.
Since $\tau$ and $(K,\bm Z)$ are independent under $\Pi$, Fubini's theorem and the stationarity of $\Pi_f$ give
\[
\Pi\Big(\tau\in S_1,\,(K,\bm Z^{(\tau c)})\in S_2\Big)=\int_{S_1}\Pi_f\big((K,\bm Z^{(\tau'c)})\in S_2\big)\,d\Pi_\tau(\tau')=\Pi_\tau(S_1)\,\Pi_f(S_2).
\]
Since such product sets determine a distribution, $\big(\tau,f(\cdot+\tau c)\big)$ has distribution $\Pi_\tau\times\Pi_f=\Pi$.
Since $\psi(\tau,f)=\psi_c\big(\tau,f(\cdot+\tau c)\big)$, by~\eqref{eq:bayess_formula},
\[
\Pi\Big(\big(\tau,f(\cdot+\tau c)\big)\in S\;\Big|\;X^{(n)}\Big)=\frac
{\int1_S\big(\tau,f(\cdot+\tau c)\big)\,p_{\psi_c(\tau,f(\cdot+\tau c))}^{(n)}\big(X^{(n)}\big)\,d\Pi(\tau,f)}
{\int p_{\psi_c(\tau,f(\cdot+\tau c))}^{(n)}\big(X^{(n)}\big)\,d\Pi(\tau,f)},
\]
and the change of variables $(\tau,f)\mapsto\big(\tau,f(\cdot+\tau c)\big)$, which leaves $\Pi$ unchanged, turns the right side into the right side in the statement.
The map does not change $\tau$, which gives the last assertion.
\qed

\paragraph{Proof of Proposition~\ref{prop:LAN}:}
We use the notation $|\cdot|_n$, $\langle\cdot,\cdot\rangle_n$ and $g_\tau$ introduced at the start of this appendix.
Let $\delta=\tau_n-\tau_0=sn^{-3/2}$, and let $n$ be large enough that $n\tau_0\ge2$ and $n\tau_n\ge1$.
Since $\theta_0=\psi_{n/2}(\tau_0,\tilde f_{0,n})$, the signal $\theta_n=\psi_{n/2}(\tau_n,f_n)$ satisfies $\theta_n-\theta_0=\lambda+\gamma+\zeta$ on $[0,n]$, where
\begin{align*}
\lambda(t)&=\delta\Big(t-\frac n2\Big)f_0'(\tau_0t),
\qquad
\gamma=n^{-1/2}g_{\tau_0},\\
\zeta&=(\tilde f_{0,n})_{\tau_n}-(\tilde f_{0,n})_{\tau_0}-\lambda+n^{-1/2}\big(g_{\tau_n}-g_{\tau_0}\big);
\end{align*}
note that $\lambda=\delta(\cdot-n/2)(\tilde f_{0,n}')_{\tau_0}$.
Then $\int_0^n\lambda\,dW=s\Delta_n$, $\int_0^n\gamma\,dW=\mathbb{W}_n(g)$, $|\lambda|_n^2=s^2\mathcal{I}_n^{(n/2)}(\tau_0,f_0)/n^3$ and $|\gamma|_n^2=n^{-1}|g_{\tau_0}|_n^2$, and, by~\eqref{eq:log_lik_ratio},
\begin{align*}
\log\frac{p_{\theta_n}^{(n)}}{p_{\theta_0}^{(n)}}\big(X^{(n)}\big)
&=s\Delta_n+\mathbb{W}_n(g)-\frac12\big(|\lambda|_n^2+|\gamma|_n^2\big)-\langle\lambda,\gamma\rangle_n\\
&\quad+\int_0^n\zeta\,dW-\langle\zeta,\lambda+\gamma\rangle_n-\frac12|\zeta|_n^2.
\end{align*}
By~\eqref{eq:information_2}, $|\lambda|_n^2\to s^2\mathcal{I}_0$, and, by the first bound in~\eqref{eq:information_1}, $|\gamma|_n^2\to\|g\|_2^2$.
The function $g_{\tau_0}$ is $\tilde g(\tau_0\,\cdot)$ for a shape $\tilde g$ with $\|\tilde g\|_2=\|g\|_2$, so, by~\eqref{eq:information_3}, $|\langle\lambda,\gamma\rangle_n|\le|s|\,\|f_0'\|_2\|g\|_2/(2n\tau_0)$, which tends to zero.
Since $\int_0^n\zeta\,dW$ is normally distributed with mean zero and variance $|\zeta|_n^2$, and $|\lambda|_n$ and $|\gamma|_n$ are bounded, it remains to show that $|\zeta|_n\to0$.

Let $k\in\mathbb{N}$, and split $\tilde f_{0,n}=p_k+e_k$, with $p_k=(\tilde f_{0,n})_{[k]}$, and $g=g_{[k]}+(g-g_{[k]})$.
By~\eqref{eq:proof_BvM_second_order} applied to $p_k$, \eqref{eq:proof_BvM_lipschitz} applied to $e_k$ and to $g_{[k]}$, the bound on $|(\cdot-n/2)g_\tau|_n$ applied to $e_k'$, and the bound $|g_\tau|_n^2\le2n\|g\|_2^2$ applied to $g-g_{[k]}$,
\[
|\zeta|_n\le\frac{s^2\|p_k''\|_2}{\sqrt{32}\,n^{1/2}}+\sqrt2\,|s|\,\|e_k'\|_2+\frac{|s|\,\|g_{[k]}'\|_2}{\sqrt2\,n^{1/2}}+2\sqrt2\,\|g-g_{[k]}\|_2.
\]
Here $\|p_k''\|_2\le2\pi k\|f_0'\|_2$, $\|e_k'\|_2=\|(f_0-f_{0,[k]})'\|_2$ and $\|g_{[k]}'\|_2\le2\pi k\|g\|_2$, by Section~\ref{sec:preliminaries}.
Hence $\limsup_{n\to\infty}|\zeta|_n\le\sqrt2\,|s|\,\|(f_0-f_{0,[k]})'\|_2+2\sqrt2\,\|g-g_{[k]}\|_2$, which tends to zero as $k\to\infty$, since $f_0'$ and $g$ are square integrable.

Finally, $\big(\Delta_n,\mathbb{W}_n(g)\big)$ is a centred normal vector with variances $\mathcal{I}_n^{(n/2)}(\tau_0,f_0)/n^3$ and $|\gamma|_n^2$ and covariance
\[
\frac1{n^2}\int_0^n\Big(t-\frac n2\Big)f_0'(\tau_0t)g_{\tau_0}(t)\,dt,
\]
which is at most $\|f_0'\|_2\|g\|_2/(2n\tau_0)$ in absolute value, by~\eqref{eq:information_3}.
The variances converge to $\mathcal{I}_0$ and $\|g\|_2^2$ and the covariance to zero, which gives the convergence in distribution.
\qed

The proof of Theorem~\ref{theo:BvM} uses a bound on the supremum of the following Gaussian processes.
Bounds on the supremum over the period of periodogram-type processes are in~\cite[Lemma~6]{Castillo:2007} and~\cite[Lemma~4]{LevyLeduc:2008}; the processes here are different, and their supremum is bounded by Dudley's entropy bound.
For $\tau_0>0$ and $s\in\mathbb{R}$, write $\tau_{s}=\tau_0+sn^{-3/2}$, and, for $i\in\mathbb{Z}^*$, let
\begin{align*}
\mathbb{G}_i(s)&=\frac1{n^{3/2}}\int_0^n\Big(t-\frac n2\Big)\int_0^1\frac{(\varphi_i')_{\tau_{rs}}(t)}{v_i}\,dr\,dW_t,\\
\mathbb{H}_i(s)&=\frac1{n^{5/2}}\int_0^n\Big(t-\frac n2\Big)^2\int_0^1\int_0^1r\,\frac{(\varphi_i'')_{\tau_{rr's}}(t)}{v_i^2}\,dr'\,dr\,dW_t.
\end{align*}

\begin{lemma}[Maximal inequality]\label{lemma:maximal}
Let $\tau_0>0$, $n\ge1$, $k\in\mathbb{N}$ and $\rho\ge1$.
Then, for $0<|i|\le k$ and $|s|\le\rho$,
\[
\int_0^n\big\{(\varphi_i)_{\tau_s}-(\varphi_i)_{\tau_0}\big\}\,dW=sv_i\,\mathbb{G}_i(s)
\quad\text{and}\quad
\mathbb{G}_i(s)-\mathbb{G}_i(0)=\frac{sv_i}{n^{1/2}}\,\mathbb{H}_i(s),
\]
almost surely, and, for an absolute constant $C$,
\[
\mathbb{E}\Big[\max_{0<|i|\le k}\sup_{|s|\le\rho}\big\{|\mathbb{G}_i(s)|+|\mathbb{H}_i(s)|\big\}\Big]\le C\big\{\log(k\rho)+1\big\}^{1/2}.
\]
\end{lemma}
\begin{proof}
For fixed $t$, the fundamental theorem of calculus gives
\begin{align*}
(\varphi_i)_{\tau_s}(t)-(\varphi_i)_{\tau_0}(t)&=sn^{-3/2}\Big(t-\frac n2\Big)\int_0^1(\varphi_i')_{\tau_{rs}}(t)\,dr,\\
(\varphi_i')_{\tau_{rs}}(t)-(\varphi_i')_{\tau_0}(t)&=rsn^{-3/2}\Big(t-\frac n2\Big)\int_0^1(\varphi_i'')_{\tau_{rr's}}(t)\,dr',
\end{align*}
and the two identities follow by linearity of the stochastic integral, since $\mathbb{G}_i(0)=n^{-3/2}\int_0^n(t-n/2)(\varphi_i')_{\tau_0}(t)\,dW_t/v_i$.

The functions $\varphi_i'/v_i$ and $\varphi_i''/v_i^2$ are bounded by $\sqrt2$ and Lipschitz with constant $\sqrt2\,v_i$, and $|\tau_{rs}-\tau_{rs'}|\le|s-s'|n^{-3/2}$ for $r\in[0,1]$.
Hence, since $|t-n/2|\le n/2$ for $t\in[0,n]$, the variances of $\mathbb{G}_i(s)$ and $\mathbb{H}_i(s)$ are at most 1, and
\[
\mathbb{E}\big\{\mathbb{G}_i(s)-\mathbb{G}_i(s')\big\}^2+\mathbb{E}\big\{\mathbb{H}_i(s)-\mathbb{H}_i(s')\big\}^2\le\frac{v_i^2(s-s')^2}{n}.
\]
These crude bounds, a uniform bound on the variances and a Lipschitz bound in $s$, are all that the entropy bound below needs.
Let $\mathcal{U}=\{1,2\}\times\{i:0<|i|\le k\}\times[-\rho,\rho]$, and let $\mathbb{X}(1,i,s)=\mathbb{G}_i(s)$ and $\mathbb{X}(2,i,s)=\mathbb{H}_i(s)$.
This is a centred Gaussian process on $\mathcal{U}$.
Let $d$ be the standard deviation of its increments.
Then $d\le2$ on $\mathcal{U}\times\mathcal{U}$, and $d\big((j,i,s),(j,i,s')\big)\le2\pi k|s-s'|$, since $v_i\le2\pi k$ and $n\ge1$.
Covering $[-\rho,\rho]$ by intervals of length $2\epsilon/(2\pi k)$ gives $N(\mathcal{U},d,\epsilon)\le4k(1+2\pi k\rho/\epsilon)$ for $\epsilon\in(0,1]$, and the diameter of $\mathcal{U}$ is at most 2.
By~\cite[Theorem~2.3.6]{Gine:2016}, applied to Gaussian processes, which satisfy its sub-Gaussian condition with respect to $d$, for every finite set $\mathcal{S}\subseteq\mathcal{U}$,
\[
\mathbb{E}\Big[\max_{x\in\mathcal{S}}|\mathbb{X}(x)|\Big]\le1+4\sqrt2\int_0^1\big\{\log\big(8k(1+2\pi k\rho/\epsilon)\big)\big\}^{1/2}\,d\epsilon\le C\big\{\log(k\rho)+1\big\}^{1/2}.
\]
The paths of $\mathbb{X}(j,i,\cdot)$ are continuous, so the supremum over $\mathcal{U}$ is the supremum over a countable set, and the bound extends to it by monotone convergence.
The lemma follows, since $|\mathbb{G}_i(s)|+|\mathbb{H}_i(s)|\le2\sup_{\mathcal{U}}|\mathbb{X}|$.
\qed\end{proof}

\begin{lemma}[Bernstein-von Mises theorem from a uniform expansion]\label{lemma:BvM_general}
Let $\tau_0>0$, and let $\Pi=\Pi_\tau\times\Pi_f$ be a prior on parameters $(\tau,f)$, where $\Pi_\tau$ has a Lebesgue density $\pi_\tau$ that is continuous and positive at $\tau_0$.
Let the posterior be
\[
\Pi\big(S\mid X^{(n)}\big)=\frac{\int_S\exp\big\{\ell_n(\tau,f)\big\}\,d\Pi(\tau,f)}{\int\exp\big\{\ell_n(\tau,f)\big\}\,d\Pi(\tau,f)},
\]
for a function $\ell_n$ that is jointly measurable in the parameters and the data and such that the denominator is positive and finite almost surely.
Let $r_n,\rho_n>0$ with $\rho_n\to\infty$ and $\rho_n/r_n\to0$, let
\[
\mathcal{T}_n=\big\{\tau:r_n|\tau-\tau_0|\le\rho_n\big\},
\]
and let $\mathcal{G}_n$ be measurable sets of parameters $f$ such that $\Pi(\mathcal{T}_n\times\mathcal{G}_n\mid X^{(n)})\to1$.
Assume that there are a constant $\mathcal{I}_0>0$ and random variables $\Delta_n$, bounded in probability, and $\eta_n\ge0$, with $\eta_n\to0$ in probability, such that, for all $\tau\in\mathcal{T}_n$ and $f\in\mathcal{G}_n$, with $s=r_n(\tau-\tau_0)$,
\begin{equation}\label{eq:BvM_general_expansion}
\Big|\ell_n(\tau,f)-\ell_n(\tau_0,f)-s\Delta_n+\frac{s^2\mathcal{I}_0}2\Big|\le\eta_n(1+s^2).
\end{equation}
Then the total variation distance between the posterior distribution of $r_n(\tau-\tau_0)$ and $N(\Delta_n/\mathcal{I}_0,1/\mathcal{I}_0)$ tends to zero, and, for every $a\in(0,1)$, the $a$-quantile $Q_n(a)$ of the marginal posterior for $\tau$ satisfies
\[
r_n\big(Q_n(a)-\tau_0\big)=\frac{\Delta_n}{\mathcal{I}_0}+\frac{z_a}{\mathcal{I}_0^{1/2}}+o_P(1),
\]
where $z_a$ is the $a$-quantile of the standard normal distribution.
All limits are in $P_{\theta_0}^{(n)}$-probability, as $n\to\infty$.
\end{lemma}
\begin{proof}
The proof follows that of~\cite[Theorem~1]{Castillo:2012}.
Let $\Pi_n(\cdot\mid X^{(n)})$ be the posterior conditioned on $\mathcal{T}_n\times\mathcal{G}_n$.
The total variation distance between $\Pi(\cdot\mid X^{(n)})$ and $\Pi_n(\cdot\mid X^{(n)})$ is at most the posterior probability of the complement of $\mathcal{T}_n\times\mathcal{G}_n$, so it suffices to prove the first assertion for $\Pi_n(\cdot\mid X^{(n)})$.
Write $\tau_s=\tau_0+s/r_n$.
By Fubini's theorem, under $\Pi_n(\cdot\mid X^{(n)})$ the variable $s=r_n(\tau-\tau_0)$ has the density proportional to
\begin{equation}\label{eq:proof_BvM_density}
1_{\{|s|\le\rho_n\}}\,\pi_\tau(\tau_s)\int_{\mathcal{G}_n}\exp\big\{\ell_n(\tau_s,f)\big\}\,d\Pi_f(f).
\end{equation}
Let
\[
\mathcal{M}_n=\int_{\mathcal{G}_n}\exp\big\{\ell_n(\tau_0,f)\big\}\,d\Pi_f(f),
\]
which is positive when the posterior probability of $\mathcal{T}_n\times\mathcal{G}_n$ is positive, and finite almost surely, by~\eqref{eq:BvM_general_expansion} and the finiteness of the denominator.
By~\eqref{eq:BvM_general_expansion}, the integral in~\eqref{eq:proof_BvM_density} lies between the values at $\vartheta=-1$ and $\vartheta=1$ of
\[
\mathcal{M}_n\exp\Big\{s\Delta_n-\frac{s^2\mathcal{I}_0}2+\vartheta\eta_n(1+s^2)\Big\}.
\]
Since $\mathcal{M}_n$ does not depend on $s$, the density of $s$ is $\varpi_n/\int\varpi_n$, where
\[
\varpi_n(s)=1_{\{|s|\le\rho_n\}}\frac{\pi_\tau(\tau_s)}{\pi_\tau(\tau_0)}\exp\Big\{s\Delta_n-\frac{s^2\mathcal{I}_0}2+\vartheta_n(s)\eta_n(1+s^2)\Big\},
\]
for some $\vartheta_n(s)\in[-1,1]$.
Let
\[
\omega_n=\sup_{\tau\in\mathcal{T}_n}\bigg|\frac{\pi_\tau(\tau)}{\pi_\tau(\tau_0)}-1\bigg|,
\]
which tends to zero, since $\pi_\tau$ is continuous and positive at $\tau_0$ and $\rho_n/r_n\to0$.
Write $\varpi(s)=\exp(s\Delta_n-s^2\mathcal{I}_0/2)$.
If $x\ge0$ and $\omega\in[0,1]$, then every number between $(1-\omega)e^{-x}$ and $(1+\omega)e^{x}$ is at distance at most $(1+\omega)e^{x}-1$ from 1, since
\[
1-(1-\omega)e^{-x}\le\omega+x\le(1+\omega)e^{x}-1.
\]
Hence, if $\omega_n\le1$,
\[
\int_{\mathbb{R}}|\varpi_n(s)-\varpi(s)|\,ds\le\int_{\mathbb{R}}\varpi(s)\big\{(1+\omega_n)e^{\eta_n(1+s^2)}-1\big\}\,ds+\int_{|s|>\rho_n}\varpi(s)\,ds.
\]
Let $\xi_n$ be a variable that, given $X^{(n)}$, has the distribution $N(\mu_n,1/\mathcal{I}_0)$, with $\mu_n=\Delta_n/\mathcal{I}_0$.
The right side of the last display divided by $\int\varpi(s)\,ds$ is
\[
\kappa_n=(1+\omega_n)e^{\eta_n}\,\mathbb{E}\big[e^{\eta_n\xi_n^2}\mid X^{(n)}\big]-1+\mathbb{P}\big(|\xi_n|>\rho_n\mid X^{(n)}\big),
\]
and, if $2\eta_n<\mathcal{I}_0$,
\[
\mathbb{E}\big[e^{\eta_n\xi_n^2}\mid X^{(n)}\big]=\Big(1-\frac{2\eta_n}{\mathcal{I}_0}\Big)^{-1/2}\exp\bigg(\frac{\eta_n\mu_n^2}{1-2\eta_n/\mathcal{I}_0}\bigg).
\]
Since $\eta_n\to0$ in probability, $\mu_n$ is bounded in probability and $\rho_n\to\infty$, we have $\kappa_n\to0$ in probability.
For non-negative integrable functions $w_1$ and $w_2$ with positive integrals,
\[
\int\bigg|\frac{w_1}{\int w_1}-\frac{w_2}{\int w_2}\bigg|\le\frac{2\int|w_1-w_2|}{\int w_2},
\]
so the total variation distance between the distribution of $s$ under $\Pi_n(\cdot\mid X^{(n)})$ and $N(\mu_n,1/\mathcal{I}_0)$ is at most $\kappa_n$, which proves the first assertion.

For the second, let $F_n$ be the distribution function of $s$ under the posterior.
By the first assertion, $\sup_x|F_n(x)-\Phi(\mathcal{I}_0^{1/2}(x-\mu_n))|\to0$ in probability.
Let $y>0$.
Since $\Phi(z_a+y)>a>\Phi(z_a-y)$, with probability tending to one
\[
F_n\Big(\mu_n+\frac{z_a+y}{\mathcal{I}_0^{1/2}}\Big)>a>F_n\Big(\mu_n+\frac{z_a-y}{\mathcal{I}_0^{1/2}}\Big),
\]
so that the $a$-quantile of $F_n$ is $\mu_n+z_a\mathcal{I}_0^{-1/2}+R_n(a)$, with $R_n(a)\to0$ in probability.
Since $s$ is an increasing function of $\tau$, $r_n(Q_n(a)-\tau_0)$ is the $a$-quantile of $F_n$.
\qed\end{proof}

\paragraph{Proof of Theorem~\ref{theo:BvM}:}
We verify the conditions of Lemma~\ref{lemma:BvM_general} in the parametrisation $\psi_{n/2}$, with $r_n=n^{3/2}$.
As in the proof of Proposition~\ref{prop:LAN}, $|\cdot|_n$, $\langle\cdot,\cdot\rangle_n$ and $g_\tau$ are as defined at the start of this appendix.
The main step is the expansion~\eqref{eq:BvM_general_expansion}.
In it, the cross term between the score for the frequency and the deviation of the shape from $\tilde f_{0,n}$ is small by~\eqref{eq:information_3}, and this is where the choice of the origin of time plays a role.
Condition~\eqref{eq:cond_BvM} controls the remainder terms, which involve the derivatives of the shapes, in deterministic terms and in a stochastic integral, uniformly over $|s|\lesssim n^{1/2}\epsilon_n$.

\emph{Reduction to the parametrisation $\psi_{n/2}$.}
By Proposition~\ref{prop:origin} with $c=n/2$, the posterior distribution of $\big(\tau,f(\cdot+\tau n/2)\big)$ is the posterior in the model parametrised by $\psi_{n/2}$, with prior $\Pi$, and it has the same marginal for $\tau$.
In this proof we work in this model, and we write $f$ for the shape in $\psi_{n/2}$, so that the true parameter is $(\tau_0,\tilde f_{0,n})$, and
\[
\ell_n(\tau,f)=\log\frac{p_{\psi_{n/2}(\tau,f)}^{(n)}}{p_{\theta_0}^{(n)}}\big(X^{(n)}\big).
\]
Let $L$ be large enough that the conclusions of Theorems~\ref{theo:checking_conditions} and~\ref{theo:checking_conditions_2} hold, let $L_1=L(1+\|f_0'\|_2/2)$ and $\rho_n=Ln^{1/2}\epsilon_n$, and let
\[
\mathcal{T}_n=\big\{\tau:|\tau-\tau_0|\le L\epsilon_n/n\big\},
\qquad
\mathcal{G}_n=\big\{(K,\bm Z):K\le\bar k_n,\;\|f-\tilde f_{0,n}\|_2\le L_1\epsilon_n\big\}.
\]
For $g\in\mathcal{W}_1$ and $c\in\mathbb{R}$, \eqref{eq:shift} and $4\sin^2(\pi kc)\le(2\pi kc)^2$ give $\|g(\cdot+c)-g\|_2^2\le c^2\sum_{k\ge1}(2\pi k)^2a_k^2=c^2\|g'\|_2^2$, where $a_k$ is the amplitude of the $k$-th harmonic of $g$.
Hence, by the triangle inequality, a parameter whose shape is $g$ in $\psi$ has shape $g(\cdot+\tau n/2)$ in $\psi_{n/2}$, with
\[
\big\|g(\cdot+\tau n/2)-\tilde f_{0,n}\big\|_2\le\|g-f_0\|_2+\frac n2|\tau-\tau_0|\,\|f_0'\|_2,
\]
which is at most $L_1\epsilon_n$ if $|\tau-\tau_0|\le L\epsilon_n/n$ and $\|g-f_0\|_2\le L\epsilon_n$.
A time shift does not change $K$, so, by Theorems~\ref{theo:checking_conditions} and~\ref{theo:checking_conditions_2}, the posterior probability of $\mathcal{T}_n\times\mathcal{G}_n$ tends to one in probability.
The posterior in this model is of the form in Lemma~\ref{lemma:BvM_general}, with this $\ell_n$, which is jointly measurable by the argument of Section~\ref{sec:estimation}, since the signals $\psi_{n/2}(\tau,\bm z^T\bm\varphi)$ are continuously differentiable in $t$ and depend continuously on $(\tau,\bm z)$, and $\mathcal{T}_n$ is the set of the lemma with $r_n=n^{3/2}$; moreover $\rho_n\to\infty$, since $n\epsilon_n^2\to\infty$, and $\rho_n/r_n=L\epsilon_n/n\to0$.
With $\Delta_n$ and $\mathcal{I}_0$ as in Section~\ref{sec:BVM:LAN}, Lemma~\ref{lemma:BvM_general} therefore gives~\eqref{eq:BvM_statement} and the expansion of the quantiles stated after Theorem~\ref{theo:BvM}, if there are random variables $\eta_n\ge0$ such that $\eta_n\to0$ in probability and, for all $|s|\le\rho_n$ and $(K,\bm Z)\in\mathcal{G}_n$,
\begin{equation}\label{eq:proof_BvM_expansion}
\Big|\ell_n(\tau_s,f)-\ell_n(\tau_0,f)-s\Delta_n+\frac{s^2\mathcal{I}_0}2\Big|\le\eta_n(1+s^2),
\end{equation}
with $\tau_s=\tau_0+sn^{-3/2}$, as in Lemma~\ref{lemma:maximal}.
The variable $\Delta_n$ is normal with mean zero and bounded variance, hence bounded in probability, and $\mathcal{I}_0>0$, since $f_0$ is not constant.
We prove~\eqref{eq:proof_BvM_expansion} below.

\emph{Credible intervals.}
Write $\mu_n=\Delta_n/\mathcal{I}_0$ and $n^{3/2}(Q_n(a)-\tau_0)=\mu_n+z_a\mathcal{I}_0^{-1/2}+R_n(a)$, so that $R_n(a)\to0$ in probability.
Then the event $Q_n(a/2)\le\tau_0\le Q_n(1-a/2)$ is the event
\[
z_{a/2}+\mathcal{I}_0^{1/2}R_n(a/2)\le-\mathcal{I}_0^{1/2}\mu_n\le z_{1-a/2}+\mathcal{I}_0^{1/2}R_n(1-a/2).
\]
The variable $-\mathcal{I}_0^{1/2}\mu_n=-\Delta_n/\mathcal{I}_0^{1/2}$ is normal with mean zero and variance converging to 1, so, by Slutsky's lemma~\cite[Lemma~2.8]{vanderVaart:1998} and the continuity of $\Phi$, the probability of this event tends to $\Phi(z_{1-a/2})-\Phi(z_{a/2})=1-a$.

\emph{Proof of~\eqref{eq:proof_BvM_expansion}.}
Let $|s|\le\rho_n$, $\delta=sn^{-3/2}$ and $(K,\bm Z)\in\mathcal{G}_n$, with shape $f$.
Let $n$ be large enough that $\rho_n\ge1$, $\rho_nn^{-3/2}\le\tau_0/2$ and $n\tau_0\ge2$, so that $n\tau_{s'}\ge1$ for $|s'|\le\rho_n$.
Let $\lambda=\delta(\cdot-n/2)(\tilde f_{0,n}')_{\tau_0}$, which is the function $\delta(t-n/2)f_0'(\tau_0t)$, so that $\int_0^n\lambda\,dW=s\Delta_n$ and $|\lambda|_n^2=s^2\mathcal{J}_n$, with $\mathcal{J}_n=\mathcal{I}_n^{(n/2)}(\tau_0,f_0)/n^3$.
By~\eqref{eq:log_lik_ratio}, with $\Lambda=f_{\tau_s}-f_{\tau_0}$ and $\Gamma=f_{\tau_0}-(\tilde f_{0,n})_{\tau_0}=\psi_{n/2}(\tau_0,f)-\theta_0$,
\[
\ell_n(\tau_s,f)-\ell_n(\tau_0,f)=\int_0^n\Lambda\,dW-\frac12|\Lambda|_n^2-\langle\Lambda,\Gamma\rangle_n.
\]
Hence the left side of~\eqref{eq:proof_BvM_expansion} is at most
\begin{equation}\label{eq:proof_BvM_terms}
\Big|\int_0^n(\Lambda-\lambda)\,dW\Big|+\frac12|\Lambda-\lambda|_n\big(|\Lambda-\lambda|_n+2|\lambda|_n\big)+|\langle\lambda,\Gamma\rangle_n|+|\Lambda-\lambda|_n|\Gamma|_n+\frac{s^2}2|\mathcal{J}_n-\mathcal{I}_0|.
\end{equation}

Split $\tilde f_{0,n}=f_n^*+q$, with $f_n^*=(\tilde f_{0,n})_{[\bar k_n]}$, and let $u=f-f_n^*$, a trigonometric polynomial of degree at most $\bar k_n$, since $K\le\bar k_n$.
As in the proof of Theorem~\ref{theo:checking_conditions}, by~\eqref{eq:truncation} and~\eqref{eq:cond}, and since $\bar k_n\ge k_n$, the shape $q$, which has the amplitudes of $f_0-f_{0,[\bar k_n]}$, satisfies $\|q\|_2\le M\{2\pi(\bar k_n+1)\}^{-\alpha}\le\epsilon_n/6$, and similarly, since $\alpha\ge1$,
\[
\|q'\|_2\le M\{2\pi(\bar k_n+1)\}^{1-\alpha}\le2\pi(\bar k_n+1)\epsilon_n/6\le2\pi\bar k_n\epsilon_n.
\]
Also $\|u\|_2\le\|f-\tilde f_{0,n}\|_2+\|q\|_2\le(L_1+1)\epsilon_n$, so $\|u'\|_2\le2\pi\bar k_n\|u\|_2\le2\pi(L_1+1)\bar k_n\epsilon_n$, and $\|(f_n^*)''\|_2\le2\pi\bar k_n\|f_0'\|_2\le2\pi\bar k_nM$, by Section~\ref{sec:preliminaries} and since $\alpha\ge1$.
By~\eqref{eq:cond_BvM} and since $n\epsilon_n^2\to\infty$, the sequence $\beta_n=\bar k_n\epsilon_n$ satisfies $\beta_n\to0$ and $n^{1/2}\epsilon_n\beta_n\to0$.
By~\eqref{eq:cond}, $\bar k_n\log n\lesssim n\epsilon_n^2$, so that also $\bar k_n^{1/2}(\log n)^{1/2}\beta_n\lesssim n^{1/2}\epsilon_n\beta_n\to0$ and $\bar k_n^{3/2}(\log n/n)^{1/2}\lesssim\beta_n\to0$.

Since
\[
\Lambda-\lambda=\Big\{(f_n^*)_{\tau_s}-(f_n^*)_{\tau_0}-\delta\Big(\cdot-\frac n2\Big)\big((f_n^*)'\big)_{\tau_0}\Big\}+\big\{u_{\tau_s}-u_{\tau_0}\big\}-\delta\Big(\cdot-\frac n2\Big)(q')_{\tau_0},
\]
the bounds~\eqref{eq:proof_BvM_second_order} and~\eqref{eq:proof_BvM_lipschitz} and the bound on $|(\cdot-n/2)g_\tau|_n$ give
\[
|\Lambda-\lambda|_n\le\frac{s^2\|(f_n^*)''\|_2}{\sqrt{32}\,n^{1/2}}+\frac{|s|\big(\|u'\|_2+\|q'\|_2\big)}{\sqrt2}\le C\big(s^2n^{-1/2}\bar k_n+|s|\beta_n\big)\le C'|s|\beta_n,
\]
for constants $C,C'$ that depend only on $L$, $L_1$ and $M$, where we used $|s|n^{-1/2}\le\rho_nn^{-1/2}=L\epsilon_n$ in the last step.
Also $|\lambda|_n=|s|\mathcal{J}_n^{1/2}$, $|\Gamma|_n\le(2n)^{1/2}\|f-\tilde f_{0,n}\|_2\le(2n)^{1/2}L_1\epsilon_n$, and $\mathcal{J}_n\to\mathcal{I}_0$, by~\eqref{eq:information_2}.
We bound the terms in~\eqref{eq:proof_BvM_terms} other than the first.
The function $\Gamma$ is $\tilde g(\tau_0\,\cdot)$ for a shape $\tilde g$ with $\|\tilde g\|_2=\|f-\tilde f_{0,n}\|_2$, so~\eqref{eq:information_3} applies to the third term, and~\eqref{eq:information_2} applies to the fifth.
With the bounds above, the second to fifth terms are at most, respectively,
\begin{gather*}
C's^2\beta_n\big(C'\beta_n+2\mathcal{J}_n^{1/2}\big),
\qquad
\frac{|\delta|n\|f_0'\|_2L_1\epsilon_n}{2\tau_0}=\frac{|s|\,\|f_0'\|_2L_1\epsilon_n}{2n^{1/2}\tau_0},\\
C\big(s^2n^{-1/2}\bar k_n+|s|\beta_n\big)(2n)^{1/2}L_1\epsilon_n\le\sqrt2\,CL_1\big(s^2\beta_n+|s|n^{1/2}\epsilon_n\beta_n\big),
\qquad
\frac{9s^2\mathcal{I}_0}{2n\tau_0}.
\end{gather*}
Each of these bounds is at most $1+s^2$ times a deterministic sequence that tends to zero.

For the first term, write $f_n^*=\sum_{|i|\le\bar k_n}x_i^*\varphi_i$ and $u=\sum_{|i|\le\bar k_n}y_i\varphi_i$, and note that $(\varphi_0)_\tau$ does not depend on $\tau$.
By Lemma~\ref{lemma:maximal}, almost surely, for all $|s|\le\rho_n$,
\[
\int_0^n(\Lambda-\lambda)\,dW=\frac{s^2}{n^{1/2}}\sum_{0<|i|\le\bar k_n}x_i^*v_i^2\,\mathbb{H}_i(s)+s\!\!\sum_{0<|i|\le\bar k_n}\!\!y_iv_i\,\mathbb{G}_i(s)-sV_n,
\]
where $V_n=n^{-3/2}\int_0^n(t-n/2)(q')_{\tau_0}(t)\,dW_t$ does not depend on $s$ or $f$, and is normally distributed with mean zero and variance at most $\|q'\|_2^2/2$.
Let $\Upsilon_n=\max_{0<|i|\le\bar k_n}\sup_{|s|\le\rho_n}\{|\mathbb{G}_i(s)|+|\mathbb{H}_i(s)|\}$.
Since
\[
\sum_{0<|i|\le\bar k_n}(x_i^*)^2v_i^4=\|(f_n^*)''\|_2^2
\quad\text{and}\quad
\sum_{0<|i|\le\bar k_n}y_i^2v_i^2=\|u'\|_2^2,
\]
the Cauchy-Schwarz inequality gives
\[
\Big|\int_0^n(\Lambda-\lambda)\,dW\Big|\le\frac{s^2}{n^{1/2}}\|(f_n^*)''\|_2(2\bar k_n)^{1/2}\Upsilon_n+|s|\,\|u'\|_2(2\bar k_n)^{1/2}\Upsilon_n+|s||V_n|.
\]
Since $\rho_n$ and $\bar k_n$ are polynomially bounded in $n$, Lemma~\ref{lemma:maximal} gives $\Upsilon_n=O_P\big((\log n)^{1/2}\big)$.
By the bounds above, the right side is at most $(1+s^2)$ times
\[
C''\Big\{\bar k_n^{3/2}n^{-1/2}+\bar k_n^{1/2}\beta_n\Big\}\Upsilon_n+|V_n|,
\]
for a constant $C''$, which tends to zero in probability, since $\bar k_n^{3/2}(\log n/n)^{1/2}\to0$, $\bar k_n^{1/2}(\log n)^{1/2}\beta_n\to0$ and $\|q'\|_2\le2\pi\beta_n\to0$.
Adding the bounds on the five terms of~\eqref{eq:proof_BvM_terms} gives~\eqref{eq:proof_BvM_expansion}, with $\eta_n$ the sum of the coefficients of $1+s^2$.
\qed

\paragraph{Proof of Corollary~\ref{cor:BvM_example}:}
By Proposition~\ref{prop:prior_example}, the conditions of Theorems~\ref{theo:checking_conditions} and~\ref{theo:checking_conditions_2} hold, and the density of $\Pi_\tau$ is continuous and positive by the assumptions of Section~\ref{sec:examples}.
The prior on the shape is stationary, by Section~\ref{sec:examples}.
With $\tilde k_n$ and $D$ as in Proposition~\ref{prop:prior_example} and its proof, with $D=1$ in case b), we have $\bar k_n\le2D\tilde k_n$ for large $n$, $\epsilon_n^2=C_0^2\tilde k_n^{-2\alpha}$ and $n=\tilde k_n^{2\alpha+1}\log n$, so that
\[
n^{1/2}\bar k_n\epsilon_n^2\le2DC_0^2\tilde k_n^{3/2-\alpha}(\log n)^{1/2},
\]
which tends to zero, since $\alpha>3/2$ and $\tilde k_n$ is a positive power of $n/\log n$.
Theorem~\ref{theo:BvM} gives the result.
\qed

\section*{Acknowledgements}
I disclose that an AI assistant (Claude, by Anthropic) was used in the preparation of the manuscript: 
to search and check the literature,
to shorten and validate the author's arguments and proofs,
and to revise and give feedback on earlier versions of the manuscript. 
Special care was taken to ensure that all results in the paper have proper attribution.
I have read the entire manuscript, and take full responsibility for its content.

\bibliographystyle{plainnat}
\bibliography{refs}

\end{document}